\documentclass[a4paper,11pt]{article}
\usepackage[latin1]{inputenc}
\usepackage[francais,english]{babel}    
\usepackage{amssymb}
\usepackage{amsmath}
\usepackage{amsthm}
\usepackage{multicol}
\usepackage{graphicx}
\usepackage{color}
\usepackage[T1]{fontenc}
\usepackage{yfonts}
\usepackage{placeins}
\usepackage{array}
\usepackage{dsfont}
\usepackage{stmaryrd}
\usepackage{verbatim}
\usepackage{caption}
\usepackage{transparent}
\usepackage{bbold}
\usepackage{hyperref}
\usepackage{enumitem}
\usepackage{tikz}
\usetikzlibrary{patterns}
\usepackage{tikz-3dplot}

\newtheorem{thm}{Theorem}[section]
\newcommand{\bthm}{\begin{thm}}
\newcommand{\ethm}{\end{thm}}

\newtheorem{thmi}{Theorem}
\newcommand{\bthmi}{\begin{thmi}}
\newcommand{\ethmi}{\end{thmi}}

\newtheorem{cori}[thmi]{Corollary}
\newcommand{\bcori}{\begin{cori}}
\newcommand{\ecori}{\end{cori}}

\newtheorem{mthm}{Theorem}
\newcommand{\bmthm}{\begin{mthm}}
\newcommand{\emthm}{\end{mthm}}

\newtheorem{mcor}[mthm]{Corollary}
\newcommand{\bmcor}{\begin{mcor}}
\newcommand{\emcor}{\end{mcor}}
\newtheorem{mconj}[mthm]{Conjecture}
\newcommand{\bmconj}{\begin{mconj}}
\newcommand{\emconj}{\end{mconj}}
\newtheorem{mpro}[mthm]{Proposition}
\newcommand{\bmpro}{\begin{mpro}}
\newcommand{\empro}{\end{mpro}}

\newtheorem*{conj}{Conjecture}
\newcommand{\bconj}{\begin{conj}}
\newcommand{\econj}{\end{conj}}

\newtheorem*{question}{Question}
\newcommand{\bq}{\begin{question}}
\newcommand{\eq}{\end{question}}

\newtheorem*{thn}{Theorem}
\newcommand{\bthn}{\begin{thn}}
\newcommand{\ethn}{\end{thn}}

\newtheorem{exo}{Exercise}
\newcommand{\bex}{\begin{exo}}
\newcommand{\eex}{\end{exo}}

\newtheorem{sol}{Solution}
\newcommand{\bsol}{\begin{sol}}
\newcommand{\esol}{\end{sol}}

\newtheorem{pro}[thm]{Proposition}
\newcommand{\bpro}{\begin{pro}}
\newcommand{\epro}{\end{pro}}

\newtheorem{cor}[thm]{Corollary}
\newcommand{\bcor}{\begin{cor}}
\newcommand{\ecor}{\end{cor}}

\newtheorem{lem}[thm]{Lemma}
\newcommand{\blem}{\begin{lem}}
\newcommand{\elem}{\end{lem}}

\theoremstyle{definition}

\newtheorem{defi}[thm]{Definition}
\newcommand{\bdf}{\begin{defi}}
\newcommand{\edf}{\end{defi}}

\newtheorem*{defis}{Definition}
\newcommand{\bdfs}{\begin{defis}}
\newcommand{\edfs}{\end{defis}}

\newtheorem*{rmk}{Remark}
\newcommand{\brk}{\begin{rmk} \upshape}
\newcommand{\erk}{\end{rmk}}

\newtheorem*{rmks}{Remarks}
\newcommand{\brks}{\begin{rmks} \upshape}
\newcommand{\erks}{\end{rmks}}

\newtheorem*{exe}{Example}
\newcommand{\bexe}{\begin{exe} \upshape}
\newcommand{\eexe}{\end{exe}}

\newtheorem*{exes}{Examples}
\newcommand{\bexes}{\begin{exes} \upshape}
\newcommand{\eexes}{\end{exes}}

\newtheorem*{pre}{Proof}
\newcommand{\bp}{\begin{pre} \upshape}
\newcommand{\ep}{\hfill \qed \end{pre}}
\newcommand{\epp}{\end{pre}}

\newcommand{\beq}{\begin{eqnarray*}}
\newcommand{\eeq}{\end{eqnarray*}}

\newcommand{\beqn}{\begin{equation}}
\newcommand{\eeqn}{\end{equation}}

\newcommand{\ben}{\begin{enumerate}}
\newcommand{\een}{\end{enumerate}}
\newcommand{\bit}{\begin{itemize} \renewcommand{\labelitemi}{$\bullet$} \renewcommand{\labelitemii}{$\star$}}
\newcommand{\eit}{\end{itemize}}

\newcommand{\bfg}{
\begin{figure}[H]
\begin{center}}
\newcommand{\efg}{
\end{center}
\end{figure}
\FloatBarrier}

\newcolumntype{M}[1]{>{\raggedright}m{#1}}

\newcommand{\R}{\mathbb{R}}

\newcommand{\N}{\mathbb{N}}
\newcommand{\Z}{\mathbb{Z}}
\newcommand{\C}{\mathbb{C}}

\newcommand{\D}{\mathbb{D}}
\renewcommand{\SS}{\mathbb{S}}
\renewcommand{\H}{\mathbb{H}}

\newcommand{\bs}{\symbol{92}}

\newcommand{\ov}{\overline}

\renewcommand{\tilde}{\widetilde}

\renewcommand{\th}{\operatorname{th}}

\newcommand{\Isom}{\operatorname{Isom}}

\renewcommand{\max}{\operatorname{max}}

\newcommand{\MCG}{\operatorname{MCG}}

\newcommand{\eps}{\varepsilon}
\newcommand{\st}{\, | \,}

\newcommand{\ra}{\rightarrow}

\newcommand{\ral}[1]{\underset{#1}{\longrightarrow}}

\newcommand{\liml}{\lim\limits}

\newcommand{\f}{\frac}

\renewcommand{\geq}{\geqslant}
\renewcommand{\leq}{\leqslant}

\newcommand{\PGL}{\operatorname{PGL}}

\newcommand{\<}{\langle}
\renewcommand{\>}{\rangle}
\newcommand{\pif}{{+\infty}}

\newcommand{\mk}{\medskip}

\makeatletter
\def\Ddots{\mathinner{\mkern1mu\raise\p@
\vbox{\kern7\p@\hbox{.}}\mkern2mu
\raise4\p@\hbox{.}\mkern2mu\raise7\p@\hbox{.}\mkern1mu}}
\makeatother

\makeatletter
\def\maketitles{%
  \null
  \thispagestyle{empty}%
  \vfill
  \begin{center}\leavevmode
    \normalfont
    {\LARGE \@title\par}%
    \vskip 1.2cm
    {\large \@author\par}%
    \vskip 1.2cm
    {\large \@subtitle\par}%
    \vskip 0.8cm
    {\large \@date\par}%
  \end{center}%
  \vfill
  \null
  \cleardoublepage
  }

\def\date#1{\def\@date{#1}}
\def\author#1{\def\@author{#1}}
\def\title#1{\def\@title{#1}}
\def\subtitle#1{\def\@subtitle{#1}}
\makeatother

\usetikzlibrary{arrows}

\title{Lattices, injective metrics and the $K(\pi,1)$ conjecture}
\author{Thomas Haettel\thanks{thomas.haettel@umontpellier.fr, IMAG, Univ Montpellier, CNRS, France, and IRL 3457, CRM-CNRS, Universit\'{e} de Montr\'{e}al, Canada.}}
\date{\today}

\begin{document}

\selectlanguage{english}

\maketitle

\begin{center}
\begin{minipage}{0.8\textwidth}
\textsc{Abstract.} Starting with a lattice with an action of $\Z$ or $\R$, we build a Helly graph or an injective metric space. We deduce that the $\ell^\infty$ orthoscheme complex of any bounded graded lattice is injective. We also prove a Cartan-Hadamard result for locally injective metric spaces. We apply this to show that any Garside group or any FC type Artin group acts on an injective metric space and on a Helly graph. We also deduce that the natural piecewise $\ell^\infty$ metric on any Euclidean building of type $\tilde{A_n}$ extended, $\tilde{B_n}$, $\tilde{C_n}$ or $\tilde{D_n}$ is injective, and its thickening is a Helly graph.

Concerning Artin groups of Euclidean types $\tilde{A_n}$ and $\tilde{C_n}$, we show that the natural piecewise $\ell^\infty$ metric on the Deligne complex is injective, the thickening is a Helly graph, and it admits a convex bicombing. This gives a metric proof of the $K(\pi,1)$ conjecture, as well as several other consequences usually known when the Deligne complex has a CAT(0) metric.
\end{minipage}
\end{center}

\let\thefootnote\relax

\tableofcontents

\section*{Introduction}

Injective metric spaces are geodesic metric spaces where every family of pairwise intersecting balls has a non-empty global intersection. Their discrete counterpart are called Helly graphs. Their use in geometric group theory is recent and growing, see notably~\cite{dress}, \cite{lang}, \cite{huang_osajda_helly}, \cite{helly_groups}, \cite{haettel_hoda_petyt}, \cite{osajda_valiunas}, \cite{haettel_injective_buildings}. Roughly speaking, CAT(0) spaces are typically locally Euclidean spaces, whereas injective metric spaces are typically locally $\ell^\infty$ metric spaces. Injective metric spaces display many nonpositive curvature features observed in CAT(0) spaces. We believe that Helly graphs and injective metric spaces may also reveal more powerful than CAT(0) spaces for some purposes.

\mk

Firstly, it appears that many nonpositively curved groups (notably most CAT(0) groups) have a nice action on an injective metric space, so that the theory encompasses a vast class of groups: CAT(0) cubical groups, hyperbolic groups, relatively hyperbolic groups, uniform lattices in semisimple Lie groups over local fields, braid groups and more generally Garside groups, Artin groups of type FC, mapping class groups and more generally hierarchically hyperbolic groups. For instance, Helly graphs admit a simple combinatorial local characterization (see Theorem~\ref{thm:clique_helly_global_helly}), which makes them potentially as powerful as CAT(0) cube complexes, whereas piecewise Euclidean CAT(0) complexes desperately lack a combinatorial characterization.

\mk

Secondly, we can still deduce many of the results that hold true for CAT(0) spaces and CAT(0) groups in the setting of injective spaces and injective groups (groups acting geometrically on injective spaces). For instance, injective metric spaces have a conical geodesic bicombing and are thus contractible, and every isometric bounded group action has a fixed point. Moreover, every injective group is semi-hyperbolic (see~\cite{alonsobridson:semihyperbolic}) and satisfies the Farrell-Jones conjecture (see~\cite{kasprowski_rueping}). Sometimes we can even deduce stronger results than in the CAT(0) setting: as a sample result, note that any group acting properly and cocompactly on a Helly graph is biautomatic (see~\cite{helly_groups}), whereas not all CAT(0) groups are biautomatic (see~\cite{leary_minasyan}). Also note that infinite, finitely generated torsion groups do not act properly on uniformly locally finite Helly graphs (see~\cite{haettel_osajda_locally_elliptic}), whereas the analogous statement is open for CAT(0) complexes.

\mk

In this article, we will pursue this philosophy. On one hand, we develop results useful to prove that some metric spaces are injective. On the other hand, we apply these results to Euclidean buildings and Deligne complexes of Euclidean Artin groups. We believe, however, that the scope of our results will not be limited to Artin groups of Euclidean type, and could concern much larger classes of groups.

\mk

We now present a very simple criterion (already appearing in a restricted form in~\cite{haettel_injective_buildings}) showing how to build a Helly graph or an injective metric space, starting with a lattice and an action of $\Z$ or $\R$.

\bmthm \label{mthm:criterion_helly_injective} (Theorem~\ref{thm:lattice_Helly_injective})
Assume that $L$ is a lattice such that each upperly bounded subset of $L$ has a join. Assume that there is an order-preserving increasing continuous action $(f_t)_{t \in H}$ of $H=\Z$ or $\R$ on $L$, such that
$$\forall x,y \in L, \exists t \in H_+, f_{-t}(x) \leq y \leq f_t(x).$$
Let us define the following metric $d$ on $L$:
$$\forall x,y \in L, d(x,y) =\inf\{t \in H_+ \st f_{-t}(x) \leq y \leq f_t(x)\}.$$
\bit
\item If $H=\Z$, the metric space $(L,d)$ is the vertex set of a Helly graph.
\item If $H=\R$, the metric space $(L,d)$ is injective.
\eit
\emthm

With Jingyin Huang, we used the same situation of a lattice with an appropriate action of $\Z$ to define another natural graph which is weakly modular, and this applies to numerous examples (see~\cite{haettel_huang_weakly_modular}).

\mk

We applied this criterion in~\cite{haettel_injective_buildings} to prove that the thickening of a Bruhat-Tits Euclidean building of extended type $\tilde{A_n}$ is a Helly graph. It turns out that the Bruhat-Tits restriction is unnecessary.

\bmthm (Theorems~\ref{thm:AN_extended_injective} and \ref{thm:building_typesBCD_injective})
The natural $\ell^\infty$ metric on any Euclidean building of type $\tilde{A_n}$ extended, $\tilde{B_n}$, $\tilde{C_n}$ or $\tilde{D_n}$ is injective. Furthermore, the thickening of its vertex set is a Helly graph.
\emthm

Since Euclidean buildings are also endowed with CAT(0) metrics, we may wonder about the importance of this result. For once, it gives a lot of new examples to the theory of injective metric spaces and Helly graphs. Another consequence is another approach to \'{S}wi\c{a}tkowski's result (see~\cite{swiatkowski_biautomatic}) that cocompact lattices in Euclidean buildings are biautomatic, in types $\tilde{A_n}$, $\tilde{B_n}$, $\tilde{C_n}$ or $\tilde{D_n}$. The proof for types $\tilde{B_n}$, $\tilde{C_n}$ and $\tilde{D_n}$ relies on a generalization for graded semilattices, see Section~\ref{sec:semilattice}.

\mk

Another immediate consequence of Theorem~\ref{mthm:criterion_helly_injective} is the following.

\bmthm (Corollary~\ref{cor:thickening_garside_helly})
The thickening of the Cayley graph of any Garside group with respect to its simple elements is a Helly graph.
\emthm

Note that, in the case of a finite type Garside group, this is due to Huang and Osajda (see~\cite{huang_osajda_helly}). However, our proof is different, and does not rely on the deep local-to-global result for Helly graphs (see Theorem~\ref{thm:clique_helly_global_helly}). Besides, it also works in the case of a Garside group with infinite set of simple elements.

\mk

One application is the study of the orthoscheme complex of a bounded, graded lattice $L$. Note that simplices of the geometric realization $|L|$ of $L$ correspond to chains in $L$, so that one can endow each simplex with metric of the standard $\ell^\infty$ orthosimplex associated to this order on vertices (see Section~\ref{sec:review} for details). We endow the geometric realization $|L|$ of $L$ with a lattice structure, and use Theorem~\ref{mthm:criterion_helly_injective} to deduce the following.

\bmthm \label{mthm:orthoscheme_complex_lattice_injective} (Theorem~\ref{thm:lattice_orthoscheme_injective}) 
Let $L$ denote a bounded, graded lattice. Then the orthoscheme complex $|L|$ of $L$, with the piecewise $\ell^\infty$ metric, is injective.
\emthm

When $L$ is a bounded, graded lattice, one may endow its orthoscheme complex with the piecewise Euclidean metric. Deciding whether it is a CAT(0) metric space is a very difficult question. It turns out that for orthoscheme complexes of posets, the CAT(0) property is more restrictive than the injective property, see Theorem~\ref{thm:orthoscheme_cat0_vs_injective}.

\mk

One famous conjecturally CAT(0) example is the dual braid complex, defined by Brady and McCammond in~\cite{brady_mccammond}. The $n$-strand braid group $B_n$ has a standard Garside structure, associated to the "half-turn" as Garside element. The braid group also enjoys a dual presentation introduced by Birman, Ko and Lee (see~\cite{birman_ko_lee}), which corresponds to a dual Garside structure (see~\cite{dehornoy_paris_gaussian}). It is associated to the rotation of a $n\th$ of a turn as Garside element. Associated to this dual Garside structure, one may consider the geometric realization and endow it with the piecewise Euclidean orthoscheme metric (see Section~\ref{sec:review}). Brady and McCammond conjecture that this so-called dual braid complex is CAT(0) for all braid groups, but it has only been proved for $n \leq 7$ (see~\cite{brady_mccammond}, \cite{b6}, \cite{jeong}). However, an immediate consequence of Theorem~\ref{mthm:orthoscheme_complex_lattice_injective} is the following.

\bmthm (Corollary~\ref{cor:garside_complex_injective}) \label{mthm:garside_complex_injective}
The Garside complex of any Garside group, endowed with the piecewise $\ell^\infty$ orthoscheme metric, is injective.
\emthm

We also obtain another proof of a result by Huang and Osajda stating that FC type Artin groups are Helly (\cite[Theorem~5.8]{huang_osajda_helly}). We also provide explicit Helly and injective models (see Theorem~\ref{thm:FC_Artin_Helly} for the precise statement).

\bmthm[Huang-Osajda] (Theorem~\ref{thm:FC_Artin_Helly})
Let $A$ denote an Artin group of type FC. Then a natural simplicial complex $X$ with vertex set $A$, with the $\ell^\infty$ metric, is injective. Moreover, a thickening of $X$ is a Helly graph. 
\emthm

\mk

As a particular case of Theorem~\ref{mthm:garside_complex_injective}, the dual braid complex, endowed with the piecewise $\ell^\infty$ orthoscheme metric, is injective for all braid groups. This also holds, more generally, for every spherical type Artin group with some Garside structure. Note that Theorem~\ref{mthm:orthoscheme_complex_lattice_injective} proves local injectivity, and one also needs a Cartan-Hadamard result for injective metric spaces in order to conclude. We therefore rely on the local-to-global result for Helly graphs to prove the following generalization of~\cite{miesch} in the non-proper setting (see Section~\ref{sec:review} for the definition of semi-uniformly locally injective). This result is clearly of independent interest.

\bmthm[Cartan-Hadamard for injective metric spaces] (Theorem~\ref{thm:cartan_hadamard_injective})
Let $X$ denote a complete, simply connected, semi-uniformly locally injective metric space. Then $X$ is injective.
\emthm

Another very promising family of examples are Deligne complexes of Artin groups (see Section~\ref{sec:review} for definitions). Note that buildings and Deligne complexes of Artin groups are closely related: in fact in his original article (see~\cite{deligne}), Deligne called "buildings for generalized braid groups" the complexes later called Deligne complexes. However, to the best of our knowledge, the close relationship between Euclidean buildings and Deligne complexes of Euclidean type Artin groups has not yet been exploited in the literature. Notably, the automorphism groups of Euclidean buildings do not possess a Garside structure, and the Deligne complexes do not have an apartment system as rich as in the building case. However, one common feature is that they locally look like a lattice, which is the key combinatorial property we are using in this article.

\mk

Associated to every Coxeter graph $\Gamma$ with vertex set $S$, we may define the Coxeter group $W(\Gamma)$, the Artin group $A(\Gamma)$ and the hyperplane complement $M(\Gamma)$ (see Section~\ref{sec:review}). The Coxeter group $W(\Gamma)$ acts naturally on $M(\Gamma)$, and $A(\Gamma)$ is the fundamental group of the quotient $W(\Gamma) \bs M(\Gamma)$. One very natural question is to decide whether it is a classifying space. This is the statement of the following conjecture.

\bmconj[The $K(\pi,1)$ conjecture]
The hyperplane complement $M(\Gamma)$ is aspherical.
\emconj

The $K(\pi,1)$ conjecture has been proved for spherical type Artin groups by Deligne in~\cite{deligne}, and for Euclidean type Artin groups by Paolini and Salvetti in~\cite{paolini_salvetti} very recently, even for the type $\tilde{D_n}$. Another approach, more closely related to the metric approach of this article, has been used by Charney and Davis in~\cite{charney_davis_kpi1} to prove the $K(\pi,1)$ conjecture for Artin groups of type FC or of $2$-dimensional type. Their proof relies on the use of a simplicial complex, called the Deligne complex $\Delta(\Gamma)$, which is the geometric realization of the poset of cosets of spherical type parabolic subgroups (see Section~\ref{sec:review}). This complex (in this form) has been defined by Charney and Davis in~\cite{charney_davis_kpi1}, where they proved that $\Delta(\Gamma)$ is homotopy equivalent to the universal cover of $W(\Gamma) \bs M(\Gamma)$. In other words, the $K(\pi,1)$ conjecture is equivalent to the contractibility of the Deligne complex $\Delta(\Gamma)$.

\mk

Charney and Davis's method for proving that the Deligne complex is contractible is to endow it with a CAT(0) metric. This works for Artin groups of type FC or of $2$-dimensional type. However, in the general case, the key question is to decide whether the Deligne complex for the braid group is CAT(0) with Moussong's metric. It is only known for the braid group up to $4$ strands only (see~\cite{charney_deligne_b4}). However, we will see that the natural piecewise $\ell^\infty$ metric is injective for all braid groups, up to taking the product with $\R$.

\bmthm (Theorem~\ref{thm:AN_extended_injective})
Let $\Delta$ denote the Deligne complex of the Artin group of Euclidean type $\tilde{A_n}$. Then the natural piecewise $\ell^\infty$ length metric on $\Delta \times \R$ is injective. Moreover, the thickening of $\Delta \times \R$ is a Helly graph.
\emthm

The proof consists in applying Theorem~\ref{mthm:orthoscheme_complex_lattice_injective} to prove that the Deligne complex is locally injective. One key combinatorial property is that the Deligne complex in type $A_n$ is essentially a lattice, see Section~\ref{sec:cut_curves}. The proof of the lattice property, through the cut-curve lattice, is due to Crisp and McCammond, copied here with their permission.

\mk

Note that one may wonder whether it is necessary to consider the direct product with $\R$. In fact, Hoda proved (see~\cite{hoda:crystallographic}) that the Euclidean Coxeter group of type $\tilde{A_n}$ is not Helly for $n \geq 2$, even though its direct product with $\Z$ is. We made a similar distinction for automorphism groups of Euclidean buildings of type $\tilde{A_n}$ in~\cite{haettel_injective_buildings}. We therefore strongly believe that there is no injective metric on the Deligne complex of type $\tilde{A_n}$ itself which is invariant under the Artin group. However, we will see in Theorem~\ref{mthm:deligne_bicombing} below that there is a convex bicombing on the Deligne complex itself.

\mk

In order to deal with the Euclidean type $\tilde{C_n}$, we first prove a generalization of Theorem~\ref{mthm:orthoscheme_complex_lattice_injective} for graded semilattices, see Section~\ref{sec:semilattice}.

\bmthm (Theorem~\ref{thm:building_typesBCD_injective})
Let $\Delta$ denote the Deligne complex of the Artin group of Euclidean type $\tilde{C_n}$. Then the natural piecewise $\ell^\infty$ length metric on $\Delta$ is injective. Moreover, the thickening of $\Delta$ is a Helly graph.
\emthm

An immediate consequence is another proof of the $K(\pi,1)$ conjecture in Euclidean types $\tilde{A_n}$ and $\tilde{C_n}$, originally due to Okonek (see~\cite{okonek}). The novelty is that it is the first metric proof. Moreover, in Charney and Davis's approach to the $K(\pi,1)$ conjecture by showing that Mousson's metric is CAT(0), the main difficulty is to prove that it is locally CAT(0) for the braid groups. It is precisely this statement that we are able to prove in the injective setting.

\bmcor[Okonek]
The Deligne complex $\Delta$ of Euclidean type $\tilde{A_n}$ or $\tilde{C_n}$ is contractible. In particular, the $K(\pi,1)$ conjecture holds in these cases.
\emcor

Moreover, several results have been proved, relying on the assumption that one may endow the Deligne complex with a piecewise Euclidean CAT(0) metric. Crisp studied the fixed-point subgroup under a symmetry group of the Artin system (see~\cite{crisp_symmetrical}). Godelle studied the centralizer and normalizer of standard parabolic subgroups (see~\cite{godelle_cat0}). Morris-Wright studied the intersections of parabolic subgroups (see~\cite{morris_wright_parabolic_FC}). It turns out almost all the arguments merely used the existence of an equivariant geodesic bicombing on the Deligne complex (see the end of Section~\ref{sec:affine_lattice} for a definition of a bicombing). Moreover, work of Descombes and Lang notably (see~\cite{descombes_lang_flats} and \cite{descombes_lang_hyperbolicity}) justify the importance of convex geodesic bicombings for themselves. We therefore state the following conjecture, which may be seen as a metric strategy for the proof of the $K(\pi,1)$ conjecture.

\bmconj \label{mconj:deligne_bicombing}
The Deligne complex of any Artin group $A$ has an $A$-invariant metric that admits a convex, consistent, reversible geodesic bicombing.
\emconj

We are able to prove this conjecture in spherical types $A_n$, $B_n$ and Euclidean types $\tilde{A_n}$, $\tilde{C_n}$, and we believe that our result represents a major step towards the general case.

\bmthm (Theorem~\ref{thm:deligne_bicombing}) \label{mthm:deligne_bicombing}
Let $\Delta$ denote the Deligne complex of the Artin group of spherical type $A_n$ or $B_n$ or Euclidean type $\tilde{A_n}$ or $\tilde{C_n}$. There exists a metric on $\Delta$, invariant under the Artin group, which admits a convex, consistent, reversible geodesic bicombing.
\emthm

If an Artin group satisfies Conjecture~\ref{mconj:deligne_bicombing}, then the $K(\pi,1)$ conjecture follows. Moreover, we may also list consequences of work of Crisp (see~\cite{crisp_symmetrical}), Godelle (see~\cite{godelle_cat0}) and Morris-Wright (see~\cite{morris_wright_parabolic_FC} and \cite{cumplido_parabolic_spherical}) that rely on the assumption that the Deligne complex has a CAT(0) metric. However, most of the arguments only use the geodesic bicombing. The following results were only essentially known for Artin groups of type FC or of $2$-dimensional type. Note that concerning Artin groups of type FC, only intersections of spherical type parabolic subgroups are known to be parabolic (see~\cite{morris_wright_parabolic_FC}). see also~\cite{moller_paris_varghese}.

\bmcor (Corollaries~\ref{cor:intersection_of_parabolics}, \ref{cor:centralizers_normalizers} and \ref{cor:symmetry_artin_system})
Let $A$ denote the Artin group of Euclidean type $\tilde{A_n}$ or $\tilde{C_n}$.
\bit
\item The intersection of any parabolic subgroups of $A$ is a parabolic subgroup.
\item $A$ satisfies Properties $(\star)$, $(\star \star)$ and $(\star \star \star)$ from~\cite{godelle_cat0}, notably: for any subset $X \subset S$, we have
$$Com_A(A_X) = N_A(A_X) =  A_X \cdot QZ_A(X),$$
where the quasi-centralizer of $X$ is $QZ_A(X) = \{g \in A \st g \cdot X = X\}$.
\item For any group $G$ of symmetries of the Artin system, the fixed point subgroup $A^G$ is isomorphic to an Artin group.
\eit
\emcor

In view of Conjecture~\ref{mconj:deligne_bicombing}, looking for other consequences of the existence of a convex bicombing on the Deligne complex may reveal fruitful.

\subsection*{Structure of the article} 

In Section~\ref{sec:review}, we review basic definitions of posets, lattices, Artin groups, Deligne complexes, injective metric spaces and Helly graphs. We also prove the Cartan-Hadamard Theorem for injective metric spaces. In Section~\ref{sec:thickening_lattice}, we prove the central simple criterion showing how to produce an injective metric space or a Helly graph starting from a lattice with an action of $\Z$ or $\R$. In Section~\ref{sec:affine_lattice}, we apply this criterion to prove that the orthoscheme complex of a bounded, graded lattice is injective. In Section~\ref{sec:type_affine_A}, we use this to prove that, for Euclidean buildings and the Deligne complex in Euclidean type $\tilde{A_n}$, the natural piecewise $\ell^\infty$ metric is injective. In Section~\ref{sec:semilattice}, we show how to adapt the criterion to a mere semilattice with some extra property. We then apply it in Section~\ref{sec:type_affine_C} to prove that, for Euclidean buildings and the Deligne complex in Euclidean type $\tilde{C_n}$, the natural piecewise $\ell^\infty$ metric is injective. In the last Section~\ref{sec:corollaries_using_bicombing}, we use the convex bicombing on the Deligne complexes to deduce many corollaries about parabolic subgroups of Artin groups.

\mk

\textbf{Acknowledgments:} We would like to thank very warmly Jingyin Huang, Jon McCammond, Damian Osajda and Luis Paris for very interesting discussions. We are indebted to Jon McCammond and John Crisp for their authorization to copy their work about the lattice of cut-curves. We also would like to thank the organizers of the ICMS Edimburgh 2021 meeting "Perspectives on Artin groups", where we had fruitful discussions. We also would like to thank J\'{e}r\'{e}mie Chalopin, Victor Chepoi, Mar\'{i}a Cumplido, Anthony Genevois, Hiroshi Hirai, Alexandre Martin and Bert Wiest for various discussions. We also thank an anonymous referee for many helpful comments.

We acknowledge support from French project ANR-16-CE40-0022-01 AGIRA.

\section{Review of posets, Artin groups and injective metric spaces} \label{sec:review}

\subsection{Posets and lattices}

We recall briefly basic definitions related to posets and lattices.

\bdf[Poset]\

Let $L$ denote a poset. A \emph{chain} of $L$ is a totally ordered non-empty subset of $L$. A \emph{maximal chain} is a chain that is maximal with respect to inclusion. A finite chain of $n+1$ elements $x_0<x_1< \dots < x_n$ is called of \emph{length $n$}. The poset $L$ is called \emph{bounded below} (resp. above) if it has a global minimum denoted $0$ (resp. global maximum denoted $1$). The poset $L$ is called \emph{bounded} if it is both bounded above and below. The poset $L$ has \emph{rank $n$} if it is bounded and all maximal chains have length $n$.
\edf

\bdf[Interval]\

Given two elements $x \leq y$ in a poset $L$, we denote the \emph{interval}
$$I(x,y) = \{z \in L \st x \leq z \leq y\}.$$
In case of possible confusion, we may also denote the interval $I_L(x,y)$ to emphasize that we are considering the interval in the poset $L$. The poset $L$ is called \emph{graded} if every interval of $L$ has a rank. Let $x$ denote an element in a graded poset $L$ that is bounded below, the \emph{rank} of $x$ is the rank of the interval $I(0,x)$.
\edf

\bdf[Lattice]\

Given elements $x,y$ in a poset $L$, if there exists a unique maximal lower bound to $\{x,y\}$, it is called the \emph{meet} of $x$ and $y$ and denoted $x \wedge y$. Similarly, if there exists a unique minimal upper bound to $\{x,y\}$, it is called the \emph{join} of $x$ and $y$ and denoted $x \vee y$. If any two elements of $L$ have a meet (resp. a join), the poset $L$ is called a \emph{meet-semilattice} (resp. join-semilattice). If any two elements of $L$ have a meet and a join, the poset $L$ is called a \emph{lattice}.
\edf

\bexes\
\bit
\item The Boolean lattice $L$ of rank $n$ is the poset of subsets of $E=\{1,\dots,n\}$, partially ordered by inclusion. The join of $A,B \in L$ is $A \cup B$ and their meet is $A \cap B$.
\item Consider a CAT(0) cube complex $X$, with a base vertex $v_0$. Order the set $V$ of vertices of $X$ by declaring that $v \leq w$ if some combinatorial geodesic from $v_0$ to $w$ passes through $v$. Then $V$ is a graded meet-semilattice, with minimum $v_0$. The meet of two vertices $v,w \in V$ is the median of $v_0,v,w$.
\item The partition lattice $L$ of $E=\{1,\dots,n\}$ is the poset of partitions of $E$, partially ordered by declaring that $A \leq B$ if every element of $A$ is contained in an element of $B$. This lattice has rank $n-1$, its minimum is the partition $\{\{1\},\{2\},\dots,\{n\}\}$ into singletons and its maximum is the partition $\{E\}$ into one element.
\eit
\eexes

We will now describe a very simple criterion due to Brady and McCammond (see~\cite{brady_mccammond}) to decide when a bounded graded poset is a lattice.

\bdf[Bowtie]\

In a poset $L$, a \emph{bowtie} consists of $4$ distinct elements $a,b,c,d$ such that $a,c$ are minimal upper bounds of $b,d$, and $b,d$ are maximal lower bounds of $a,c$. 
\edf

\bpro{\cite[Proposition~1.5]{brady_mccammond}} \label{pro:bowtie}
Let $L$ denote a bounded graded poset. Then $L$ is a lattice if and only if $L$ does not contain a bowtie.
\epro

\subsection{Coxeter groups, Artin groups and Deligne complexes} \label{sec:intro_artin}

We recall the definitions of Coxeter groups, Artin groups, and their associated Deligne complexes.

\mk

For every finite simple graph $\Gamma$ with vertex set $S$ and with edges labeled by some integer in $\{2,3,\dots\}$, one associates the Coxeter group $W(\Gamma)$ with the following presentation:
$$W(\Gamma) = \<S \st \forall \{s,t\} \in \Gamma^{(1)}, \forall s \in S, s^2=1, [s,t]_m=[t,s]_m \mbox{ if the edge $\{s,t\}$ is labeled $m$}\>,$$
where $[s,t]_m$ denotes the word $ststs\dots$ of length $m$.

\mk

The associated Artin group $A(\Gamma)$ is defined by a similar presentation:
$$A(\Gamma) = \<S \st \forall \{s,t\} \in \Gamma^{(1)}, [s,t]_m=[t,s]_m \mbox{ if the edge $\{s,t\}$ is labeled $m$}\>.$$
The groups $A(\Gamma)$ are also called Artin-Tits groups, since they have been defined by Tits in~\cite{tits_artin}.

Note that only the relations $s^2=1$ have been removed, so that there is a natural surjective morphism from $A(\Gamma)$ to $W(\Gamma)$. Also note that when $m=2$, then $s$ and $t$ commute, and when $m=3$, then $s$ and $t$ satisfy the classical braid relation $sts=tst$.

\mk

The knowledge of general Artin groups is extremely limited (see notably~\cite{mccammond_mysterious}, \cite{charney_problems}, \cite{godelle_paris}). In particular, we do not know whether the word problem is solvable, nor whether they are torsion-free.

\mk

Most results about Artin-Tits groups concern particular classes. The Artin group $A(\Gamma)$ is called:
\bit
\item of \emph{spherical type} if its associated Coxeter group $W(\Gamma)$ is finite, i.e. may be realized as a reflection group of a sphere.
\item of \emph{Euclidean type} if its associated Coxeter group $W(\Gamma)$ may be realized as a reflection group of a Euclidean space.
\item of \emph{FC type} if for any complete subset $T \subset S$ the parabolic subgroup $A_T=\<T\>$ is spherical.
\eit

\begin{table}
\begin{center}
\begin{tabular}{|c|c|c|c|}
\hline
& Spherical type & & Euclidean type \\

\hline

$A_n, n \geq 2$
&
\begin{tikzpicture}
\def \p {0.05}
\def \op {1}
\def \gris {black!10}

\draw[fill] (-3,0) circle (\p) node(s1) {};
\draw[fill] (-2,0) circle (\p) node(s2) {};
\draw[fill] (-1,0) circle (\p) node(s3) {};
\draw[fill] (1,0) circle (\p) node(s''n) {};
\draw[fill] (2,0) circle (\p) node(s'n) {};
\draw[fill] (3,0) circle (\p) node(sn) {};

\node (dots) at (0,0) {\bfseries $\dots$};

\draw [-] (s1) edge (s2) (s2) edge (s3) (s''n) edge (s'n) (s'n) edge (sn);

\end{tikzpicture}

& $\tilde{A_n}$
&
\begin{tikzpicture}
\def \p {0.05}
\def \op {1}
\def \gris {black!10}

\draw[fill] (-3,0) circle (\p) node(s1) {};
\draw[fill] (-2,0) circle (\p) node(s2) {};
\draw[fill] (-1,0) circle (\p) node(s3) {};
\draw[fill] (1,0) circle (\p) node(s''n) {};
\draw[fill] (2,0) circle (\p) node(s'n) {};
\draw[fill] (3,0) circle (\p) node(sn) {};
\draw[fill] (0,0.5) circle (\p) node(s0) {};

\node (dots) at (0,0) {\bfseries $\dots$};

\draw [-] (s1) edge (s2) (s2) edge (s3) (s''n) edge (s'n) (s'n) edge (sn) (sn) edge (s0) (s0) edge (s1);

\end{tikzpicture} \\

\hline

$B_n, n \geq 2$
&
\begin{tikzpicture}
\def \p {0.05}
\def \op {1}
\def \gris {black!10}

\draw[fill] (-3,0) circle (\p) node(s1) {};
\draw[fill] (-2,0) circle (\p) node(s2) {};
\draw[fill] (-1,0) circle (\p) node(s3) {};
\draw[fill] (1,0) circle (\p) node(s''n) {};
\draw[fill] (2,0) circle (\p) node(s'n) {};
\draw[fill] (3,0) circle (\p) node(sn) {};

\node (dots) at (0,0) {\bfseries $\dots$};

\draw [-] (s1) edge (s2) (s2) edge (s3) (s''n) edge (s'n) (s'n) edge (sn);

\node (label) at (2.5,0.2) {$4$};

\end{tikzpicture}

& $\tilde{B_n}$
&
\begin{tikzpicture}
\def \p {0.05}
\def \op {1}
\def \gris {black!10}

\draw[fill] (-3,0.25) circle (\p) node(s1) {};
\draw[fill] (-3,-0.25) circle (\p) node(s'1) {};
\draw[fill] (-2,0) circle (\p) node(s2) {};
\draw[fill] (-1,0) circle (\p) node(s3) {};
\draw[fill] (1,0) circle (\p) node(s''n) {};
\draw[fill] (2,0) circle (\p) node(s'n) {};
\draw[fill] (3,0) circle (\p) node(sn) {};

\node (dots) at (0,0) {\bfseries $\dots$};

\draw [-] (s'1) edge (s2) (s1) edge (s2) (s2) edge (s3) (s''n) edge (s'n) (s'n) edge (sn);

\node (label) at (2.5,0.2) {$4$};

\end{tikzpicture}\\

\hline

$C_n=B_n, n \geq 2$
&
\begin{tikzpicture}
\def \p {0.05}
\def \op {1}
\def \gris {black!10}

\draw[fill] (-3,0) circle (\p) node(s1) {};
\draw[fill] (-2,0) circle (\p) node(s2) {};
\draw[fill] (-1,0) circle (\p) node(s3) {};
\draw[fill] (1,0) circle (\p) node(s''n) {};
\draw[fill] (2,0) circle (\p) node(s'n) {};
\draw[fill] (3,0) circle (\p) node(sn) {};

\node (dots) at (0,0) {\bfseries $\dots$};

\draw [-] (s1) edge (s2) (s2) edge (s3) (s''n) edge (s'n) (s'n) edge (sn);

\node (label) at (2.5,0.2) {$4$};

\end{tikzpicture}

& $\tilde{C_n}$
&
\begin{tikzpicture}
\def \p {0.05}
\def \op {1}
\def \gris {black!10}

\draw[fill] (-3,0) circle (\p) node(s1) {};
\draw[fill] (-2,0) circle (\p) node(s2) {};
\draw[fill] (-1,0) circle (\p) node(s3) {};
\draw[fill] (1,0) circle (\p) node(s''n) {};
\draw[fill] (2,0) circle (\p) node(s'n) {};
\draw[fill] (3,0) circle (\p) node(sn) {};

\node (dots) at (0,0) {\bfseries $\dots$};

\draw [-] (s1) edge (s2) (s2) edge (s3) (s''n) edge (s'n) (s'n) edge (sn);

\node (label) at (-2.5,0.2) {$4$};
\node (label) at (2.5,0.2) {$4$};

\end{tikzpicture}\\

\hline

$D_n, n \geq 4$
&
\begin{tikzpicture}
\def \p {0.05}
\def \op {1}
\def \gris {black!10}

\draw[fill] (-3,0) circle (\p) node(s1) {};
\draw[fill] (-2,0) circle (\p) node(s2) {};
\draw[fill] (-1,0) circle (\p) node(s3) {};
\draw[fill] (1,0) circle (\p) node(s''n) {};
\draw[fill] (2,0) circle (\p) node(s'n) {};
\draw[fill] (3,0.25) circle (\p) node(sn) {};
\draw[fill] (3,-0.25) circle (\p) node(sn') {};

\node (dots) at (0,0) {\bfseries $\dots$};

\draw [-] (s1) edge (s2) (s2) edge (s3) (s''n) edge (s'n) (s'n) edge (sn) (s'n) edge (sn');

\end{tikzpicture}

& $\tilde{D_n}$
&
\begin{tikzpicture}
\def \p {0.05}
\def \op {1}
\def \gris {black!10}

\draw[fill] (-3,0.25) circle (\p) node(s1) {};
\draw[fill] (-3,-0.25) circle (\p) node(s'1) {};
\draw[fill] (-2,0) circle (\p) node(s2) {};
\draw[fill] (-1,0) circle (\p) node(s3) {};
\draw[fill] (1,0) circle (\p) node(s''n) {};
\draw[fill] (2,0) circle (\p) node(s'n) {};
\draw[fill] (3,0.25) circle (\p) node(sn) {};
\draw[fill] (3,-0.25) circle (\p) node(sn') {};

\node (dots) at (0,0) {\bfseries $\dots$};

\draw [-] (s'1) edge (s2) (s1) edge (s2) (s2) edge (s3) (s''n) edge (s'n) (s'n) edge (sn) (s'n) edge (sn');

\end{tikzpicture}\\

\hline

\end{tabular}
\caption{Spherical and Euclidean irreducible diagrams of type $A$, $B$, $C$ and $D$}
\label{tab:classification}
\end{center}
\end{table}

We recall in Table~\ref{tab:classification} the classification of the four infinite families of spherical and Euclidean irreducible diagrams, see~\cite{bourbaki_lie_456} for the full classification. We only present those because we will only consider these types in this article. Note that we use in this table the convention of Dynkin diagrams: vertices that are not joined by an edge commute, and we drop the label $3$ from edges.

\mk

Artin groups are closely related to hyperplane complements, which can be presented in a simple way in spherical and Euclidean types. Fix a Coxeter group $W=W(\Gamma)$ of spherical type or Euclidean type acting by reflections on a sphere $\SS^{n-1}$ or a Euclidean space $\R^n$. In case $W$ is of spherical type, consider the associated linear action on $\R^n$. A conjugate of an element of the standard generating set $S$ is called a \emph{reflection}. Let ${\cal R}$ denote the set of reflections in $W$. Consider the family of affine hyperplanes of $\R^n$
$${\cal H} = \{H_r \st r \in {\cal R}\},$$
where $H_r \subset \R^n$ denotes the fixed point set of the reflection $r$.

\mk

The complement of the complexified hyperplane arrangement is
$$M(\Gamma) = \C^n \bs \bigcup_{r \in {\cal R}} (\C \otimes H_r).$$
Note that $W$ acts naturally on $M$, and we have the following (see~\cite{vanderlek}):
$$\pi_1(W(\Gamma) \bs M(\Gamma)) \simeq A(\Gamma).$$
So the Artin group $A(\Gamma)$ appears as the fundamental group of (a quotient of) the complement of a complexified hyperplane arrangement. One very natural question is to decide whether it is a classifying space. This is the statement of the following conjecture.

\bconj[$K(\pi,1)$ conjecture]
The space $M(\Gamma)$ is aspherical.
\econj

This conjecture has been proved for spherical type Artin groups by Deligne in~\cite{deligne}, and for Euclidean type Artin groups by Paolini and Salvetti in~\cite{paolini_salvetti} very recently, even for the type $\tilde{D_n}$. Another approach, more closely related to the content of this article, has been used by Charney and Davis in~\cite{charney_davis} to prove the $K(\pi,1)$ conjecture for Artin groups of type FC or of $2$-dimensional type. Their proof relies on the use of a simplicial complex, called the Deligne complex, and endow it with a particular metric to show that it is contractible.

\mk

We will now recall the definition of the Deligne complex of an Artin group $A=A(\Gamma)$. A \emph{standard parabolic} subgroup of $A$ is the subgroup $A_T=\<T\>$ generated by a subset $T$ of $S$, the standard generating set of $A$. A \emph{parabolic} subgroup denotes any conjugate of a standard parabolic subgroup. Let us denote
$${\cal S}_f = \{T \subset S \st W_T \emph{ is finite}\}.$$
The \emph{Deligne complex} $\Delta=\Delta(\Gamma)$ is the order complex of the set of cosets of parabolic subgroups
$$\{gA_T \st g \in A,  T \in {\cal S}_f\},$$
where the partial order is given by the inclusion $gA_T \subset g'A_{T'}$ of cosets.

\mk

One key property of the Deligne complex is that it has the same homotopy type as the universal cover of the hyperplane complement:

\bthm \cite[Theorem~1.5.1]{charney_davis_kpi1} \label{thm:deligne_equivalent_hyperplane_arrangement}
The Deligne complex $\Delta(\Gamma)$ is homotopy equivalent to the universal cover of the quotient of the hyperplane complement $W(\Gamma) \bs M(\Gamma)$.
\ethm

In particular, the $K(\pi,1)$ conjecture amounts to proving that the Deligne complex is contractible.

\subsection{Injective metrics and Helly graphs}

We briefly recall basic definitions of injective metric spaces and Helly graphs. We will also state the local-to-global result for Helly graphs and deduce the analogous Cartan-Hadamard result for injective metric spaces.

\mk

A geodesic metric space is called \emph{injective} (one may also say \emph{hyperconvex}, or \emph{absolute $1$-Lipschitz retract}) if any family of pairwise intersecting closed balls has a non-empty global intersection, the so-called \emph{Helly property}. We refer the reader to~\cite{lang} for a presentation of injective metric spaces.

\bexes\
\bit
\item The normed vector space $(\R^n,\ell^\infty)$ is injective for all $n \geq 1$. In fact, it is up to isometry the only injective norm on $\R^n$ (see~\cite{nachbin}).
\item Any tree is injective.
\item Any finite-dimensional CAT(0) cube complex, endowed with the standard piecewise $\ell^\infty$ metric, is injective (see~\cite{miesch_CCC} and \cite{bowditch_median_injective}).
\eit
\eexes

\mk

A connected graph $X$ is called \emph{Helly} if any family of pairwise intersecting combinatorial balls of $X$ has a non-empty global intersection, in other words the combinatorial balls satisfy the Helly property. We refer the reader to~\cite{helly_groups} for a presentation of Helly graphs and Helly groups. Many examples of Helly graphs come from thickening of complexes, which we define now (see~\cite{helly_groups,huang_osajda_helly}).

\bdf[Thickening] \label{def:thickening}
Let $X$ denote a cell complex. The \emph{thickening} of $X$ (with respect to the cell structure) is the graph with vertex set $X^{(0)}$, with an edge between two vertices if and only if they are contained in a common cell of $X$.
\edf

\bexes\
\bit
\item For each $n \geq 1$, the graph with vertex set $\Z^n$, with an edge between $v,w$ if $d_\infty(v,w) = 1$, is a Helly graph.
\item Any simplicial tree is a Helly graph.
\item The thickening of the vertex set of a CAT(0) cube complex is a Helly graph.
\eit
\eexes

\mk

In order to endow a simplicial complex with a potentially injective (or CAT(0)) metric, it is natural to ask for metric simplices which may tile the Euclidean space $\R^n$. One choice is to consider the barycentric subdivision of the standard cubical tiling of $\R^n$, whose simplices are orthosimplices, which we now formally define.

\begin{figure}
\begin{center}
\begin{tikzpicture}
\def \p {0.05}
\def \op {0.5}
\def \gris {black!10}
\draw[fill] (0,0) circle (\p) node(0) {};
\draw[fill] (3,0) circle (\p) node(1) {};
\draw[fill] (3,3) circle (\p) node(2) {};
\draw[fill] (3,3) + (20:3) circle (\p) node(3) {};

\draw[black,fill opacity=\op,fill=\gris] (0.center) -- (1.center) -- (3.center) -- (0.center);
\draw[black,fill opacity=\op,fill=\gris] (0.center) -- (2.center) -- (3.center) -- (0.center);
\draw[black,fill opacity=\op,fill=\gris] (0.center) -- (1.center) -- (2.center) -- (0.center);
\draw[black,fill opacity=\op,fill=\gris] (3.center) -- (1.center) -- (2.center) -- (3.center);
\node at ([yshift=-0.5cm]0) {\bfseries $v_0=(0,0,0)$};
\node at ([yshift=-0.5cm]1) {\bfseries $v_1=(1,0,0)$};
\node at ([yshift=0.7cm]2) {\bfseries $v_2=(1,1,0)$};
\node at ([yshift=0.5cm]3) {\bfseries $v_3=(1,1,1)$};

\draw (2.5,0) -- (2.5,0.5) -- (3,0.5);
\draw (2.7,0) -- (2.95,0.35) -- (3.25,0.35);
\draw (3,2.5) -- (3.45,2.65) -- (3.45,3.15);
\draw (2.6,2.6) -- (3.1,2.85) -- (3.4,3.15);

\end{tikzpicture}
\end{center}
\caption{The standard $3$-dimensional orthosimplex.}
\label{fig:orthosimplex}
\end{figure}
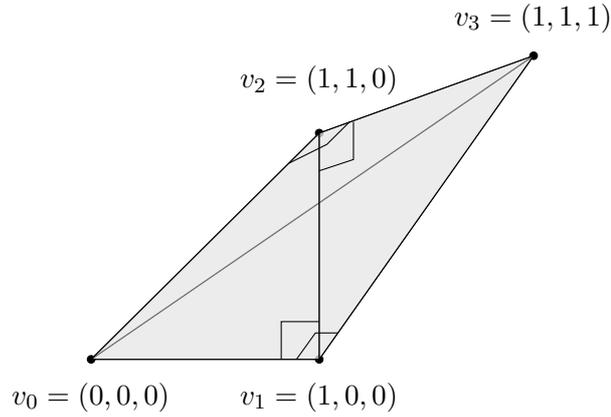

\bdf[Orthosimplex]
The \emph{standard orthosimplex} of dimension $n$ is the simplex of $\R^n$ with vertices $(0,\dots,0),(1,0,\dots,0),\dots,(1,1,\dots,1)$ (see Figure~\ref{fig:orthosimplex}). One may endow the simplex with the standard $\ell^p$ metric on $\R^n$, for any $p \in [1,\infty]$. Throughout this article (except in Theorem~\ref{thm:orthoscheme_cat0_vs_injective}), we will only consider the $\ell^\infty$ metric, called the $\ell^\infty$ orthosimplex. Note that any $n$-simplex with a total order on its vertices $v_0<v_1<\dots<v_n$ may be identified uniquely with the $\ell^\infty$ orthosimplex of dimension $n$, where each $v_i$ is identified with the vertex $(1,\dots,1,0,\dots,0)$ with $i$ ones and $n-i$ zeros. Also note that reversing the total order on the vertices gives rise to an isometry of the orthosimplex. 
\edf

\bdf[Simplicial complex with ordered simplices and maximal edges] \label{def:ordered_simplices_maximal_edges}
A simplicial complex $X$ is called \emph{with ordered simplices} if each simplex of $X$ has a total order on its vertex set, which is consistent with respect to inclusions of simplices of $X$. Moreover, we say that $X$ has \emph{maximal edges} if, given any simplex $\sigma$ of $X$, there exist adjacent vertices $a,b$ in $X$ such that $\sigma \cup \{a,b\}$ is a simplex of $X$, and such that any vertex $c$ of $X$ adjacent to both $a$ and $b$ satisfies $a \leq c \leq b$. Such an edge $\{a,b\}$ is then called a maximal edge of $X$.\edf

\bdf[Orthoscheme complex]
Let $X$ denote a finite-dimensional simplicial complex with ordered simplices. Then one may endow each simplex of $X$ with the associated $\ell^\infty$ orthoscheme metric. Then the geometric realization of $|X|$, endowed with the length metric associated with the $\ell^\infty$ orthoscheme metric on each simplex, is called the \emph{$\ell^\infty$ orthoscheme complex} of $X$.
\edf

As a particular case, if $L$ is a poset with a bound on the length of its chains, then the geometric realization of $|L|$ satisfies the assumptions, so we may talk about the $\ell^\infty$ orthoscheme complex of $L$. The following is an immediate adaptation of~\cite[Theorem~7.13]{bridson_haefliger}.

\bthm \label{thm:orthoscheme_complete_length}
Let $X$ denote a finite-dimensional simplicial complex with ordered simplices. Then its $\ell^\infty$ orthoscheme complex is a complete length space.
\ethm

\bp
Note that since $X$ has dimension $n$, the $\ell^\infty$ orthoscheme complex of $X$ has finitely many isometry types of cells: the standard orthoscheme $k$-simplices, where $k \leq n$. The proof of~\cite[Theorem~7.13]{bridson_haefliger} adapts without change to this situation.
\ep

One key property in the study of Helly graphs is the following local-to-global statement (see~\cite{chalopin_chepoi_hirai_osajda}). Recall that a \emph{clique} of a graph is a complete subgraph. A graph is called \emph{clique-Helly} if its family of maximal cliques satisfies the Helly property (i.e. any family of pairwise intersecting maximal cliques has non-empty intersection). The triangle complex of a simplicial graph $X$ is the simplicial $2$-complex whose $1$-skeleton is $X$, and whose $2$-simplices correspond to triangles in $X$.

\bthm \cite[Theorem~3.5]{chalopin_chepoi_hirai_osajda} \label{thm:clique_helly_global_helly}
Let $X$ denote a graph. Then $X$ is Helly if and only if $X$ is clique-Helly and its triangle complex is simply connected.
\ethm

In order to transfer this local-to-global property to injective metric spaces, we will need the following technical lemma. Say that a metric space is \emph{$\eps$-coarsely injective} (for some $\eps \geq 0$) if, for any families $(x_i)_{i \in I}$ in $X$ and $(r_i)_{i \in I}$ in $\R_+$ such that $\forall i,j \in I, d(x_i,x_j) \leq r_i+r_j$, we have $\bigcap_{i \in I} B_X(x_i,r_i+\eps) \neq \emptyset$. Note that when $\eps=0$, we recover the definition of an injective metric space.

\blem \label{lem:epsilon_injective}
Let $X$ denote a complete metric space that is $\eps$-coarsely injective for every $\eps>0$. Then $X$ is injective.
\elem

\bp We will prove that $X$ is injective: consider a family $(B_X(x_i,r_i))_{i \in I}$ of pairwise intersecting balls in $X$. We know that $\forall i,j \in I, d(x_i,x_j) \leq r_i+r_j$. For any $\eps>0$, let us denote $A_\eps= \bigcap_{i \in I} B_X(x_i,r_i+\eps)$, which is non-empty by assumption $\eps$-coarse injectivity.

Fix $0 < \eps \leq \eps'$, we will prove that the Hausdorff distance between $A_\eps$ and $A_{\eps'}$ is at most $\eps+\eps'$. Note that $A_\eps \subset A_{\eps'}$. Fix $x_0 \in A_{\eps'}$, we will prove that $d(x_0,A_\eps) \leq \eps+\eps'$. Let $I_0=I \sqcup \{0\}$, and let $r_0=\eps'$. Consider the families $(x_i)_{i \in I_0}$ in $X$ and $(r_i)_{i \in I_0}$ in $\R_+$. For each $i,j \in I_0$, we know that $d(x_i,x_j) \leq r_i+r_j$: indeed, for any $i \in I$, we have $x_0 \in B_X(x_i,r_i+\eps')$. By $\eps$-coarse injectivity, we deduce that the intersection $\bigcap_{i \in I_0} B_X(x_i,r_i+\eps)$ is not empty. In particular, the ball $B_X(x_0,\eps+\eps')$ intersects $A_\eps= \bigcap_{i \in I} B_X(x_i,r_i+\eps)$. This implies that $d(x_0,A_\eps) \leq \eps+\eps'$. So we have proved that the Hausdorff distance between $A_\eps$ and $A_{\eps'}$ is at most $\eps+\eps'$.

For each $n \in \N$, consider by induction $x_n \in A_{2^{-n}}$ such that $\forall n \geq 0, d_X(x_{n+1},x_n) \leq 2^{-n}+2^{-(n+1)} \leq 2^{-n+1}$. For each $0 \leq n \leq m$, we have $d_X(x_n,x_m) \leq 2^{-n+2}$, hence the sequence $(x_n)_{n \in \N}$ is a Cauchy sequence in $X$. Since $X$ is complete, it has a limit $y \in X$. For each $n \in \N$, we have $y \in A_{2^{-n}}$, so for each $i \in I$ we have $d_X(y,x_i) \leq r_i+2^{-n}$. We deduce that, for each $i \in I$, we have $d_X(y,x_i) \leq r_i$. In other words, $y$ belongs to the intersection $\bigcap_{i \in I} B(x_i,r_i)$: we have proved that $X$ is injective.
\ep

Say that a metric space is \emph{uniformly locally injective} if there exists $\eps>0$ such that each ball of radius $\eps$ is injective. Say that a metric space is \emph{semi-uniformly locally injective} if there exists $\eps>0$ such that each ball of radius $\eps$ is uniformly locally injective. For instance, any locally compact, locally injective metric space is semi-uniformly locally injective.

\mk

As a concrete example, if $X$ denotes the injective hull of the hyperbolic plane $\H^2_\R$ and $\Gamma$ is a non-uniform lattice in $\PGL(2,\R)=\Isom(\H^2_\R)$, then the quotient $\Gamma \bs X$ is semi-uniformly locally injective, but not uniformly locally injective. There are similar examples in higher rank, see~\cite{haettel_injective_buildings}.

There are also non-locally compact examples: let $X$ denote a metric simplicial graph, such that the systole at each vertex of $X$ is bounded below. Then $X$ is locally uniformly locally injective. However, if the systole of $X$ is $0$, then $X$ is not uniformly locally injective.

\mk

We now prove a Cartan-Hadamard theorem for such injective metric spaces, relying on the local-to-global property for Helly graphs, Theorem~\ref{thm:clique_helly_global_helly}. Note that this statement generalizes~\cite[Theorem~1.2]{miesch} without the local compactness assumption.

\bthm[Cartan-Hadamard for injective metric spaces] \label{thm:cartan_hadamard_injective}
Let $X$ denote a complete, simply connected, semi-uniformly locally injective metric space. Then $X$ is injective.
\ethm

\bp
We will first prove the statement when $X$ is uniformly locally injective. Fix $\eps>0$ small enough such that balls in $X$ of radius at most $2\eps$ are injective. Consider the graph $\Gamma_\eps$ with vertex set $X$, and with an edge between $x,y \in X$ if $d(x,y) \leq \eps$. Since $X$ is geodesic, $\Gamma_\eps$ is a connected graph. Also note that, for any $x \in X$ and $n \in \N$, we have
$$B_{\Gamma_\eps}(x,n) = B_X(x,n\eps).$$

\mk

We will prove that, for each $\eps>0$, the graph $\Gamma_\eps$ is a Helly graph by applying Theorem~\ref{thm:clique_helly_global_helly}.

\mk

We first prove that the family of combinatorial $1$-balls in $\Gamma_\eps$ satisfies the Helly property. Fix a family of vertices $(x_i)_{i \in I}$ of $\Gamma_\eps$ such that $\forall i,j \in I, d_{\Gamma_\eps}( x_i,x_j) \leq 2$. We want to prove that these balls intersect in $\Gamma_\eps$.

The family of metric balls $(B_X(x_i,\eps))_{i \in I}$ in $X$ pairwise intersects: since such balls have the Helly property by assumption on $X$, so we deduce that there exists $y \in X$ such that $\forall i \in I, d_X(x_i,y) \leq \eps$. In other words, the vertex $y \in \Gamma_\eps$ lies in the intersection of all combinatorial $1$-balls $(B_{\Gamma_\eps}(x_i,1))_{i \in I}$.

\mk

We now deduce that $\Gamma_\eps$ is clique-Helly: fix a family of pairwise intersecting maximal cliques $(\sigma_i)_{i \in I}$ of $\Gamma_\eps$. Then the family of combinatorial $1$-balls centered at each vertex of each clique $\sigma_i$, for $i \in I$, pairwise intersects: according to the previous paragraph, we deduce that there exists a vertex $y \in \Gamma_\eps$ adjacent to each vertex of each clique $\sigma_i$, for $i \in I$. Since each clique $\sigma_i$ is maximal, we deduce that $y$ belongs to the intersection of all $\sigma_i$, for $i \in I$. The graph $\Gamma_\eps$ is clique-Helly.

\mk

We now prove that the triangle complex of $\Gamma_\eps$ is simply connected. Fix a combinatorial loop $\ell$ in $\Gamma_\eps$. Since $X$ is simply connected, there exists a disk $D$ in $X$ bounding $\ell$. Consider a triangulation $T$ of $D$ such that triangles have diameter for $d_X$ at most $\eps$. Then the vertex set of each triangle of $T$ is a clique in $\Gamma_\eps$, therefore $\ell$ is null-homotopic in the triangle complex of $\Gamma_\eps$. So the triangle complex of $\Gamma_\eps$ is simply connected.

\mk

The graph $\Gamma_\eps$ is clique-Helly and has a simply connected triangle complex, so according to Theorem~\ref{thm:clique_helly_global_helly}, we deduce that $\Gamma_\eps$ is a Helly graph.

\mk

Note that, for any $\eps>0$, we have $d_X \leq \eps d_{\Gamma_\eps} \leq d_X + \eps$, and balls for the metric $d_{\Gamma_\eps}$ with integral radius satisfy the Helly property.

We will show that the metric space $X$ is $\eps$-coarsely injective: consider families $(x_i)_{i \in I}$ in $X$ and $(r_i)_{i \in I}$ in $\R_+$ such that $\forall i,j \in I, d_X(x_i,x_j) \leq r_i+r_j$. For each $i \in I$, let $n_i \in \N$ such that $n_i\eps \leq r_i < (n_i+1)\eps$. For each $i,j \in I$, since the balls $B_X(x_i,r_i)$ and $B_X(x_j,r_j)$ intersect, we deduce that $d_{\Gamma_\eps}(x_i,x_j) \leq n_i+n_j+2$. So, in the Helly graph $\Gamma_\eps$, the balls $(B_{\Gamma_\eps}(x_i,n_i+1))_{i \in I}$ pairwise intersect: we deduce that there exists $y \in \Gamma_\eps$ such that $\forall i \in I, d_{\Gamma_\eps}(y,x_i) \leq n_i+1$. In particular, for any $i \in I$, we have $d_X(y,x_i) \leq (n_i+1)\eps \leq r_i+\eps$. Hence $y$ belongs to each ball $B_X(x_i,r_i+\eps)$, for $i \in I$: this proves that $X$ is $\eps$-coarsely injective.

Since this holds for any small enough $\eps>0$, according to Lemma~\ref{lem:epsilon_injective}, we conclude that $X$ is injective.

\mk

We now turn to the general case, when $X$ is only semi-uniformly locally injective: there exists $\eps>0$ such that each ball of radius $\eps$ in $X$ is uniformly locally injective. According to the uniformly locally injective case, we deduce that each ball of radius $\eps$ is injective: this means that $X$ is uniformly locally injective. According to the uniformly locally injective case again, we deduce that $X$ is injective.
\ep

We now see that, under a mild assumption on a simplicial complex with ordered simplices, saying that the $\ell^\infty$ orthoscheme realization is injective is equivalent to saying that the thickening of the $1$-skeleton is Helly. Note that, if we refer to the definition of thickening as in Definition~\ref{def:thickening}, the corresponding cell structure is not the simplicial one, but a coarser cell structure whose cells correspond to intervals.

\bthm \label{thm:injective_orthoscheme_implies_helly}
Let $X$ denote a finite-dimensional simplicial complex with ordered simplices, and with maximal edges. Let $\Gamma$ denote the graph with vertex set $X^{(0)}$, and with an edge between $x,y \in X^{(0)}$ if there exist $a,b \in X^{(0)}$ and ordered triangles $a \leq x \leq b$ and $a \leq y \leq b$ in $X$. If the $\ell^\infty$ orthoscheme complex of $X$ is injective, then the thickening $\Gamma$ of $X$ is a Helly graph. \ethm

\bp
We see that maximal cliques in $\Gamma$ correspond to intervals $I_{ab}=\{x \in X^{(0)} \st a \leq x \leq b \mbox{ is an ordered triangle in $X$}\}$, for any maximal edge $a \leq b$ in $X$. Let $m_{ab} \in |X|$ denote the midpoint of the maximal edge $a \leq b$. By assumption, any simplex of $X$ containing $a$ and $b$ has for vertex set a chain from $a$ to $b$. Then $B_{|X|}(m_{ab},\f{1}{2})$ is a subcomplex of $|X|$ with vertex set $B_{|X|}(m_{ab},\f{1}{2}) \cap X^{(0)} = I_{ab}$.

\mk

We will prove that $\Gamma$ is clique-Helly: let $(I_{a_ib_i})_{i \in I}$ denote a family of pairwise intersecting maximal cliques in $\Gamma$. Since $|X|$ is injective, there exists $z \in \bigcap_{i \in I} B_{|X|}(m_{a_ib_i},\f{1}{2})$. Since each such ball is a subcomplex of $|X|$, we may assume that $z$ is a vertex of $X$. We deduce that $z$ belongs to each clique $I_{a_ib_i}$, for $i \in I$. So $\Gamma$ is clique-Helly.

\mk

We will now prove that the triangle complex of $\Gamma$ is simply connected. Let $\ell$ denote a combinatorial loop in the $1$-skeleton of $\Gamma$. Up to homotopy in the triangle complex of $\Gamma$, we may assume that $\ell$ lies in the $1$-skeleton of $X$. Since $|X|$ is injective, it is contractible, so the $2$-skeleton of $X$ is simply connected. As the $2$-skeleton of $X$ is contained in the triangle complex of $\Gamma$, we conclude that $\ell$ is null-homotopic in the triangle complex of $\Gamma$.

\mk

According to Theorem~\ref{thm:clique_helly_global_helly}, we deduce that $\Gamma$ is Helly.
\ep

\section{The thickening of a lattice} \label{sec:thickening_lattice}

We will explain a very simple construction of Helly graphs and injective metric spaces, starting with a lattice endowed with an action of the group $\Z$ or $\R$.

\mk

Assume that $L$ is a lattice, such that each upperly bounded subset of $L$ has a join. Assume that there is an order-preserving, increasing, continuous (with respect to the order topology on $L$), action $(f_t)_{t \in H}$ of $H=\Z$ or $\R$ on $L$, such that
$$\forall x,y \in L, \exists t \in H_+, f_{-t}(x) \leq y \leq f_t(x).$$
Let us define the following metric $d$ on $L$:
$$\forall x,y \in L, d(x,y) =\inf\{t \in H_+ \st f_{-t}(x) \leq y \leq f_t(x)\}.$$

\bthm \label{thm:lattice_Helly_injective}\

If $H=\Z$, then $(L,d)$ is a Helly graph with the combinatorial distance.

If $H=\R$, then $(L,d)$ is injective.
\ethm

\bp
We start by proving that $d$ is indeed a metric on $L$. If $x,y \in L$ and $t \in H_+$ are such that $f_{-t}(x) \leq y \leq f_t(x)$, then by applying $f_t$ and $f_{-t}$ we deduce that $f_{-t}(y) \leq x \leq f_t(y)$: the metric $d$ is symmetric.

\mk

We will now prove the triangle inequality: let $x,y,z \in L$, and for $\eps>0$, consider $t,s \in H_+$ such that $d(x,y) \leq t < d(x,y)+\eps$ and $d(y,z) \leq s < d(y,z)+\eps$. We have $f_{-t}(x) \leq y \leq f_t(x)$ and $f_{-s}(y) \leq z \leq f_s(y)$. Hence $f_{-t-s}(x) \leq f_{-s}(y) \leq z \leq f_s(y) \leq f_{s+t}(x)$, so $d(x,z) \leq t+s \leq d(x,y)+d(y,z)+2\eps$. This holds for any $\eps>0$, hence $d(x,z) \leq d(x,y)+d(y,z)$.

\mk

We will now prove that the metric is positive: assume that $x,y \in L$ are distinct, we will prove that $d(x,y)>0$.

\bit
\item Assume first that $x,y$ are comparable, for instance $x<y$. Since $\{z \in L \st z < y\}$ is open and contains $x$, by continuity of the action there exists $\eps>0$ such that, for any $|t| \leq \eps$, we have $f_t(x) <y$. In particular $d(x,y) \geq \eps$.
\item Assume now that $x,y$ are not comparable. Since $U=\{z \in L \st x \wedge y < z < x \vee y\}$ is open and contains $x$, by continuity of the action there exists $\eps>0$ such that, for any $|t| \leq \eps$, we have $f_t(x) \in U$. For any $0 \leq t < \eps$, since $x \leq f_t(x) < x \vee y$, we have $f_t(x) \not\geq y$. In particular, $d(x,y) \geq \eps$.
\eit

\mk

So we have proved that $d$ is a metric on $L$. We will now prove that balls in $L$ satisfy the Helly property.

\mk

We will first consider the case $H=\Z$. Let us define a graph $\hat{L}$ with vertex set $L$, and with an edge between $x, y \in L$ if $f_{-1}(x) \leq y \leq f_1(x)$: we will prove that the graph $\hat{L}$ is Helly. We will first prove, by induction on $k \geq 0$, that for any $x \in L$, the ball $B_{\hat{L}}(x,k)$ in the graph $\hat{L}$ coincides with the interval $I(f_{-k}(x),f_k(x))=\{y \in L \st f_{-k}(x) \leq y \leq f_k(x)\}$.

\mk

For $k \leq 1$ it is the definition of the edges of $\hat{L}$, so fix $k \geq 2$ and assume that the statement holds for $k-1$. Fix $y \in I(f_{-k}(x),f_k(x))$, we will prove that $y \in B_{\hat{L}}(x,r)$. Since $y \geq f_{-k}(x)$, we deduce that $f_1(y) \geq f_{-k+1}(x)$, and also since $y \leq f_k(x)$ we deduce that $f_{-1}(y) \leq f_{k-1}(x)$. So we have that both $f_1(y)$ and $f_{k-1}(x)$ are superior to both $f_{-1}(y)$ and $f_{-k+1}(x)$: since $L$ is a lattice, there exists some element $z \in L$ in the intersection $I(f_{-1}(y),f_1(y)) \cap I(f_{-k+1}(x),f_{k-1}(x))$. In particular, $y$ and $z$ are adjacent in $\hat{L}$, and by induction we know that $d_{\hat{L}}(z,x) \leq k-1$, so $d_{\hat{L}}(x,y) \leq k$. Conversely, it is clear that the ball $B_{\hat{L}}(x,k)$ is included in the interval $I(f_{-k}(x),f_k(x))$. So we have $B_{\hat{L}}(x,k) = I(f_{-k}(x),f_k(x))$.

\mk

By assumption on the lattice $L$, we deduce that the graph $\hat{L}$ is connected, and furthermore that $d_{\hat{L}}=d$.

\mk

Remark now that intervals in the lattice $L$ satisfy the Helly property. Fix any collection $(I(x_i,y_i))_{i \in I}$ of pairwise intersecting intervals of $L$. Fix $j_0 \in I$. For any $i \in I$, we have $x_i \leq y_{j_0}$, so the set $\{x_i, i \in I\}$ is upperly bounded. By assumption, the set $\{x_i, i \in I\}$ has a join $z \in L$, such that $z \leq y_{j_0}$. This holds for any $j_0 \in I$, so $z$ belongs to the intersection of all intervals $(I(x_i,y_i))_{i \in I}$.

\mk

Hence the graph $\hat{L}$ is connected, and its balls satisfy the Helly property: it is a Helly graph.

\mk

We now turn to the case $H=\R$, and we will prove that $(L,d)$ is an injective metric space. First note that balls in $(L,d)$ are intervals in $L$, so according to the previous argument, we know that balls in $(L,d)$ satisfy the Helly property. In order to prove that $(L,d)$ is injective, according to the definition of hyperconvex metric spaces (see for instance~\cite{lang}), it is sufficient to prove that if $x,y \in L$ and $r,s \geq 0$ are such that $d(x,y) \leq r+s$, then the balls $B(x,r)$ and $B(y,s)$ intersect. In other words, it is enough to prove that $(L,d)$ is weakly geodesic, i.e. for any $x,y \in L$ and any $0 \leq r \leq d(x,y)$, there exists $z \in B(x,r) \cap B(y,d(x,y)-r)$.

\mk

For each $k \in \N \bs \{0\}$, let us consider the action of $\f{1}{k} \Z \subset \R$ on $L$, and the associated Helly graph distance
$$\forall x,y \in L, d_k(x,y) = \inf\left\{t \in \f{1}{k} \N \st f_{-t}(x) \leq y \leq f_t(x)\right\}.$$
Fix $x,y \in L$, and $0 \leq r \leq d(x,y)$. For each $k \in \N \bs \{0\}$, there exists $z_k \in L$ such that $d_k(x,z_k) \leq r+\f{1}{k}$ and $d_k(z_k,y) \leq d_k(x,y)-r+\f{1}{k} \leq d(x,y)-r+\f{2}{k}$. So the intervals $I_k=I(f_{-r-\f{1}{k}}(x),f_{r+\f{1}{k}}(x))$ and $J_k=I(f_{-d(x,y)+r-\f{2}{k}}(y),f_{d(x,y)-r+\f{2}{k}}(y))$ in $L$ intersect. If $k \leq k'$, we know that $I_{k'} \subset I_k$ and $J_{k'} \subset J_k$. We deduce that the family of intervals $\{I_k\}_{k \in \N \bs \{0\}} \cup \{J_k\}_{k \in \N \bs \{0\}}$ pairwise intersects. By the Helly property for intervals, the global intersection is non-empty: let us denote $z$ some element in the intersection.

\mk

We know that $\liml_{k \ra \pif} d_k(x,z)=d(x,z) \leq r$ and $\liml_{k \ra \pif} d_k(y,z)=d(y,z) \leq d(x,y)-r$, so we have proved that $(L,d)$ is a weakly geodesic metric space. Since balls in $(L,d)$ satisfy the Helly property, we conclude that $(L,d)$ is injective.
\ep

An immediate consequence concerns Garside groups (see~\cite{garside}, \cite{mccammond_intro_garside} and \cite{haettel_huang_weakly_modular}). Recall that a group $G$ is called \emph{Garside} if there exists a subset $S \subset G$ and an element $\delta \in G$ such that the following hold:
\bit
\item $S$ spans the group $G$.
\item For each element in the semigroup $\<S\>^+$, there is a bound on the length of factorizations over $S$. 
\item The element $\delta$ belongs to the semigroup $\<S\>^+$, and $\delta$ is \emph{balanced}: the set of prefixes of $\delta$ coincides with the set of suffixes of $\delta$.
\item The poset of prefixes of $\delta$ in $\<S\>^+$ is a lattice.
\eit
Note that $G$ is endowed with two natural orders, the left order $\leq_L$ and the right order $\leq_R$:
\bit
\item $g \leq_L h$ if and only if $g^{-1}h \in \<S\>^+$.
\item $g \leq_R h$ if and only if $hg^{-1} \in \<S\>^+$.
\eit
Authors sometimes add the requirement that $S$ is finite, in which case $G$ may also be called a Garside group of finite type. Also note that, given a Garside group $(G,S,\delta)$, for any $g \in G$, there exists $n \in \N$ such that $\delta^{-n} \leq_L g \leq_L \delta^n$. 

\mk

Fix a Garside group $(G,S,\delta)$. Let $X$ denote the graph with vertex set $G$, with an edge between $g,h \in G$ if $g\delta^{-1} \leq_L h \leq_L g\delta$. The graph $X$ is called the \emph{thickening} of $G$. In relation to Definition~\ref{def:thickening}, it corresponds to the cell complex with vertex set $G$, whose maximals cells are translates $(g[\Delta^{-1},\Delta]_{\leq_L})_{g \in G}$ of the interval $[\Delta^{-1},\Delta]$.

\bcor \label{cor:thickening_garside_helly}
The thickening of any Garside group is a Helly graph.
\ecor

\bp
The left order on $G$ is a lattice order (see~\cite{garside}). Furthermore, consider the action of $\Z$ on $G$ by right multiplication by $\delta$. For any $g \in G$, there exists $n \in \N$ such that $\delta^{-n} \leq_L g \leq_L \delta^n$. So this action satisfies the assumptions of Theorem~\ref{thm:lattice_Helly_injective}: we deduce that the graph $X$ is Helly.
\ep

Note that this applies, in particular, to Garside groups of infinite type, such as crystallographic Garside groups (see~\cite{mccammond_sulway}). However, we do not have yet an application of this simply transitive action of a Garside group on a locally infinite Helly graph.

In the case of Garside groups of finite type, we recover a particularly simple proof of the following result by Huang and Osajda (see~\cite{huang_osajda_helly}). In particular, our proof does not rely on the deep result that local-to-global result for Helly graphs (Theorem~\ref{thm:clique_helly_global_helly}, see~\cite{chalopin_chepoi_hirai_osajda}).

\bcor[Huang-Osajda]
Any Garside group of finite type is a Helly group.
\ecor

In particular, this leads to a particularly simple proof that braids groups are Helly, relying only on some Garside structure.

\section{The affine version of a lattice} \label{sec:affine_lattice}

In this section, we will prove Theorem~\ref{thm:lattice_orthoscheme_injective} stating that the orthoscheme complex of a bounded graded lattice, endowed with the orthoscheme $\ell^\infty$ metric, is injective. In order to do so, we will apply results from the previous section, and endow the geometric realization of a lattice with a partial order, which is a lattice.

\mk

Assume that $L$ is a bounded, graded lattice of rank $n$. Let $\leq_L$ denote the order on $L$. Let $H$ denote either a cyclic subgroup of $(\R,+)$, or $H=\R$. We will define a new poset $M_H$, which will be called the \emph{affine version} of $L$ over $H$. If there is no ambiguity about $H$, we will simply denote $M=M_H$. Let $C(L)$ denote the set of maximal chains $c_{0,1}=0 <_L c_{1,2} <_L \dots c_{n-1,n} <_L c_{n,n+1}=1$ in $L$.  We will use the convention that the element denoted $c_{i,i+1}$ has rank $i$.

Let us consider the subspace
$$\sigma = \{u \in H^n \st u_1 \leq u_2 \leq \dots \leq u_n\}$$
of $H^n$.

\mk

For each maximal chain $c \in C(L)$, let $\sigma_c$ denote a copy of $\sigma$.

Let us consider the space
$$M = \bigcup_{c \in C(L)} \sigma_c / \sim,$$
where for each $c,c' \in C(L)$, if we denote $I=\{1 \leq i \leq n-1 \st c_{i,i+1} \neq c'_{i,i+1}\}$, we identify $\sigma_c$ and $\sigma_{c'}$ along the subspaces
$$\{u \in \sigma_c \st \forall i \in I, u_i=u_{i+1}\} \simeq \{u \in \sigma_{c'} \st \forall i \in I, u_i = u_{i+1}\}.$$
We can describe elements of $M$ as a quotient of the space $M_0 = C(L) \times \sigma$.

\bexe
One illustrating example is the following: consider the boolean lattice $L$ of rank $n$, i.e. the lattice of subsets of the finite set $\{1,\dots,n\}$, with the inclusion order. 
Maximal chains in $L$ correspond to permutations of $\{1,\dots,n\}$. The space $M_H$ may be identified with $H^n$, where for each permutation $w$ of $\{1,\dots,n\}$, the subspace $\sigma_w$ is
$$\sigma_w=\{x \in H^n \st x_{w(1)} \leq x_{w(2)} \leq \dots \leq x_{w(n)}\}.$$
\eexe

\mk

If $c \in C(L)$ and $u \in \sigma$, let us denote by $[c,u]$ the equivalence class of $(c,u) \in M_0$ in $M$.

\mk

For each $c \in C(L)$, let us endow $\sigma_c$ with the partial order from $H^n \subset \R^n$: $u \leq v$ if $\forall 1 \leq i \leq n, u_i \leq v_i$. Let us endow $M$ with the induced partial order: we have $\alpha \leq \beta$ in $M$ if there exists a sequence $\alpha_0=\alpha, \alpha_1,\dots,\alpha_m=\beta$ in $M$ such that, for each $0 \leq i \leq m-1$, there exists $c \in C(L)$ and $x \leq y$ in $\sigma_c$ such that $\alpha_i=[c,x]$ and $\alpha_{i+1}=[c,y]$.

Let $\leq_M$ denote the order on $M$. We will prove the following.

\bthm \label{thm:affine_poset_lattice}
If $H$ is a discrete subgroup of $\R$, the poset $M_H$ is a lattice.
\ethm

Before proving Theorem~\ref{thm:affine_poset_lattice}, we will gather some preliminary results. Without loss of generality, assume that $H=\Z$. To simplify notations, we will let $M=M_H$.

\mk

First notice that if $(a,u) \sim (b,v)$, then $u=v$. Therefore the second projection $M_0 \ra \sigma$ defines a projection $\pi: M \ra \sigma$. Fix $\beta=[b,v],\gamma=[c,w] \in M$, and fix $\alpha \leq_M \beta,\gamma$ in $M$. We will prove that $\beta,\gamma$ have a join, and $\alpha$ will play an auxiliary role.

\mk

We say that $\beta$ is \emph{elementarily superior} to $\alpha$ if there are representatives $\alpha=[a,u]$ and $\beta=[b,v]$ such that $a=b$ and there exist $1 \leq i \leq j \leq n$ such that:
\bit
\item $u_i=u_{i+1}=\dots=u_j$,
\item $v_i=v_{i+1}=\dots=v_j=u_j+1$, and
\item $\forall k \not\in [i,j], u_k=v_k$.
\eit

\blem \label{lem:elementary_superior}
Fix $\alpha=[a,u] \in M$. Any element $\beta=[b,v]$ of $M$ elementarily superior to $\alpha$ is uniquely determined by:
\bit
\item Integers $1 \leq i \leq j \leq n$, such that $u_i=u_{i+1}= \dots = u_j$, and $u_j<u_{j+1}$ if $j < n$.
\item Some element $b_{i-1,i}$ of rank $i-1$ in the interval $I(a_{i_0-1,i_0},a_{j,j+1})$, where $i_0 \in [1,i]$ is minimal such that $u_{i_0}=u_i$.
\eit
Let us denote $\beta=\alpha[i,j,b_{i-1,i}]$. 
\elem

\bp
Let us define the element $\beta=[b,v]$ in $M$ by:
\beq \forall k \not\in [i,j], v_k&=&u_k \\
\forall k \in [i,j], v_k &=& u_{j}+1 \\
\forall k \leq i_0-1, b_{k,k+1}&=&a_{k,k+1} \\
\forall k \geq j, b_{k,k+1}&=&a_{k,k+1}.\eeq
Note that $v$ is non-decreasing, hence $v \in \sigma$. Furthermore, since $v_{i_0}=v_{i_0+1}=\dots=v_{i-1}$ and $v_i = v_{i+1} = \dots = u_j+1 \leq v_{j+1}$, it is enough to define $b_{k,k+1}$ for $k \leq i_0-1$, $k=i-1$ and $k \geq j$. Hence $\beta$ is well-defined, and it is elementarily superior to $\alpha$. It is clear that $\beta$ is the only such element in $M$.
\ep

\blem \label{lem:sequence_elementary}
Given $\alpha \leq_M \beta$ in $M$, there exist $m \geq 0$ and a sequence $\beta_0=\alpha,\beta_1,\dots,\beta_m=\beta$ for which, for each $0 \leq i \leq m-1$, $\beta_{i+1}$ is elementarily superior to $\beta_i$.
\elem

\bp
According to the definition of the order on $M$, it is sufficient to prove the statement for $\alpha=[a,u]$ and $\beta=[b,v]$ such that $a=b$. We have $u \leq v$ in $\Z^n \subset \R^n$. Then consider a sequence $u^0=u \leq u^1 \leq \dots \leq u^m=v$ such that, for each $0 \leq k \leq m-1$, there exist $1 \leq i_k \leq j_k \leq n$ such that:
\bit
\item $u^k_{i_k}=u^k_{i_k+1}=\dots=u^k_{j_k}$,
\item $u^{k+1}_{i_k}=u^{k+1}_{i_k+1}=\dots=u^{k+1}_{j_k}=u^k_{j_k}+1$, and
\item $\forall \ell \not\in [i_k,j_k], u^k_\ell=u^{k+1}_\ell$.
\eit
For each $0 \leq k \leq m-1$, the element $\beta_{k+1}=[a,u^{k+1}]$ is elementarily superior to $\beta_k=[a,u^k]$. 
\ep

\blem \label{lem:two_orders_coincide}
We have $\alpha=[a,u] \leq \beta=[b,v]$ in $M$ if and only if $u \leq v$ and for every $0 \leq j \leq n-1$ such that $u_j < u_{j+1}$, if we denote by $i \in \{0,\dots,j\}$ the minimal element such that $v_{i+1} \geq u_{j+1}$, we have $b_{i,i+1} \leq_L a_{j,j+1}$.
\elem

\bp
Let us denote by $\prec$ this relation. We first show that $\prec$ is transitive and antisymmetric.

\mk

Assume that $\alpha=[a,u] \preccurlyeq \beta=[b,v]$ and $u=v$, we will prove that $\alpha = \beta$. Fix $0 \leq j \leq n-1$, and assume that $a_{j,j+1} \neq b_{j,j+1}$. Since these two elements of $L$ have the same rank $j$, we deduce that they are not comparable: $b_{j,j+1} \not\leq_L a_{j,j+1}$ and $a_{j,j+1} \not\leq_L b_{j,j+1}$. We will prove that $u_j=u_{j+1}$: by contradiction, if $u_j < u_{j+1}$, then the minimal $i \in \{0,\dots,j\}$ such that $v_{i+1} \geq u_{j+1}$ is equal to $j$, so $b_{j,j+1} \leq_L a_{j,j+1}$. So $u_j=u_{j+1}$. This proves that $\alpha=\beta$.

\mk

This implies in particular that $\prec$ is antisymmetric: indeed assume that $\alpha \preccurlyeq \beta$ and $\beta \preccurlyeq \alpha$. Then $u \leq v$ and $v \leq u$, so $u=v$, hence $\alpha=\beta$.

\mk

We now prove that $\prec$ is transitive: assume that $\alpha=[a,u] \preccurlyeq \beta=[b,v] \preccurlyeq \gamma=[c,w]$. Then $u \leq v \leq w$, so $u \leq w$. Let $0 \leq k \leq n-1$ such that $u_k < u_{k+1}$, let $i \in \{0,\dots,k\}$ be minimal such that $w_{i+1} \geq u_{k+1}$: we want to prove that $c_{i,i+1} \leq_L a_{k,k+1}$. Let $j \in \{0,\dots,k\}$ be minimal such that $v_{j+1} \geq u_{k+1}$: we know that $b_{j,j+1} \leq_L a_{k,k+1}$. Since $j$ is minimal, we know that either $j=0$ or $v_j < u_{k+1}$.

If $j=0$, then since $i \leq j$ we know that $i=0$ and so $c_{i,i+1} = 0 = b_{j,j+1} \leq_L a_{k,k+1}$.

If $v_j < u_{k+1}$ then $v_j < v_{j+1}$, and let $i' \in \{0,\dots,j\}$ be minimal such that $w_{i'+1} \geq v_{j+1}$. We know that $c_{i',i'+1} \leq_L b_{j,j+1}$. But since $w_i < u_{k+1} \leq v_{j+1}$, we have $i \leq i'$, so $c_{i,i+1} \leq_L c_{i',i'+1} \leq_L b_{j,j+1} \leq_L a_{k,k+1}$.

Hence $\alpha \preccurlyeq \gamma$.

\mk

We will now prove that if $\beta=[b,v]$ is elementarily superior to $\alpha=[a,u]$, we have $\alpha \preccurlyeq \beta$. Let us denote $\beta=\alpha[i,j,b_{i-1,i}]$. First note that $u \leq v$. If $0 \leq k \leq n-1$ is such that $u_k<u_{k+1}$, then either $k \leq i-1$ or $k \geq j$. Let $\ell \in \{0,\dots,k\}$ denote the minimal element such that $v_{\ell+1} \geq u_{k+1}$, we will prove that $b_{\ell,\ell+1} \leq_L a_{k,k+1}$.

\bit
\item If $k \leq i-1$, then $\ell=k$ and $b_{\ell,\ell+1}=a_{k,k+1}$.
\item If $k \geq j$ and $v_{k+1} = v_j$, then $\ell=i-1$ and $b_{\ell,\ell+1}=b_{i-1,i} \leq_L a_{j,j+1} \leq_L a_{k,k+1}$.
\item If $k \geq j$ and $v_{k+1} > v_j$, then $\ell \geq j$ and $b_{\ell,\ell+1} =a_{\ell,\ell+1} \leq_L a_{k,k+1}$.
\eit

\mk

According to Lemma~\ref{lem:sequence_elementary}, we deduce that if $\alpha \leq_M \beta$, then $\alpha \preccurlyeq \beta$.

\mk

We will now prove that if $\alpha=[a,u] \preccurlyeq \beta=[b,v]$ with $\alpha \neq \beta$, there exists $\gamma$ elementarily superior to $\alpha$ such that $\gamma \preccurlyeq \beta$. If $u=v$, we have seen that $\alpha=\beta$, so let us assume that $u < v$. Let $1 \leq i \leq n$ be minimal such that $u_i < v_i$. Let $j \in \{i,\dots,n\}$ be maximal such that $u_i=u_j$.

Let $i_0 \in \{1,\dots,i\}$ be minimal such that $u_{i_0}=u_i$. Since for all $k \leq i-1$ we have $u_k=v_k$, there exist representatives $\alpha=[a,u]$ and $\beta=[b,v]$ such that, for any $1 \leq k \leq i-1$, we have $a_{k-1,k} = b_{k-1,k}$.

\mk

We will also prove that we may assume that $a_{i-1,i}=b_{i-1,i}$.

Assume first that $u_{i-1}<u_i$. Let $p \in \{0,\dots,i-1\}$ be minimal such that $v_{p+1} \geq u_i$. Since $v_{i-1}=u_{i-1}<u_i$, we deduce that $p=i-1$. Since $\alpha \preccurlyeq \beta$, we have $b_{i-1,i} \leq_L a_{i-1,i}$. As both elements have the same rank $i-1$ in $L$, we have $a_{i-1,i}=b_{i-1,i}$.

Assume now that $u_{i-1}=u_i$. In the case $j<n$, we have $u_j < u_{j+1}$: let $p \in \{0,\dots,j\}$ be minimal such that $v_{p+1} \geq u_{j+1}$. 
Since $v_{i-1}=u_{i-1}<u_{j+1}$, we deduce that $p \geq i-1$. Since $\alpha \preccurlyeq \beta$, we have $b_{i-1,i} \leq_L b_{p,p+1} \leq_L a_{j,j+1}$. Since $u_{i-1}=u_i=\dots=u_j$, we may choose a representative of $\alpha$ such that $a_{i-1,i}=b_{i-1,i}$. In the case $j=n$, we have $u_{i-1}=u_i= \dots =u_n$, we may also choose a representative of $\alpha$ such that $a_{i-1,i}=b_{i-1,i}$.

We have proved that we may assume that $a_{i-1,i}=b_{i-1,i}$.

\mk

In the case $j=n$, since $u_i= \dots =u_n$, we may choose a representative of $\alpha$ such that $a=b$. In this case, the element $\gamma=\alpha[i,n,b_{i-1,i}]$ is elementarily superior to $\alpha$, and $\gamma \preccurlyeq \beta$.

\mk

For the rest of the proof, assume that $j<n$. 

We know thar $b_{i-1,i } \geq_L b_{i_0-1,i_0}=a_{i_0-1,i_0}$. Since $u_j < u_{j+1}$ and $\alpha \preccurlyeq \beta$, if we denote by $i' \in \{0,\dots,j\}$ the minimal element such that $v_{i'+1} \geq u_{j+1}$, we have $b_{i',i'+1} \leq_L a_{j,j+1}$. Since $v_{i-1} =u_{i-1} < u_{j+1}$, we have $i' >i-2$, so $b_{i-1,i} \leq_L b_{i',i'+1} \leq_L a_{j,j+1}$.

So we know that $b_{i-1,i} \in I(a_{i_0-1,i_0},a_{j,j+1})$: we can define $\gamma=\alpha[i,j,b_{i-1,i}]=[w,c]$. We will now prove that $\gamma \preccurlyeq \beta$. Fix $0 \leq k \leq n-1$ such that $w_k < w_{k+1}$, and denote by $\ell \in \{0,\dots,k\}$ the minimal element such that $v_{\ell+1} \geq w_{k+1}$, we will prove that $b_{\ell,\ell+1} \leq_L c_{k,k+1}$. Recall that either $k \leq i-1$ or $k \geq j$.

\bit
\item If $k \leq i-2$, then $w_{k+1}=u_{k+1}$. Since $\alpha \preccurlyeq \beta$, we know that $b_{\ell,\ell+1} \leq_L a_{k,k+1} = c_{k,k+1}$.
\item If $k=i-1$, then $w_{k+1}=w_i = u_i+1 > u_i$. Then $v_{i-1} = u_{i-1} \leq u_i < w_{k+1}$, so $\ell=i-1=k$. Hence $b_{\ell,\ell+1} = b_{i-1,i} = a_{i-1,i} = c_{i-1,i} = c_{k,k+1}$.
\item If $k \geq j$, then $w_{k+1}=u_{k+1}$. Since $\alpha \preccurlyeq \beta$, we know that $b_{\ell,\ell+1} \leq_L a_{k,k+1}$. We also have $a_{k,k+1}=c_{k,k+1}$, so $b_{\ell,\ell+1} \leq_L c_{k,k+1}$.
\eit

So we conclude that $\gamma \preccurlyeq \beta$.

\mk

Now fix any $\alpha=[a,u] \preccurlyeq \beta=[b,v]$. By induction on $\|v-u\|_1$, we see that there is a bound on sequences of elementarily superior elements starting from $\alpha$ which all are $\preccurlyeq \beta$. Therefore we conclude that $\alpha \leq_M \beta$.

\mk

In conclusion, the two orders $\leq_M$ and $\prec$ coincide.
\ep

\mk

Given $\alpha \leq_M \beta$ in $M$, let us denote by $D(\alpha,\beta)$ the minimal number $m \geq 0$ such that there exists a sequence $\beta_0=\alpha,\beta_1,\dots,\beta_m=\beta$ for which, for each $0 \leq i \leq m-1$, $\beta_{i+1}$ is elementarily superior to $\beta_i$.

\mk

We will prove, by induction on $D(\alpha,\beta)+D(\alpha,\gamma)$, that $\beta$ and $\gamma$ have a join.

\mk

\blem \label{lem:join_distance_1}
Assume that $\alpha,\beta,\gamma \in M$ are such that $\alpha \leq_M \beta$, $\alpha \leq_M \gamma$ and $D(\alpha,\beta)=D(\alpha,\gamma)=1$. Then $\beta,\gamma$ have a join $\delta$ such that $D(\beta,\delta) \leq 1$ and $D(\gamma,\delta) \leq 1$.
\elem

\bp

Consider representatives $\alpha=[a,u]$, $\beta=[b,v]$ and $\gamma=[c,w]$ of $\alpha$, $\beta$ and $\gamma$ respectively. According to Lemma~\ref{lem:elementary_superior}, there exist $1 \leq i \leq j \leq n$ and $b_{i,i+1}$ such that $\beta=\alpha[i,j,b_{i,i+1}]$, and there exist $1 \leq i' \leq j' \leq n$ and $c_{i',i'+1}$ such that $\gamma=\alpha[i',j',c_{i',i'+1}]$.

\mk

{\bf First case:} assume that the intervals $[i,j]$ and $[i',j']$ are disjoint, for instance $j < i'$. Let us define $\delta = \beta[i',j',c_{i',i'+1}]$: we will see that $\delta$ is well-defined and that $\delta=\gamma[i,j,b_{i,i+1}]$.

\mk

First note that $v_{j'} = u_{j'} < u_{j'+1} = v_{j'+1}$. Furthermore, let $i'_0 \in [1,i']$ be minimal such that $u_{i'_0}=u_{i'}$. Since $u_j < u_{j+1}$, we know that $j+1 \leq i'_0$. So we have $v_{i'_0} = u_{i'_0} = u_{i'} = v_{i'}$. So, if we denote $i''_0 \in [1,i']$ the minimal integer such that $v_{i''_0} = v_{i'}$, we have $i''_0 \leq i'_0$.

By definition we have $c_{i',i'+1} \in I(a_{i'_0-1,i'_0},a_{j',j'+1})$. Since $j < j'$, we have $b_{j',j'+1}=a_{j',j'+1}$, so $c_{i',i'+1} \leq_L b_{j',j'+1}$. And as $b_{i''_0-1,i''_0} \leq_L b_{i'_0-1,i'_0} = a_{i'_0-1,i'_0}$, we deduce that $c_{i',i'+1} \geq_L a_{i'_0-1,i'_0} \geq_L b_{i''_0-1,i''_0}$. So we have $c_{i',i'+1} \in I(b_{i''_0-1,i''_0},b_{j',j'+1})$. Hence $\delta=\beta[i',j',c_{i',i'+1}]$ is well-defined, and it is elementarily superior to $\beta$.

\mk

Following the same argument, the element $\delta'=\gamma[i,j,b_{i,i+1}]$ is well-defined. We will prove that $\delta=\delta'$. Let $i_0 \in [1,i]$ be minimal such that $u_{i_0}=u_{i}$. According to the proof of Lemma~\ref{lem:elementary_superior}, we can see that $\delta=\delta'=[d,x]$ are explicitly equal to the following:
\beq \forall k \not\in [i,j] \cup [i',j'], x_k&=&u_k \\
\forall k \in [i,j], x_k &=& u_{j}+1 \\
\forall k \in [i',j'], x_k &=& u_{j'}+1\\
\forall k \leq i_0-1, d_{k,k+1}&=&a_{k,k+1} \\
d_{i,i+1} &=& b_{i,i+1} \\
\forall k \in [j,i'_0-1], d_{k,k+1}&=&a_{k,k+1} \\
d_{i',i'+1} &=& c_{i',i'+1} \\
\forall k \geq j', d_{k,k+1}&=&a_{k,k+1}.\eeq

So the element $\delta=\delta'$ is elementarily superior to both $\beta$ and $\gamma$. 

\mk

We will now prove that $\delta$ is the minimal element of $M$ superior to both $\beta$ and $\gamma$. Fix $\theta=[e,y] \in M$ any element superior to both $\beta$ and $\gamma$, we will prove that $\delta \leq \theta$. Since $y \geq v$ and $y \geq w$, we deduce that $y \geq x$. Fix any $0 \leq k \leq n-1$ such that $x_k < x_{k+1}$, and let $\ell \in \{0,\dots,k\}$ denote the minimal element such that $y_{\ell+1} \geq x_{k+1}$. We will prove that $e_{\ell,\ell+1} \leq_L d_{k,k+1}$.

\bit
\item Assume that $k \leq j-1$. Then $x_{k+1}=v_{k+1}$, and since $\beta \leq \theta$, we deduce by Lemma~\ref{lem:two_orders_coincide} that $e_{\ell,\ell+1} \leq_L b_{k,k+1}=d_{k,k+1}$.
\item Assume that $k \geq j$. Then $x_{k+1}=w_{k+1}$, and since $\gamma \leq \theta$, we deduce by Lemma~\ref{lem:two_orders_coincide} that $e_{\ell,\ell+1} \leq_L c_{k,k+1}=d_{k,k+1}$.
\eit

According to Lemma~\ref{lem:two_orders_coincide}, we deduce that $\delta \leq_M \theta$: $\delta$ is the minimal element of $M$ superior to both $\beta$ and $\gamma$. Hence $\delta = \beta \vee_M \gamma$. Furthermore, we have noticed that $\delta$ is elementarily superior to both $\beta$ and $\gamma$, so $D(\beta,\delta)=D(\gamma,\delta)=1$.

\mk
 
{\bf Second case:} assume now that the intervals $[i,j]$ and $[i',j']$ intersect. Without loss of generality, assume that $i \leq i'$. Since $u_j < u_{j+1}$, we deduce that $j=j'$. Let $1 \leq i_0 \leq n$ be minimal such that $u_{i_0}=u_i$. The elements $b_{i-1,i}$ and $c_{i'-1,i'}$ both belong to the interval $I(a_{i_0-1,i_0},a_{j,j+1})$.

\mk

If $b_{i-1,i} = c_{i'-1,i'}$, then $\beta=\gamma$ and they have a trivial join $\delta=\beta=\gamma$. So we may assume that $b_{i-1,i} \neq c_{i'-1,i'}$.

\mk

If $b_{i-1,i} <_L c_{i'-1,i'}$, we have $\beta \leq \gamma$, so $\beta$ and $\gamma$ have a join $\delta=\gamma$ which satisfies $D(\beta,\delta)=1$ and $D(\gamma,\delta)=0$. Let us assume now that $b_{i-1,i} \not\leq_L c_{i'-1,i'}$.

\mk

Consider the meet $g= b_{i-1,i} \wedge_L c_{i'-1,i'} \in L$: its rank $r-1$ is such that $i_0 \leq r < i,i'$. Let us define $\delta=\alpha[r,j,g] \in M$. We see that $\delta = \beta[r,i-1,g] = \gamma[r,i'-1,g]$, so $\delta$ is elementarily superior to $\beta$ and $\gamma$.

\mk

We will now prove that $\delta=[d,x]$ is the minimal element of $M$ superior to both $\beta$ and $\gamma$. Fix $\theta=[e,y] \in M$ any element superior to both $\beta$ and $\gamma$, we will prove that $\delta \leq \theta$.
\mk

We will first prove that $x \leq y$. Since $\beta,\gamma \leq_M \theta$, we deduce that $v,w \leq y$. In particular, for any $m < r$ or $m \geq i$, we have $y_m \geq b_m=x_m$. And for $r \leq m \leq i-1$, we have $y_m \geq a_m=x_m-1$. Assume by contradiction that there exists $m \in \{r,\dots,i-1\}$ such that $y_m = x_m-1$, and choose such $m$ maximal. Since $x_r=x_{r+1}=\dots=x_j$, we have $y_{m+1} \geq x_{m+1} = x_i=v_i=x_{i'}=w_{i'}$. Since $\beta \leq_L \theta$ and $\gamma \leq_L \theta$, according to Lemma~\ref{lem:two_orders_coincide}, we know that $e_{m,m+1} \leq_L b_{i-1,i}$ and $e_{m,m+1} \leq_L c_{i'-1,i'}$. In particular, we deduce that $e_{m,m+1} \leq_L b_{i-1,i} \wedge_L c_{i'-1,i'} = g$. Note that the rank of $e_{m,m+1}$ is $m$, whereas the rank of $g$ is $r-1$. Since $m > r-1$, this is a contradiction. Hence $x \leq y$.

\mk

Fix any $0 \leq k \leq n-1$ such that $x_k < x_{k+1}$, and let $\ell \in \{0,\dots,k\}$ denote the minimal element such that $y_{\ell+1} \geq x_{k+1}$. We will prove that $e_{\ell,\ell+1} \leq_L d_{k,k+1}$.

\bit
\item Assume that $k \leq r-2$. Then $x_{k+1}=v_{k+1}$, and since $\beta \leq \theta$ we deduce by Lemma~\ref{lem:two_orders_coincide} that $e_{\ell,\ell+1} \leq_L b_{k,k+1}=d_{k,k+1}$.
\item Assume that $r-1 \leq k \leq j-1$. Then $x_{k+1}=v_i=w_{i'}$. Since $\beta \leq \theta$, we deduce by Lemma~\ref{lem:two_orders_coincide} that $e_{\ell,\ell+1} \leq_L b_{i-1,i}$. And since $\gamma \leq \theta$, we also deduce that $e_{\ell,\ell+1} \leq_L c_{i'-1,i'}$. Hence $e_{\ell,\ell+1} \leq_L b_{i-1,i} \wedge_L c_{i'-1,i'} = g = d_{r-1,r} \leq_L d_{k,k+1}$.
\item Assume that $k \geq j$. Then $x_{k+1}=v_{k+1}$, and since $\beta \leq \theta$ we deduce by Lemma~\ref{lem:two_orders_coincide} that $e_{\ell,\ell+1} \leq_L b_{k,k+1}=d_{k,k+1}$.
\eit

According to Lemma~\ref{lem:two_orders_coincide}, we deduce that $\delta \leq_M \theta$: $\delta$ is the minimal element of $M$ superior to both $\beta$ and $\gamma$. Hence $\delta = \beta \vee_M \gamma$. Furthermore, we have noticed that $\delta$ is elementarily superior to both $\beta$ and $\gamma$, so $D(\beta,\delta)=D(\gamma,\delta)=1$.
\ep
 
\blem \label{lem:join_distance_arbitrary}
Assume that $\alpha,\beta,\gamma \in M$ are such that $\alpha \leq_M \beta$, $\alpha \leq_M \gamma$, $D(\alpha,\beta)=m$ and $D(\alpha,\gamma)=m'$ for some $m,m' \in \N$. Then $\beta,\gamma$ have a meet $\delta$ such that $D(\beta,\delta) \leq m'$ and $D(\gamma,\delta) \leq m$.
\elem

\bp
We proceed by induction on $m+m'$: when $m+m' \leq 2$, the statement holds by Lemma~\ref{lem:join_distance_1}. Now fix $k \geq 3$, and assume that the statement holds when $m+m' <k$. Fix $m,m'$ such that $m+m'=k$, and without loss of generality assume that $m \geq 2$. Choose $\beta_0=\alpha,\beta_1,\dots,\beta_m=\beta$ an elementary sequence from $\alpha$ to $\beta$, with $m=D(\alpha,\beta)$. We have $D(\alpha,\beta_1)+D(\alpha,\gamma)=1+m'<k$, so by induction there exists $\delta'=\beta_1 \vee_M \gamma$ with $D(\beta_1,\delta') \leq m'$ and $D(\gamma,\delta') \leq 1$. Since $D(\beta_1,\beta)+D(\beta_1,\delta') \leq m-1+m'<k$, by induction there exists $\delta=\beta \vee_M \delta'$ such that $D(\beta,\delta) \leq m'$ and $D(\delta',\delta) \leq m-1$. So we deduce that $D(\gamma,\delta) \leq m$.

\mk

We will now prove that $\delta$ is the meet of $\beta$ and $\gamma$. We have $\delta \geq_M \beta$ and $\delta \geq_M \delta' \geq_M \gamma$. Furthermore, consider any $\theta \in M$ such that $\theta \geq_M \beta$ and $\theta \geq_M \gamma$. As $\beta \geq_M \beta_1$, we deduce that $\theta \geq_M \beta_1 \vee_M \gamma = \delta'$. And we deduce that $\theta \geq_M \beta \vee_M \delta'= \delta$. So we have proved that $\delta = \beta \vee_M \gamma$.
\ep

\bp[of Theorem~\ref{thm:affine_poset_lattice}]
Fix any $\beta=[b,v],\gamma=[c,w] \in M$. Let $k \geq 0$ such that $v_n-k < w_1$. Let $u=(v_1-k,v_2-k,\dots,v_n-k)$. Then $\alpha=[b,u] \in M$ is inferior to $\beta$, and we will see that it is also inferior to $\gamma$. Indeed let $\gamma'=[c,w']$, where $w'=(w_1,w_1,\dots,w_1)$. Since $w' \leq w$, we have $\gamma' \leq_M \gamma$. On the other hand, since $\gamma'=[b,w']$ and $u \leq w'$, we have $\alpha \leq_M \gamma'$, so $\alpha \leq_M \gamma$.

We can now apply Lemma~\ref{lem:join_distance_arbitrary} to deduce that $\beta$ and $\gamma$ have a meet in $M$. By symmetry of the construction, $\beta$ and $\gamma$ also have a join in $M$. So $M$ is a lattice.
\ep

If $H=\R$, the affine version $M_\R$ of $L$ over $\R$ is a gluing of subspaces $\sigma \subset \R^n$. We may therefore endow $M_\R$ with the piecewise length metric $d_\R$ induced by the standard $\ell^\infty$ metric on each $\sigma \subset \R^n$.

\mk

Let us define an action of $\R$ on $M_\R$ as follows:
\beq \R \times M_\R & \ra & M_\R \\
(t,[a,u]) & \mapsto & t \cdot [a,u]=[a,(u_1+t,u_2+t,\dots,u_n+t)].\eeq
This action is well-defined, preserves the order $\leq_M$, is increasing and continuous. Moreover, we have the following property.

\blem \label{lem:cofinal_action}
For any $\alpha,\beta \in M_\R$, there exists $t>0$ such that $(-t) \cdot \alpha \leq_M \beta \leq_M t \cdot \alpha$.
\elem

\bp
Consider representatives $\alpha=[a,u]$ and $\beta=[b,v]$ of $\alpha$ and $\beta$ respectively. Let $t>0$ such that $v_n \leq u_1+t$ and $u_n \leq v_1+t$. Then, if we denote $\gamma=[a,(v_n,\dots,v_n)]=[b,(v_n,\dots,v_n)]$, we have $\beta \leq_M \gamma \leq_M t \cdot \alpha$, hence $\beta \leq_M t \cdot \alpha$. Similarly, we have $(-t) \cdot \alpha \leq_M \beta$.
\ep

\bthm \label{thm:affine_lattice_R_injective}
If $H=\R$, the poset $M_\R$ is a lattice, and the metric space $(M_\R,d_\R)$ is injective.
\ethm

\bp
Remark that, for each $\theta>0$, the space $M_{\theta \Z}$ may be realized naturally as a closed subspace of $M_\R$. Furthermore, the sequence of closed subsets $M_{\theta \Z}$ of $M_\R$ converges to $M_\R$ as $\theta \ra 0$. We will use this convergence to prove that $M_\R$ is a lattice. We will then prove that the assumptions of Theorem~\ref{thm:lattice_Helly_injective} are satisfied.

\mk

We will now prove that any pair of elements in $M_\R$ have a join. Fix $\alpha,\beta \in M_\R$, we will define a common upper bound $\gamma$ for $\alpha$ and $\beta$.

Consider maximal chains $a,b$ in $L$ such that $\alpha \in \sigma_a$ and $\beta \in \sigma_b$. For any $\theta>0$, note that for any $x \in \R$, there exists $x_\theta \in \theta\Z$ such that $x_\theta-\theta \leq x \leq x_\theta+\theta$. So we may consider $\alpha_\theta \in M_{\theta \Z} \cap \sigma_a$ and $\beta_\theta \in M_{\theta \Z} \cap \sigma_b$ such that $(-\theta) \cdot \alpha_\theta \leq_M \alpha \leq_M \theta \cdot \alpha_\theta$, and similarly $(-\theta) \cdot \beta_\theta \leq_M \beta \leq_M \theta \cdot \beta_\theta$.

According to Theorem~\ref{thm:affine_poset_lattice}, the poset $M_{\theta \Z}$ is a lattice: consider $\gamma_\theta = \alpha_\theta \vee_{M_{\theta \Z}} \beta_\theta$. Let $C_{a,b} \subset L$ denote the smallest subset of $L$ containing $a$, $b$, and which is stable under meet. Since $L$ is bounded and graded, $C_{a,b}$ is finite. According to the proof of Theorem~\ref{thm:affine_poset_lattice}, we see that, for every $\theta>0$, there exists a maximal chain $c_\theta \subset C_{a,b}$ such that $\gamma_\theta \in \sigma_{c_\theta}$. Since $C_{a,b}$ is finite, $\sigma$ is locally compact and $(\gamma_\theta)_{\theta > 0}$ is bounded, there exists a sequence $\theta_k \ral{k \ra \pif} 0$  such that the sequence $(\gamma_{\theta_k})_{k \in \N}$ converges to some $\gamma \in M_\R$. Note that, for any $\theta>0$, we have $\gamma_\theta \geq_M \alpha_\theta \geq_M (-\theta) \cdot \alpha$. Since the sequence $(\gamma_{\theta_k})_{k \in \N}$ converges to $\gamma$, and the sequence $((-\theta_k) \cdot \alpha)_{k \in \N}$ converges to $\alpha$ by continuity of the action, we deduce that $\gamma \geq_M \alpha$. Similarly $\gamma \geq_M \beta$.

\mk

So $\gamma$ is a common upper bound for $\alpha$ and $\beta$. We will now prove that $\gamma$ is a minimal upper bound, which will prove that $\gamma$ is the join of $\alpha$ and $\beta$. Let us consider an upper bound $\delta \in M_\R$ of $\alpha$ and $\beta$, we will prove that $\gamma \leq_M \delta$. For any $\theta>0$, fix $\delta_\theta \in M_{\theta \Z}$ such that $(-\theta) \cdot \delta_\theta \leq_M \delta \leq_M \theta \cdot \delta_\theta$. In particular $(2\theta) \cdot \delta_\theta \geq_M \theta \cdot \delta \geq_M \theta \cdot \alpha \geq_M \alpha_\theta$, and similarly $(2\theta) \cdot \delta_\theta \geq_M \beta_\theta$. We deduce that $(2\theta) \cdot \delta_\theta \geq_M \alpha_\theta \vee_{M_{\theta \Z}} \beta_\theta = \gamma_\theta$. Considering the limit along $(\theta_k)_{k \in \N}$ as $k \ra \pif$, we deduce that $\delta \geq_M \gamma$.

\mk

So we have proved that $\alpha$ and $\beta$ have a join $\gamma$ in $M_\R$. By symmetry of the construction, they also have a meet, so $M_\R$ is a lattice.

\mk

We now turn to the assumptions of Theorem~\ref{thm:lattice_Helly_injective}: we will first prove that every upperly bounded subset of $M_\R$ has a join. Since $M_\R$ is a lattice, it is enough to prove that every bounded, increasing sequence is convergent. Fix an increasing sequence $(\alpha_k)_{k \geq 0}$ in $M_\R$, bounded above by some $\alpha \in M_\R$.

\mk

We will prove an intermediate result concerning $\ell^1$ metrics. Let us endow $M_\R$ with the length metric $d_1$ associated to the standard $\ell^1$ metric on each sector $\sigma_c \subset \R^n$. Let us also denote $d_1$ the standard metric on $\R^n$.
We claim that if $\beta=[b,v] \leq_M \gamma=[c,w]$, then $d_1(\beta,\gamma)=d_1(v,w)$. First notice that the second projection $M_\R \ra \R$ is $1$-Lipschitz with respect to the metrics $d_1$, hence we have $d_1(\beta,\gamma) \geq d_1(v,w)$. By definition of the order relation $\leq_M$, there exists a sequence $\beta_0=\beta \leq_M \beta_1 \leq_M \dots \leq_M \beta_p=\gamma$ such that, for each $0 \leq i \leq p-1$, the points $\beta_i,\beta_{i+1}$ lie in a common sector $\sigma_{c_i}$. Let us write representatives $\beta_i=[b_i,v_i]$, with $b_i \in L$ and $v_i \in \sigma$, for $0 \leq i \leq p$. We deduce that for each $0 \leq i \leq p-1$ we have $d_1(\beta_i,\beta_{i+1})=d_1(v_i,v_{i+1})$, hence
$$d_1(\beta,\gamma) \leq \sum_{i=0}^{p-1} d_1(\beta_i,\beta_{i+1}) \leq \sum_{i=0}^{p-1} d_1(v_i,v_{i+1}) = d_1(v,w).$$
So we have proved that if $\beta=[b,v] \leq_M \gamma=[c,w]$, then $d_1(\beta,\gamma)=d_1(v,w)$.

\mk

We now return to the increasing sequence $(\alpha_k)_{k \geq 0}$ in $M_\R$, bounded above by some $\alpha \in M_\R$. For each $k \in \N$, let us consider a representative $\alpha_k=[a_k,u_k]$ of $\alpha_k$, and a representative $\alpha=[a,u]$ of $\alpha$. Since the sequence $(\alpha_k)_{k \geq 0}$ is increasing in $M_\R$, we deduce that the sequence $(u_k)_{k \geq 0}$ is increasing in $\R^n$, and bounded above by $u$. Since increasing sequences in $\R^n$ are geodesics for the metric $d_1$, we deduce that for each $0 \leq j \leq k$, we have $d_1(u_j,u_k)+d_1(u_k,u)=d_1(u_j,u)$.

\mk

Then, for each $0 \leq j \leq k$, according to the claim about the metric $d_1$, we have $d_1(\alpha_j,\alpha_k)+d_1(\alpha_k,\alpha)=d_1(\alpha_j,\alpha)$. In particular, the sequence $(\alpha_k)_{k \geq 0}$ is a Cauchy sequence in $(M_\R,d_1)$. Since $M_\R$ has finitely many shapes, the proof of~\cite[Theorem~7.13]{bridson_haefliger} applies to show that the metric space $(M_\R,d_1)$ is complete, hence the sequence $(\alpha_k)_{k \geq 0}$ converges in $M_\R$. Equivalently, we can apply~\cite[Theorem~7.13]{bridson_haefliger} to the metric space $M_\R$ endowed with the length metric $d_2$ associated to the standard $\ell^2$ metric on each sector $\sigma_c \subset \R^n$. Since $d_1$ and $d_2$ are biLipschitz, this also implies that the metric space $(M_\R,d_1)$ is complete.

So every upperly bounded subset of $M_\R$ has a join.

\mk

We have proved that $M_\R$ is a lattice, such that each upperly bounded subset has a join. There is an increasing action of $\R$ on $M_\R$ satisfying the assumptions of Theorem~\ref{thm:lattice_Helly_injective} according to Lemma~\ref{lem:cofinal_action}. As a consequence, we deduce that the metric space $(M_\R,d)$ is injective, with respect to the metric
$$\forall x,y \in M_\R, d(x,y) =\inf\{t \geq 0 \st (-t) \cdot x \leq y \leq t \cdot x\}.$$
Note that the metric $d$ is geodesic, and it restricts on each $\sigma_c \subset M_\R$, for $c \in C(L)$, to the natural $\ell^\infty$ metric on $\sigma_c \subset \R^n$. Therefore $d$ coincides with the length metric $d_\R$.

\mk

So we conclude that $(M_\R,d_\R)$ is injective.
\ep

We have seen in the introduction that the existence of a bicombing may be extremely useful, notably in the case of Deligne complexes of Artin groups. Let us recall that a \emph{geodesic bicombing} on a metric space $X$ is a map $\sigma : X \times X \times [0,1] \ra X$ such that, for all $x,y \in X$, the map $t \in [0,1] \mapsto \sigma(x,y,t)$ is a constant speed geodesic from $x$ to $y$.

\mk

The bicombing $\sigma$ is called:
\bit
\item \emph{reversible} if $\forall x,y \in X, \forall t \in [0,1], \sigma(x,y,t) = \sigma(y,x,1-t)$,
\item \emph{consistent} if $\forall x,y \in X, \forall r,s,t \in [0,1], \sigma(\sigma(x,y,r),\sigma(x,y,s),t) = \sigma(x,y,(1-t)r+ts)$,
\item \emph{conical} if $\forall x,x',y,y' \in X, \forall t \in [0,1], d(\sigma(x,y,t),\sigma(x',y',t)) \leq (1-t)d(x,x')+td(y,y')$ and
\item \emph{convex} if $\forall x,x',y,y' \in X$, the map $t \in [0,1] \mapsto d(\sigma(x,y,t),\sigma(x',y',t))$ is convex.
\eit

Note that any consistent, conical bicombing is convex.

\bthm \label{thm:affine_lattice_R_bicombing}
The metric space $(M_\R,d_\R)$ has a unique convex, consistent, reversible geodesic bicombing.
\ethm

\bp
Let us denote $X=M_\R$ for simplicity. Given any $x,y \in X$, let us define
$$D(x,y) = \inf\{t \in \R \st x \leq t \cdot y\} \in \R.$$
Since the action of $\R$ on $M_\R$ is continuous, this infimum is attained, hence $x \leq D(x,y) \cdot y$. Note that this quantity is not symmetric with respect to $x$ and $y$, and we have
$$\forall x,y \in X, d(x,y) = \max(|D(x,y)|,|D(y,x)|).$$

\mk

We will start by defining a conical bicombing $\sigma$ on $X$ with a nice property which we call lower consistency. We then show that this is sufficient to bypass the use of properness of $X$ in~\cite[Theorem 1.4]{basso_bicombings} applied to $\sigma$.

\mk

Fix $x,y \in X$, $t \in [0,1]$, and let $D=d(x,y)$. For any $a \in X$, since $-D \cdot x \leq y \leq D \cdot x$, we have $|D(x,a)-D(y,a)| \leq D$, and so
\beq -tD \cdot x & \leq & (-tD(x,a)+tD(y,a)) \cdot x \\
& \leq & (-tD(x,a)+tD(y,a)) \cdot (D(x,a) \cdot a) \\
& \leq & ((1-t)D(x,a)+tD(y,a)) \cdot a. \eeq
Since every non-empty subset of $X$ with a lower bound has a meet, we may thus define
$$\sigma(x,y,t) = \bigwedge_{a \in X} ((1-t)D(x,a)+tD(y,a)) \cdot a.$$

\mk

We will first prove that it defines a geodesic bicombing, i.e. that $d(x,\sigma(x,y,t))=tD$. We have proved that, for every $a \in X$, we have $-tD \cdot x \leq ((1-t)D(x,a)+tD(y,a)) \cdot a$, hence $-tD \cdot x \leq \sigma(x,y,t)$. Conversely, when $x=a$, we know that
\beq \sigma(x,y,t) &\leq& ((1-t)D(x,x)+tD(y,x)) \cdot x \\
& \leq & tD(y,x) \cdot x \\
& \leq & tD \cdot x.\eeq
We conclude that $d(x,\sigma(x,y,t)) \leq tD$. By symmetry, we also have $d(\sigma(x,y,t),y) \leq (1-t)D$. Since $d(x,y)=D$, we conclude that $d(x,\sigma(x,y,t)) = tD$. So $\sigma$ is a geodesic bicombing. It is clear that $\sigma$ is reversible.

\mk

We will now prove that $\sigma$ is conical. Fix $x,y,z \in X$, and $t \in [0,1]$. For any $a \in X$, we have $|D(y,a)-D(z,a)| \leq d(y,z)$, hence $D(y,a) \leq D(z,a)+d(y,z)$. We deduce that
\beq \sigma(x,y,t) &=& \bigwedge_{a \in X} ((1-t)D(x,a)+tD(y,a)) \cdot a \\
& \leq & \bigwedge_{a \in X} ((1-t)D(x,a)+tD(z,a) +td(y,z)) \cdot a \\
& \leq & td(y,z) \cdot \sigma(x,z,t).\eeq
By symmetry, we also have $\sigma(x,z,t) \leq td(y,z) \cdot \sigma(x,y,t)$, hence $d(\sigma(x,y,t),\sigma(x,z,t)) \leq td(y,z)$. So the bicombing $\sigma$ is conical.

\mk

We will now prove that $\sigma$ is what we will call \emph{lower consistent}, which is one part of the inequality of the consistency equality. For each $x,y \in X$ and $s,t \in [0,1]$, we will prove that
$$\sigma(x,\sigma(x,y,t),s) \leq \sigma(x,y,st).$$
Let us denote $z=\sigma(x,y,t)$ and $w=\sigma(x,\sigma(x,y,t),s)$: we want to prove that $w \leq \sigma(x,y,st)$. For each $a \in X$, we have
\beq D(w,a) & \leq & (1-s)D(x,a)+sD(z,a) \\
& \leq & (1-s)D(x,a) + s((1-t)D(x,a)+tD(y,a) \\
& \leq & (1-st)D(x,a)+stD(y,a).\eeq
Hence we deduce that $w \leq D(w,a) \cdot a \leq (1-st)D(x,a) \cdot a+stD(y,a) \cdot a$. Since this holds for any $a \in X$, we conclude that $w \leq \sigma(x,y,st)$. Hence $\sigma$ is lower consistent.

\mk

According to~\cite[Lemma~5.2]{basso_bicombings}, given any $x,y \in X$ and $n \geq 1$, there exist unique elements $\sigma_{xy}(n,i)$, for $0 \leq i \leq n$, such that $\sigma_{xy}(n,0)=x$, $\sigma_{xy}(n,n)=y$ and
$$ \forall 1 \leq i \leq n-1, \sigma_{xy}(n,i) = \sigma\left(\sigma_{xy}(n,i-1),\sigma_{xy}(n,i+1),\f12\right).$$
Note that even though Lemma~5.2 is stated for a proper metric space, the uniqueness part only requires that $\sigma$ is a conical bicombing. And also the remark after the proof tells us that the only property needed for the existence is that the space is complete. We will actually give a proof below for the existence part.

\mk

Fix $n \geq 1$, we will prove that the elements $\sigma_{xy}(n,i)$ exist. For each $0 \leq i \leq n$, let us denote $x_i^0=\sigma(x,y,\f{i}{n})$. For each $k \in \N$, let us define $x_0^k=x$ and $x_n^k=y$. For each $k \in \N$ and $1 \leq i \leq n-1$, let us define inductively $x_i^{k+1}=\sigma(x_{i-1}^k,x_{i+1}^k,\f12)$. Since $\sigma$ is lower consistent, we see by induction that, for each $0 \leq i \leq n$, the sequence $(x_i^k)_{k \in \N}$ is nonincreasing in $X$. Moreover, since each $x_i^k$ lies on a geodesic from $x$ to $y$, we have $x_i^k \geq (-d(x,y)) \cdot x$ for every $k \in \N$ and $0 \leq i \leq n$.

\mk

Hence for each $0 \leq i \leq n$, since the sequence $(x_i^k)_{k \in \N}$ has a lower bound, we may define its meet $\sigma_{xy}(n,i)=\bigwedge_{k \geq 0} x_i^k$. In fact, the sequence $(x_i^k)_{k \geq 0}$ actually converges to $\sigma_{xy}(n,i)$. By continuity of $\sigma$, we deduce that the elements $(\sigma_{xy}(n,i))_{0 \leq i \leq n}$ satisfy the required property.

\mk

Remark that~\cite[Theorem 1.4]{basso_bicombings} is stated for a proper metric space, but the properness assumption is used in precisely two arguments: first in ~\cite[Lemma~5.2]{basso_bicombings} to prove the existence of the elements $\sigma_{xy}(n,i)$, which we obtained using specific properties of $X$. 

Properness of $X$ is used again, although not explicitly stated, in the proof of~\cite[Theorem 1.4]{basso_bicombings} to ensure the pointwise convergence of a sequence of bicombings with respect to some ultrafilter. Instead of using ultrafilters to ensure convergence, we will rather use the lower consistency of the bicombing.

\mk

For each $n \geq 1$, Basso states in~\cite[Lemma~5.2]{basso_bicombings} that the function $\sigma^{(n)}:X \times X \times [0,1] \ra X$ defined by
$$\sigma^{(n)}(x,y,(1-\lambda)\f{i}{n}+\lambda\f{i+1}{n})=\sigma(\sigma_{xy}(n,i),\sigma_{xy}(n,i+1),\lambda),$$
for all $x,y \in X$, $\lambda \in [0,1]$ and $0 \leq i \leq n-1$, is a conical bicombing.

\mk

First note that, since $\sigma$ is lower consistent, and by uniqueness of the points $\sigma_{xy}(n,i)$, we have for all $x,y \in X$, $n \geq 1$, $0 \leq i \leq n-1$ and $p \geq 1$ that $\sigma_{xy}(np,ip) \leq \sigma_{xy}(n,i)$.  For each $x,y \in X$ and $t \in [0,1]$, let us define
$$\gamma(x,y,t) = \bigwedge_{n \geq 1} \sigma_{xy}(n,\lceil tn\rceil).$$
According to the previous property, we deduce that
$$\gamma(x,y,t)= \lim_{n \ra +\infty} \sigma_{xy}(n!,\lceil tn!\rceil).$$
Since each $\sigma^{(n)}$ is a conical bicombing, one also deduces that $\gamma$ is conical.

\mk 

We will prove that $\gamma$ is lower consistent: let $x,y \in X$ and $s,t \in [0,1]$, we have
\beq \gamma(x,\gamma(x,y,t),s) &\leq& \lim_{n \ra +\infty} \gamma(x,\sigma_{xy}(n!,\lceil tn!\rceil),s) \\ 
&\leq& \lim_{n,m \ra +\infty} \sigma_{xy}(n!m!,\lceil s \lceil tn!\rceil m!\rceil) \\
&\leq& \lim_{n,m \ra +\infty} \sigma_{xy}(n!m!,\lceil stn!m!\rceil) \\
&\leq& \gamma(x,y,st),\eeq
as $d(\sigma_{xy}(n!m!,\lceil s \lceil tn!\rceil m!\rceil),\sigma_{xy}(n!m!,\lceil stn!m!\rceil)) \leq \f{\lceil sm!\rceil}{n!m!}d(x,y) \ra 0$ as $n \ra +\infty$. So we deduce that $\gamma$ is a reversible, conical, lower consistent geodesic bicombing such that $\gamma \leq \sigma$.

\mk

As a consequence, if we start with a reversible, conical, lower consistent geodesic bicombing $\sigma'$ which is minimal, we have $\gamma'=\sigma'$, hence $\sigma'$ satisfies the following consistency property: $\forall x,y \in X, \forall s,t \in [0,1], \sigma'(x,\sigma'(x,y,t),s) = \sigma'(x,y,st)$. Since $\sigma'$ is reversible, we deduce that $\sigma'$ is actually consistent. Since $\sigma'$ is conical, it is also convex.

According to~\cite[Theorem~1.2]{descombes_lang_hyperbolicity}, since $X$ has finite combinatorial dimension, we conclude that $\sigma'$ is the unique convex consistent reversible geodesic bicombing of $X$.
\ep

\bthm \label{thm:lattice_orthoscheme_injective}
Let $L$ denote a bounded, graded lattice. The orthoscheme realization $|L|$ of $L$, endowed with the piecewise $\ell^\infty$ metric, is injective. Moreover, $|L|$ has a unique convex reversible consistent geodesic bicombing.
\ethm

\bp
Consider the affine version $M=M_\R$ of $L$ over $\R$. For some maximal chain $c \in C(L)$, consider the elements $0_M=[(0,\dots,0),c] \in M$, $\mu_M = [(\f{1}{2},\dots,\f{1}{2}),c] \in M$ and $1_M=[(1,\dots,1),c] \in M$: note that $0_M$, $\mu_M$ and $1_M$ do not depend on $c$.

\mk

Note that the interval $I(0_M,1_M)$ coincides with the ball $B(\mu_M,\f{1}{2})$ in $M$ for the metric $d_\R$. According to Theorem~\ref{thm:affine_lattice_R_injective}, $(M,d_\R)$ is injective. So the ball $B(\mu_M,\f{1}{2})=I(0_M,1_M)$ is injective.

\mk

We remark that the interval $I(0_M,1_M)$ of $M$, endowed with the metric $d_\R$, is isometric to the orthoscheme realization of $L$, endowed with the piecewise $\ell^\infty$ metric. Indeed, for each maximal chain $c \in C(L)$, notice that $\sigma_C \cap I(0_M,1_M)$ identifies with the standard orthoscheme $\{x \in \R^n \st 0 \leq x_1 \leq \dots \leq x_n \leq 1\}$, where $n$ denotes the rank of $L$. We conclude that the orthoscheme realization of $L$ is injective.

\mk

According to Theorem~\ref{thm:affine_lattice_R_bicombing}, we know that $M_\R$ has a unique convex reversible consistent geodesic bicombing $\sigma$. Since the ball $|L| = B(\mu_M,\f{1}{2})$ is stable under $\sigma$, we deduce that $|L|$ has a convex reversible consistent geodesic bicombing. Since $|L|$ has combinatorial dimension at most $n$, according to~\cite[Theorem~1.2]{descombes_lang_hyperbolicity}, we deduce that $\sigma$ is the only convex bicombing on $|L|$.
\ep

Note that knowing when the orthoscheme complex of a lattice, endowed with the piecewise Euclidean metric, is CAT(0) is a very subtle question (see the introduction). We can prove the following.

\bthm \label{thm:orthoscheme_cat0_vs_injective}
Let $L$ denote a bounded, graded poset, let $|L|$ denote the geometric realization of $L$, and let $d_p$ denote the $\ell^p$ orthoscheme metric on $|L|$, for $p \in \{2,\infty\}$. If $(|L|,d_2)$ is CAT(0), then $(|L|,d_\infty)$ is injective. The converse is false.
\ethm

\bp
Assume that $L$ is a bounded, graded poset such that $(|L|,d_\infty)$ is not injective. According to Theorem~
\ref{thm:lattice_orthoscheme_injective}, we deduce that $L$ is not a lattice. According to Proposition~\ref{pro:bowtie}, there exists a bowtie in $L$: consider $x_1,x_2,x_3,x_4 \in L$ such that $x_1$ and $x_3$ are maximal elements inferior to both $x_2$ and $x_4$, and such that $x_2$ and $x_4$ are minimal elements superior to both $x_1$ and $x_3$. In the link $L_{[0,1]}$ of the diagonal edge $[0,1]$ in $(|L|,d_2)$, consider the piecewise geodesic loop $\ell$ going through $x_1,x_2,x_3,x_4,x_1$.

\mk

According to~\cite[Proposition~4.8]{brady_mccammond}, each geodesic segment $[x_i,x_{i+1}]$ has length smaller than $\f{\pi}{2}$, hence $\ell$ has length smaller than $2\pi$. On the other hand, according to~\cite[Lemma~7.2]{brady_mccammond} the loop $\ell$ is locally geodesic in $L_{[0,1]}$. Hence $L_{[0,1]}$ is not CAT(1), and $(|L|,d_2)$ is not CAT(0).

\mk

We will now prove that the converse does not hold. Consider the Coxeter symmetric group $W=\frak{S}_4$, with standard generators $S=\{s_1,s_2,s_3\}$ and standard Garside longest element $\Delta=s_1s_2s_3s_1s_2s_1$. Let $L$ denote the poset $W$, with order relation "being a left prefix for a shortest representative in $S$". Then $L$ is a lattice, and so $(|L|,d_\infty)$ is injective according to Theorem~\ref{thm:lattice_orthoscheme_injective}.

\mk

On the other hand, consider the piecewise geodesic loop $\ell$ in the link $L_{[0,1]}$ of the diagonal edge $[0,1]$ in $(|L|,d_2)$ going through the vertices $s_1$, $s_1s_2s_1$, $s_2$, $s_2s_3s_2$, $s_3$, $s_3s_1$, $s_1$. According to~\cite[Proposition~4.8]{brady_mccammond}, its length is
$$2 \arccos\left(\sqrt{\f{1}{2} \f{4}{5}}\right) + 4 \arccos \left( \sqrt{\f{1}{3} \f{3}{5}}\right) \simeq 0.987(2\pi) < 2\pi.$$
Since $L_{[0,1]}$ has the homotopy type of a circle, $\ell$ is not nullhomotopic in $L_{[0,1]}$. So $L_{[0,1]}$ is not CAT(1), and $(|L|,d_2)$ is not CAT(0).
\ep

We also deduce an immediate consequence concerning the Garside complex of a Garside group. Fix a Garside group $(G,S,\delta)$. Let $X$ denote the \emph{Garside complex} of $G$, i.e. the simplicial complex with vertex set $G$, and with simplices corresponding to chains $g_1 <_L g_2 <_L \dots <_L g_n$ such that $g_n \leq_L g_1 \delta$. Since simplices of $X$ have a total order on their vertices, $X$ is a simplicial complex with ordered simplices, so we may endow $X$ with the piecewise $\ell^\infty$ orthoscheme metric.

\bcor \label{cor:garside_complex_injective}
The Garside complex of any Garside group, endowed with the piecewise $\ell^\infty$ orthoscheme metric, is injective.
\ecor

This applies, in particular, to the dual braid complex studied by Brady and McCammond in~\cite{brady_mccammond}. This complex, endowed with the piecewise orthoscheme Euclidean metric, is conjectured to be CAT(0), but it is only known for a very small number of strands (see~\cite{brady_mccammond}, \cite{b6}, \cite{jeong}). On the other hand, if we endow the dual braid complex with the piecewise orthoscheme $\ell^\infty$ metric, we see that it is injective.

\section{Application to Euclidean buildings and the Deligne complex of type $\tilde{A_n}$} \label{sec:type_affine_A}

Let us consider the Euclidean Coxeter group $W \simeq \frak{S}_n \ltimes \Z^{n-1}$ of type $\tilde{A_{n-1}}$. Its Coxeter complex may be identified with
$$\Sigma = \{x \in \R^n \st x_1+\dots+x_n=0\}.$$
Up to homothety, we may choose the following affine hyperplanes to define $\Sigma$:
$$\left\{x_i-x_j =k \st 1 \leq i \neq j \leq n, k \in \Z\right\},$$
so that maximal simplices of $\Sigma$ identify with the $W$-orbit of
$$K = \{x \in \Sigma \st x_1 \leq x_2 \leq \dots \leq x_n \leq x_1+1\}.$$
The vertex set of $\Sigma$ identifies with $\left(\f{1}{n}\Z\right)^n \cap \Sigma$. Since the simplex $K$ is a strict fundamental domain for the action of $W$ on $\Sigma$, one may define a $W$-invariant type function $\tau$ on the vertex set of $\Sigma$:
\beq \tau : \Sigma^{(0)} &\ra& \Z/n\Z \\
x \in W \cdot v_i & \mapsto & i,\eeq
where $v_i=(\f{n-i}{n},\dots,\f{n-i}{n},\f{-i}{n},\dots,\f{-i}{n})$ is the vertex of $K$ whose first $i$ coordinates equal $\f{n-i}{n}$ and the $n-i$ last coordinates equal $\f{-i}{n}$.

This type function is such that adjacent vertices have distinct types.

\mk

Let us define the \emph{extended} Coxeter complex
$$\hat{\Sigma} = \Sigma \times \R = \R^n,$$
where the action of the standard generators $w_1,\dots,w_n$ of $W$ on $\hat{\Sigma}$ is given by
\beq \forall 1 \leq i \leq n-1, w_i \cdot (x_1,\dots,x_i,x_{i+1},\dots,x_n) &=& (x_1,\dots,x_{i+1},x_i,\dots,x_n) \\
\mbox{and } w_n \cdot (x_1,\dots,x_n) &=& (x_n-1,x_2,\dots,x_{n-1},x_1+1).\eeq
Also note that $\hat{\Sigma}$ has a natural simplicial complex structure, with vertex set $\Z^n$, with maximal simplices corresponding to the $W$-orbits of the simplices
$$\{k \leq x_i \leq x_{i+1} \leq x_{i+2} \leq \dots \leq x_{i+n-1} \leq k+1 \st k \in \Z, i \in \Z/n\Z\},$$
where indices in $\R^n$ are considered modulo $n$. Note that the action of the Coxeter group $W$ preserves $\tau$.

\mk

A fundamental domain for the action of $W$ on $\Sigma$ is the simplex
$$K = \{x \in \Sigma \st x_1 \leq x_2 \leq \dots \leq x_n \leq x_1+1\},$$
and a fundamental domain for the action of $W$ on $\hat{\Sigma}$ is the column
$$\hat{K} = \{x \in \hat{\Sigma} \st x_1 \leq x_2 \leq \dots \leq x_n \leq x_1+1\}.$$
Note that such columns have been studied by Brady-McCammond in~\cite{brady_mccammond}, and by Dougherty, McCammond and Witzel in~\cite{dougherty_mccammond_witzel}.

\mk

We will endow $\hat{\Sigma}$ with the standard $\ell^\infty$ metric from $\R^n$.

\mk

Let us consider a simplicial complex $X$ such that:
\bit
\item either $X$ is a Euclidean building of type $\tilde{A_{n-1}}$
\item or $X$ is the Deligne complex of the Euclidean Artin group $A$ of type $\tilde{A_{n-1}}$, with a coarser simplicial structure.
\eit
We will define the \emph{extended} version of $X$, denoted $\hat{X}$. It is a simplicial complex whose geometric realization is homeomorphic to $X \times \R$.

\mk

Before giving the precise definition of $\hat{X}$, let us first recall briefly the definition of the Deligne complex $\Delta$ of the Euclidean Artin group $A$ of type $\tilde{A_{n-1}}$ (see also Section~\ref{sec:intro_artin}). Let $S \simeq \Z/n\Z$ denote the standard generating set of $A$. Consider the set
$$\{gA_T \st g \in A, T \subsetneq S\},$$
endowed with the following partial order: $gA_T \leq g'A_{T'}$ if $gA_T \subset g'A_{T'}$. Then the Deligne complex $\Delta$ of $A$ is the geometric realization of this poset. We will define a coarser simplicial structure $X$ on $\Delta$. Note that, for any minimal vertex $gA_\emptyset \in \Delta$ (where $g \in A$), the $1$-neighbourhood of $gA_\emptyset$ in $\Delta$ is
$$\{gA_T \st T \subsetneq S\},$$
which is precisely the barycentric subdivision of the simplex with vertex set
$$\{gA_T \st T \subsetneq S, |S \bs T|=1\}$$
consisting of (the $g$-translates of) all maximal proper standard parabolic subgroups of $A$.

\mk

Therefore we may define the simplicial complex $X$ with vertex set $\{gA_T \st g \in A, T \subsetneq S, |S \bs T|=1\}$, and such that $g_1A_{T_1},\dots,g_kA_{T_k}$ span a simplex in $X$ if and only if $ \bigcap_{i=1}^k g_iA_{T_i} \neq \emptyset$.
Then $\Delta$ identifies with the barycentric subdivision of $X$, and also the geometric realizations of $\Delta$ and $X$ are homeomorphic.

\mk

Note that, in both cases (Euclidean building or Deligne complex), there is a well-defined type function $\tau : X^{(0)} \ra \Z/n\Z$ such that adjacent vertices have different types.
\bit
\item In the case $X$ is a Euclidean building of type $\tilde{A_{n-1}}$, each apartment is identified with the Coxeter complex $\Sigma$, we may define the type of a vertex $v$ of $X$ to be its type in any apartment containing $v$. Since the Coxeter group $W$ preserves the type, this definition does not depend on the choice of apartment.
\item In the case $X$ is the Deligne complex of type $\tilde{A_{n-1}}$, one may either use the projection onto the Coxeter complex, or use a direct definition: if $gA_T$ is a vertex of $X$, with $g \in A$ and $T =S \bs \{i\}$, its type is $i$.
\eit

\mk

The complex $X$ may also be defined as a particular gluing of copies of the Coxeter complex $\Sigma$. Roughly speaking, $\hat{X}$ will be the same gluing of copies of $\hat{\Sigma}$.

\mk

More precisely, let $\hat{X} = X \times \R$, and we will define a simplicial complex structure on $\hat{X}$. The vertex set of $\hat{X}$ is $\{(x,i) \in X^{(0)} \times \Z \st \tau(x)=i+n\Z\}$, where $X^{(0)}$ denotes the vertex set of $X$. Using the type function $\tau$, remark that vertices of maximal simplices of $X$ have a well-defined cyclic ordering in $\Z/n\Z$. The maximal simplices of $\hat{X}$ are
$$\{(x_i,kn+i),(x_{i+1},kn+i+1),\dots,(x_n,kn+n),(x_1,kn+n+1),(x_2,kn+n+2),\dots,(x_i,kn+n+i)\},$$
for any $k \in \Z$, any $1 \leq i \leq n$ and any maximal simplex $(x_1,x_2,\dots,x_n)$ of $X$ with $\forall 1 \leq i \leq n, \tau(x_i)=i$. We endow each such simplex with the $\ell^\infty$ orthoscheme metric for the given ordering. Endow $\hat{X}$ with the associated length metric.

\mk

Note that the translation action on the $\R$ factor of $\hat{X} = X \times \R$ defines an isometric action denoted $\theta$. Also note that $\theta$ does not preserve the simplicial structure of $\hat{X}$, but only its restriction to the subgroup $n\Z$ of $\R$.

\mk

If $X$ is a Euclidean building of type $\tilde{A_{n-1}}$, then $\hat{X}$ is called a Euclidean building of extended type $\tilde{A_{n-1}}$ (see~\cite{bruhat_tits}, and also~\cite{haettel_injective_buildings}).

\mk

If $X$ is the Deligne complex of the Euclidean Artin group $A$ of type $\tilde{A_{n-1}}$, we will give another description of $\hat{X}$. Recall that the classical Deligne complex $X$ may be defined as
$$X = (A \times K)/\sim,$$
where $(g,x) \sim (g',x')$ if $x=x'$ and, if the stabilizer of $x$ in $\Sigma$ equals $W_T$ for some $T \subsetneq S$, then $g^{-1}g' \in A_T$.
Then the extended Deligne complex may also be defined similarly as
$$\hat{X} = (A \times \hat{K})/\sim,$$
where $(g,x) \sim (g',x')$ if $x=x'$ and, if the stabilizer of $x$ in $\hat{\Sigma}$ equals $W_T$ for some $T \subsetneq S$, then $g^{-1}g' \in A_T$.

\mk

Fix any vertex $x \in X$, and let $L_{x,0}$ denote the set of vertices of $X$ adjacent to $x$. Without loss of generality, we may assume that $x$ has type $\tau(x)=0$. Note that vertices in $L_{x,0}$ have a type in $\Z/n\Z \bs \{\tau(x)\} = \{1,\dots,n\}$, which is an interval. Since maximal simplices of $X$ have a natural cyclic ordering in $\Z/n\Z$, there is a natural induced order on $L_{x,0}$ that is consistent with the type function $\tau : L_{x,0} \ra \{1,\dots,n\}$. Consider the poset $L_x=L_{x,0} \cup \{0,1\}$, where $0$ and $1$ are defined to be the minimum and the maximum of $L_x$ respectively.

\bpro \label{pro:locally_isometric_poset}
Consider any $p \in \hat{X}$ whose projection $p_X$ onto $X$ is contained in the open star of the vertex $x$. Then $\hat{X}$ is locally isometric at $p$ to a neighbourhood of a point in the $\ell^\infty$ orthoscheme realization of the poset $L_x$.
\epro

\bp
Remember that there is a diagonal, isometric action $\theta$ of $\R$ on $\hat{X}$, whose quotient is $X$. So, up to isometry, we may assume that $p$ lies in the open star of the diagonal edge joining the vertices $(x,0)$ and $(x,n)$ of $\hat{X}$. Note that each simplex of $\hat{X}$ containing $p$ corresponds to a chain in $L_x$ containing $(x,0)$ (identified with $0 \in L_x$) and $(x,n)$ (identified with $1 \in L_x$). As a consequence, a neighbourhood of $p$ in $\hat{X}$ is isometric to a neighbourhood of a point in the orthoscheme realization of $L_x$.
\ep

We now prove that, in the case of a Euclidean building, the poset $L_x$ is a lattice.

\bpro
Let $L_0$ denote the vertex set of a (possibly non-thick) spherical building of type $A_{n-1}$, for some $n \geq 1$, and let $L=L_0 \cup \{0,1\}$. In other words, $L$ is the poset of linear subspaces of a projective space of dimension $n-1$. Then $L$ is a bounded, graded lattice of rank $n-1$.
\epro

\bp
The lattice property is obvious: if $S,S'$ are linear subspaces of a projective space $X$, then $S \wedge S'= S \cap S'$ and $S \vee S'$ is the linear subspace spanned by $S \cup S'$.
\ep

As a consequence, if $X$ is a Euclidean building of type $\tilde{A_{n-1}}$, for any vertex $x \in X$, the poset $L_x$ is a bounded graded lattice.

\mk

We now turn to the case of the Deligne complex of type $\tilde{A_{n-1}}$. We defer the proof of the following result to Section~\ref{sec:cut_curves}.

\bthm[Crisp-McCammond]\ \label{thm:crisp_mccammond_lattice}

Let $X$ denote the Deligne complex of type $\tilde{A_{n-1}}$. For any vertex $x \in X$, the poset $L_x$ is a bounded, graded lattice.
\ethm

We may now deduce the main result.

\bthm \label{thm:AN_extended_injective}
Any Euclidean building of extended type $\tilde{A_{n-1}}$, or the extended Deligne complex of type $\tilde{A_{n-1}}$, endowed with the piecewise $\ell^\infty$ metric, is injective.
\ethm

\bp
According to Proposition~\ref{pro:locally_isometric_poset} and Theorem~\ref{thm:crisp_mccammond_lattice}, we know that $\hat{X}$ is locally isometric to the $\ell^\infty$ orthoscheme complex of a bounded graded lattice. According to Theorem~\ref{thm:lattice_orthoscheme_injective}, we deduce that $\hat{X}$ is uniformly locally injective.

If $X$ is a Euclidean building, $X$ and $\hat{X}$ are contractible, and in particular simply connected. If $X$ is a Deligne complex, according to Theorem~\ref{thm:deligne_equivalent_hyperplane_arrangement}, $\hat{X}$ is simply connected.

According to Theorem~\ref{thm:orthoscheme_complete_length}, $\hat{X}$ is complete. We deduce with Theorem~\ref{thm:cartan_hadamard_injective} that $\hat{X}$ is injective. 
\ep

We can also deduce that the thickening is a Helly graph. Note that this thickening, described more precisely in Theorem~\ref{thm:injective_orthoscheme_implies_helly}, corresponds to a coarser cell structure on the building or the Deligne complex.

\bcor
The thickening of the vertex set of any Euclidean building of extended type $\tilde{A_{n-1}}$, or the extended Deligne complex of type $\tilde{A_{n-1}}$, is a Helly graph.
\ecor

\bp
Let $X$ denote either a Euclidean building of extended type $\tilde{A_{n-1}}$ or the extended Deligne complex of type $\tilde{A_{n-1}}$. We will see that $X$ satisfies the assumptions of Theorem~\ref{thm:injective_orthoscheme_implies_helly}.

\mk

Since simplices of $X$ have a well-defined order, we see that $X$ is a simplicial complex with ordered simplices. For each maximal simplex $\sigma$ of $X$, the minimal and maximal vertices of $\sigma$ form a maximal edge in $X$. Therefore, $X$ has maximal edges as in Definition~\ref{def:ordered_simplices_maximal_edges}. According to Theorem~\ref{thm:injective_orthoscheme_implies_helly}, we deduce that the thickening of $X$ is a Helly graph.
\ep

We will now deduce a bicombing on $X$, considered as the quotient of $(\hat{X},d_{\hat{X}})$ by the diagonal action $\theta$ of $\R$ on $\hat{X}$. Let us define the quotient metric $d_X$ on $X$:
$$\forall x,y \in X, d_X(x,y) = \inf_{t \in \R} d_{\hat{X}}(x,\theta(t) \cdot y).$$

\bthm \label{thm:AN_bicombing}
Any Euclidean building $X$ of type $\tilde{A_{n-1}}$, or the Deligne complex $X$ of type $\tilde{A_{n-1}}$, has a metric $d_X$ that admits a convex, consistent, reversible geodesic bicombing $\sigma$. Moreover, $d_X$ and $\sigma$ are invariant under the group of type-preserving automorphisms of $X$.
\ethm

\bp
We will prove it locally, i.e. for the star $Y$ of a vertex $v$ of $X$. According to Theorem~\ref{thm:lattice_orthoscheme_injective}, there exists a unique convex, consistent, reversible bicombing $\hat{\sigma}$ on $\hat{Y}$. For each $x,y \in Y$, choose lifts $\hat{x},\hat{y} \in \hat{Y}$ such that $d_{\hat{X}}(\hat{x},\hat{y}) = d_X(x,y)$. For each $t \in [0,1]$, let us define $\sigma(x,y,t) = \theta(\R) \cdot \hat{\sigma}(\hat{x},\hat{y},t) \in Y$.

\mk

We will first see that $\sigma$ is well-defined: indeed fix $x,y \in Y$, and consider two pairs of lifts $\hat{x},\hat{x}',\hat{y},\hat{y'} \in \hat{Y}$ such that $d_{\hat{X}}(\hat{x},\hat{y}) = d_{\hat{X}}(\hat{x}',\hat{y}') = d_X(x,y)$. Note that, up to the action of $\theta$, we may assume that $\hat{x}=\hat{x}'$. If $\hat{y} \neq \hat{y}'$, let us denote $a \in \R \bs \{0\}$ such that $\hat{y}'=\theta(a) \cdot \hat{y}$. We then have $d_{\hat{X}}(\hat{x},\theta(\f{a}{2}) \cdot \hat{y}) < d(x,y)$, which is a contradiction. So $\hat{y} = \hat{y}'$, and $\sigma$ is well-defined.

Moreover, it is easy to see that $\sigma$ is a reversible, consistent, geodesic bicombing.

\mk

We will now see that $\sigma$ is convex. Consider $x,x',y,y' \in Y$, and choose lifts $\hat{x},\hat{x'},\hat{y},\hat{y}' \in \hat{Y}$ such that $d_{\hat{X}}(\hat{x},\hat{x}') = d_X(x,x')$ and $d_{\hat{X}}(\hat{y},\hat{y}') = d_X(y,y')$. Fix any $t \in [0,1]$ and $s \in \R$. We have
\beq d_X(\sigma(x,y,t),\sigma(x',y',t)) &=& \inf_{s \in \R} d_{\hat{X}}(\hat{\sigma}(\hat{x},\hat{y},t),\theta(s) \cdot \hat{\sigma}(\hat{x},\hat{y},t)) \\
& \leq & d_{\hat{X}}(\hat{\sigma}(\hat{x},\hat{y},t),\hat{\sigma}(\hat{x},\hat{y},t)) \\
& \leq & (1-t)d_{\hat{X}}(\hat{x},\hat{x}')+td_{\hat{X}}(\hat{y},\hat{y}') \\
& \leq & (1-t)d_X(x,x')+td_X(y,y'), \eeq
so $\sigma$ is a conical bicombing. Since $\sigma$ is also consistent, it is a convex bicombing.

\mk

We have seen that $(X,d_X)$ is locally convexly reversibly consistently bicombable. According to~\cite[Theorem~1.1]{miesch}, we deduce that $(X,d_X)$ has a unique global convex reversible geodesic consistent bicombing $\sigma$ that is consistent with local bicombing. Since the local bicombing only depends on the local combinatorics of $X$ and the type, we deduce that $\sigma$ is invariant under type-preserving automorphisms of $X$.
\ep

Note that we strengthened this result in~\cite{haettel_NPC_CUB}, by proving that this convex geodesic bicombing is actually unique. Our result also extends to simplicial complexes which are much more general than buildings.

\section{The lattice of maximal parabolic subgroups of braid groups} \label{sec:cut_curves}

This whole section is unpublished work of John Crisp and Jon McCammond, copied here with their permission. The results concerning the lattice of cut-curves are contained in Bessis's work (see~\cite{bessis_free}) about Garside structures on free groups. For consistency, we choose to follow Crisp and McCammond's presentation instead.

\subsection{The lattice of cut-curves}

\begin{figure}
\begin{center}
\begin{tikzpicture}
\def \p {0.05}
\def \op {1}
\def \gris {black!10}

\draw (0,0) circle (3);

\draw[fill] (0,3) circle (\p) node(a) {};
\node (a) at ([yshift=0.5cm]a) {\bfseries $a$};
\draw[fill] (0,-3) circle (\p) node(b) {};
\node (b) at ([yshift=-0.5cm]b) {\bfseries $b$};
\draw[fill] (-3,0) circle (\p) node(p-) {};
\node (p-) at ([xshift=-0.5cm]p-) {\bfseries $p^-$};
\draw[fill] (3,0) circle (\p) node(p+) {};
\node (p+) at ([xshift=0.5cm]p+) {\bfseries $p^+$};

\draw[fill] (-2,0) circle (\p) node(1) {};
\node (1) at ([yshift=-0.5cm]1) {\bfseries $p_1$};
\draw[fill] (-1,0) circle (\p) node(2) {};
\node (2) at ([yshift=-0.5cm]2) {\bfseries $p_2$};
\node (dots) at (0,0) {\bfseries $\dots$};
\draw[fill] (1,0) circle (\p) node(n') {};
\node (n') at ([yshift=-0.5cm]n') {\bfseries $p_{n-1}$};
\draw[fill] (2,0) circle (\p) node(n) {};
\node (n) at ([yshift=-0.5cm]n) {\bfseries $p_{n}$};

\draw [blue, thick] plot [smooth, tension=1] coordinates { (-3,0) (-1.5,2) (-0.5,0) (-1,-1) (-1.5,0.5) (-2.5,0) (0,-2)  (2,1) (3,0)};

\node (c) at (0,-2.5) {\color{blue} \bfseries $c$};

\end{tikzpicture}
\end{center}
\caption{An example of cut-curve $c$ in the punctured disk $\D_*$}
\label{fig:cut_curve}
\end{figure}
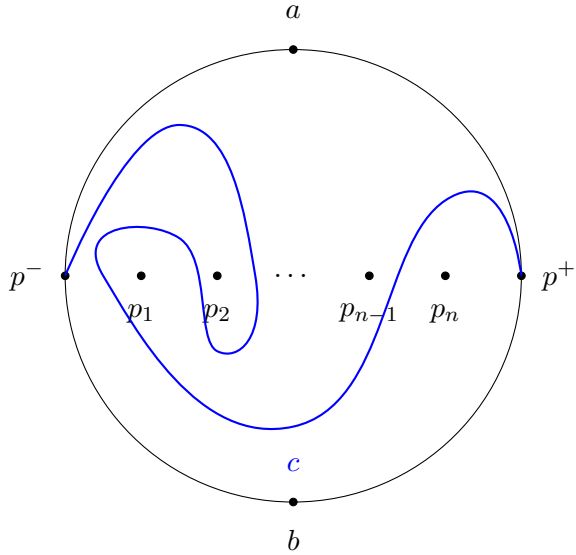

Let $\D^2$ denote the unit disk in $\R^2$, and let $\{p_1,\dots,p_n\}$ denote a set of $n$ distinct points in $\D^2$, as shown in Figure~\ref{fig:cut_curve}. We write $\D_*=\D^2 \bs \{p_1,\dots,p_n\}$. Let $a=(0,1)$ and $b=(0,-1)$ denote the top and bottom points of the boundary $\partial \D^2$.

\bdf
By a \emph{cut-curve} or \emph{curve} on $\D_*$ we shall mean a smoothly embedded curve $c$ in $\D_*$ which meets $\partial \D=\partial \D^2$ precisely at its endpoints, and which separates the boundary points $a$ and $b$.
\edf

We shall say that two cut-curves are isotopic if they are isotopic in $\D_*$ relative to $\{a,b\}$. We denote by $[c]$ the isotopy class of a curve $c$, and write $\cal C$ for the set of all isotopy classes of cut-curves in $\D_*$.

\mk

Observe that any cut-curve $c$ separates $\D_*$ into two regions, an upper region containing $a$ and a lower region containing $b$, and induces, in particular, a partition of the points $\{p_1,\dots,p_n\}$ into two sets. In general we shall say that the contents of the region containing $a$ lie \emph{above} $c$ and the contents of the region containing $b$ lie \emph{below} $c$. For each curve $c$, we write $\deg(c)$ for the number of points $p_i$ which lie below $c$. Clearly this number is invariant under isotopy, and so defines a degree function on $\cal C$ by $\deg([c])=\deg(c)$.

\mk

Let $c_1,c_2$ be two curves. We say that $c_1$ and $c_2$ are \emph{in minimal position} with respect to one another if they do not cobound any disk regions ("bigons") in $\D_*$. This includes triangular regions against the boundary. Any such disk regions can always be removed by modifying just one of the two curves in its isotopy class without changing the other. That is to say that, for any two curves $c_1,c_2$ we can always find $c'_1$ such that $[c'_1]=[c_1]$ and $c'_1$ is in minimal position with respect to $c_2$.

\mk

Here is another argument about minimal position using hyperbolic geometry. Such an argument has been used in~\cite{bessis_free}. If $n \geq 2$, one may also endow $\D_*$ with a fixed complete hyperbolic metric. Fix $p^+=(1,0),p^-=(-1,0) \in \partial \D$ points on each connected component of $\partial \D \bs \{a,b\}$. Then, for each isotopy class $[c]$, we may consider the unique geodesic line in $\D_*$ in the isotopy class $[c]$ with endpoints $p^+,p^- \in \partial \D$. Then, for any $[c_1],[c_2] \in {\cal C}$, the geodesic representatives are in minimal position.

\bdf
Let $c_1,c_2$ be two curves. We say that $[c_1] \leq [c_2]$ if $c_1$ is isotopic to a curve which lies below $c_2$. It is easily checked that this defines a partial order on the set $\cal C$ of cut-curves classes.
\edf

Note that if $c_1$ and $c_2$ are in minimal position with respect to one another then $[c_1] \leq [c_2]$ if and only if $c_1$ lies below $c_2$ (in particular they are disjoint). Note also that the function $\deg : {\cal C} \ra \{1,\dots,n\}$ is a strict order preserving map, a grading on the poset $({\cal C},\leq)$.

\bthm{\cite[Theorem~2.6]{bessis_free}}
The graded poset $({\cal C},\leq)$ is a lattice.
\ethm

\bp
Let $x=[c_1]$, $y=[c_2]$ be arbitrary elements of ${\cal C}$. Suppose that the curves $c_1$ and $c_2$ are in minimal position with respect to one another. Now consider how the union of $c_1$ and $c_2$ cut $\D_*$ up into connected regions. There is a unique lowermost such region $R$ which lies below both $c_1$ and $c_2$ (and contains the point $b$), and another uppermost such region $R'$ which lies above both $c_1$ and $c_2$ (and contains the point $a$). Let $c$, resp. $c'$, denote the curves which skirt along the part of the boundary of $R$, resp. $R'$, lying in the interior of $\D_*$ (i.e. not along $\partial \D$). Then we claim that $x \wedge y=[c]$ and $x \vee y = [c']$.

\mk

We check just the first of these two claims. Suppose that $c_0$ represents a common lower bound for $x$ and $y$. Then, by a sequence of bigon removing isotopies (or by considering the geodesic representative), we may choose the representative $c_0$ to be in minimal position with respect to both $c_1$ and $c_2$. Since $c_0$ represents a common lower bound, $c_0$ lies below both $c_1$ and $c_2$, hence lies below the curve $c$.
\ep

Recall that the $n$-strand braid group $B_n$ is isomorphic to the mapping class group $\MCG(\D_*,\partial \D)$. As a consequence, $B_n$ acts naturally on the set ${\cal C}$ of isotopy classes of cut-curves.

\blem
The action of $B_n$ on ${\cal C}$ preserves the order and the degree.
\elem

\bp
Let $[c] \in {\cal C}$ be a curve of degree $1 \leq k \leq n-1$ and $g \in B_n$. Then $k$ points among $\{p_1,\dots,p_n\}$ lie below $c$. Since the action of $B_n$ fixes the boundary of $\D_*$ pointwise, we know that $k$ points among $\{p_1,\dots,p_n\}$ lie below $g(c)$. Hence $\deg g([c])=\deg [c]$.

\mk

Similarly, consider two curves $[c_1],[c_2] \in {\cal C}$ such that $[c_1] \leq [c_2]$. Consider $c_1,c_2$ in minimal position, so that $c_1$ is below $c_2$. For any $g \in B_n$, we have that $g(c_1)$ is below $g(c_2)$, hence $g([c_1]) \leq g([c_2])$. We conclude that the action of $B_n$ preserves the order on ${\cal C}$.
\ep

Note that $B_n$ acts transitively on the set of cut-curves with fixed degree.

\subsection{Cosets in braid groups}

We will show how the cut-curve lattice $({\cal C},\leq)$ can be reinterpreted in purely algebraic terms.

\bdf
Recall that the braid group $B_n$ is generated by the set of standard generators $S=\{\sigma_1,\sigma_2,\dots,\sigma_{n-1}\}$. For each $1 \leq k \leq n-1$, let $P_k$ denote the maximal parabolic subgroup of $B_n$ generated by $S \bs \{\sigma_k\}$. Thus $P_k$ is isomorphic to a product $B_k \times B_{n-k}$. Note also that these subgroups are distinct for distinct values of $k$. We define the following augmented collection of cosets in $B_n$:
$${\cal B} = \{gP_k \st g \in B_n, 1 \leq k \leq n-1\} \cup \{0,1\}.$$
For $b \in {\cal B}$, we write $\deg(b)=k$ if $b=gP_k$ for $1 \leq k \leq n-1$, and $\deg(0)=0$, $\deg(1)=n$.
We define an order relation $\leq_{\cal B}$ on $\cal B$ as follows. First, define $0 \leq b \leq 1$ for all $b \in {\cal B}$. Otherwise, for $g_1,g_2 \in B_n$ and $1 \leq k_1,k_2 \leq n-1$, we write $g_1P_{k_1} \leq_{\cal B} g_2P_{k_2}$ if $k_1 \leq k_2$ and $g_1P_{k_1} \cap g_2P_{k_2} \neq \emptyset$.
\edf

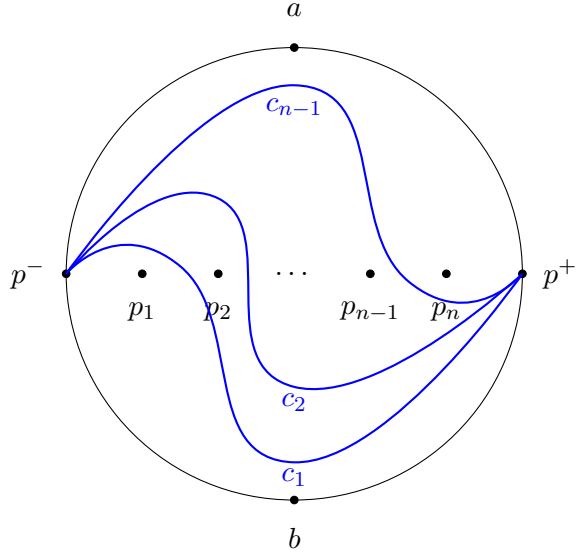
\begin{figure}
\begin{center}
\begin{tikzpicture}
\def \p {0.05}
\def \op {1}
\def \gris {black!10}

\draw (0,0) circle (3);

\draw[fill] (0,3) circle (\p) node(a) {};
\node (a) at ([yshift=0.5cm]a) {\bfseries $a$};
\draw[fill] (0,-3) circle (\p) node(b) {};
\node (b) at ([yshift=-0.5cm]b) {\bfseries $b$};
\draw[fill] (-3,0) circle (\p) node(p-) {};
\node (p-) at ([xshift=-0.5cm]p-) {\bfseries $p^-$};
\draw[fill] (3,0) circle (\p) node(p+) {};
\node (p+) at ([xshift=0.5cm]p+) {\bfseries $p^+$};

\draw[fill] (-2,0) circle (\p) node(1) {};
\node (1) at ([yshift=-0.5cm]1) {\bfseries $p_1$};
\draw[fill] (-1,0) circle (\p) node(2) {};
\node (2) at ([yshift=-0.5cm]2) {\bfseries $p_2$};
\node (dots) at (0,0) {\bfseries $\dots$};
\draw[fill] (1,0) circle (\p) node(n') {};
\node (n') at ([yshift=-0.5cm]n') {\bfseries $p_{n-1}$};
\draw[fill] (2,0) circle (\p) node(n) {};
\node (n) at ([yshift=-0.5cm]n) {\bfseries $p_{n}$};

\draw [blue, thick] plot [smooth, tension=1] coordinates { (-3,0) (-1.5,0.1) (0,-2.5) (3,0)};
\draw [blue, thick] plot [smooth, tension=1] coordinates { (-3,0) (-1,1) (0,-1.5) (3,0)};
\draw [blue, thick] plot [smooth, tension=1] coordinates { (-3,0) (0,2.5) (1.5,-0.1) (3,0)};

\node (c1) at (0,-2.7) {\color{blue} \bfseries $c_1$};
\node (c2) at (0,-1.7) {\color{blue} \bfseries $c_2$};
\node (cn) at (0,2.2) {\color{blue} \bfseries $c_{n-1}$};

\end{tikzpicture}
\end{center}
\caption{The base maximal chain $\alpha=(c_1,c_2,\dots,c_{n-1})$ of ${\cal C}$}
\label{fig:base_cut_curve}
\end{figure}

For simplicity, we shall henceforth identify $B_n$ with the mapping class group of $\D_*$ (relative to $\partial \D$). Let us consider the base maximal chain $\alpha=(c_1,c_2,\dots,c_{n-1})$ of ${\cal C}$ depicted in Figure~\ref{fig:base_cut_curve}.

\bdf
Define a map $\Phi:{\cal C} \ra {\cal B}$ by setting $\Phi(c)=0$ if $\deg(c)=0$, $\Phi(c)=1$ if $\deg(c)=n$, and otherwise
$$ \Phi(c) = \{g \in B_n \st g(\alpha) \mbox{ contains } c\}.$$
\edf

\blem \label{lem:image_of_phi}
If $c \in \alpha$ and $1 \leq \deg(c)=k \leq n-1$, then $\Phi(c)=P_k$.
\elem

\bp
The curve $c$ is a base curve such that $p_1,\dots,p_k$ lie below $c$ and $p_{k+1},\dots,p_n$ lie above $c$. Let $R^-$ denote the connected component of $\D^* \bs c$ below $c$, and $R^+$ denote the connected component of $\D^* \bs c$ above $c$. The stabilizer of $c$ in the mapping class group $B_n=\MCG(\D_*, \partial \D)$ is the direct product of $\MCG(R^-,\partial R^-)=\<\sigma_1,\dots,\sigma_{k-1}\>$ and $\MCG(R^+,\partial R^+)=\<\sigma_{k+1},\dots,\sigma_{n-1}\>$. Hence the stabilizer of $c$ in $B_n$ equals $P_k$.
\ep

\blem \label{lem:phi_surjective_degree}
The map $\Phi:{\cal C} \ra {\cal B}$ is well-defined, surjective, and respects degrees.
\elem

\bp
Fix $c \in {\cal C}$ with $1 \leq \deg(c)=k \leq n-1$. Fix some $g \in B_n$ such that the curve $c$ lies in $g(\alpha)$. Let $c_0=g^{-1}(c)$, the degree $k$ curve in $\alpha$. For any $h \in B_n$, we have
$$ h \in \Phi(c) \Leftrightarrow c \in h(\alpha) \Leftrightarrow g(c_0)=h(c_0).$$
According to Lemma~\ref{lem:image_of_phi}, this is equivalent to $gP_k=hP_k$. Hence $\Phi(c)=gP_k$ is well-defined. Moreover it is clear that $\Phi$ is surjective, and that for all $c \in {\cal C}$ we have $\deg(\Phi(c))=\deg(c)$.
\ep

\blem
The map $\Phi$ is injective, so it is a bijection between ${\cal C}$ and ${\cal B}$.
\elem

\bp
Let $c_1,c_2 \in {\cal C}$ such that $\Phi(c_1)=\Phi(c_2)$. We may suppose that $c_1$ lies in $g_1(\alpha)$ and $c_2$ in $g_2(\alpha)$, for some $g_1,g_2 \in B_n$. According to Lemma~\ref{lem:phi_surjective_degree}, we deduce that $c_1$ and $c_2$ have the same degree $0 \leq k \leq n$. Since ${\cal C}$ and ${\cal B}$ have unique elements of degree $0$ and $n$, we may restrict to the case $1 \leq k \leq n$.

\mk

Since $\Phi(c_1)=\Phi(c_2)$, we deduce that $g_1P_k=g_2P_k$, so ${g_1}^{-1}g_2 \in P_k$. Let $h \in P_k$ such that $g_2=g_1h$. Let $c_0$ denote the degree $k$ curve of $\alpha$. Then $c_1=g_1(c_0)$ and $c_2=g_2(c_0)$. Since $h \in P_k$, $h$ fixes $c_0$ so $c_1=c_2$.
\ep

\bthm[Crisp-McCammond] \label{thm:cut_curves_parabolic}
The bijection $\Phi : {\cal C} \ra {\cal B}$ is an order isomorphism between $({\cal C},\leq)$ and $({\cal B},\leq_{\cal B})$. As a consequence, $({\cal B},\leq_{\cal B})$ is a lattice.
\ethm

\bp
To prove the theorem, we need to show that for $c_1,c_2 \in {\cal C}$, we have $c_1 \leq c_2$ if and only if $\Phi(c_1) \leq_{\cal B} \Phi(c_2)$. Write $k_i = \deg(c_i)$ for $i=1,2$, and suppose, without loss of generality, that $0 \leq k_1 \leq k_2 \leq n$. Then $\Phi(c_1) \leq_{\cal B} \Phi(c_2)$ if and only if $\Phi(c_1) \cap \Phi(c_2) \neq \emptyset$ if and only if there exists a maximal chain $\alpha' \subset {\cal C}$ which contains both $c_1$ and $c_2$, if and only if $c_1 \leq c_2$ (since we already know that $\deg(c_1) \leq \deg(c_2)$).
\ep

We deduce the proof of Theorem~\ref{thm:crisp_mccammond_lattice} that each local poset $L_x$ in the Deligne complex of type $\tilde{A}_{n-1}$, is a lattice.

\bp[of Theorem~\ref{thm:crisp_mccammond_lattice}]
Let $x$ denote a vertex of the Deligne complex $X$ of type $\tilde{A}_{n-1}$. Without loss of generality, we may assume that $x$ corresponds to the maximal proper parabolic subgroup $A(A_{n-1})=\<\sigma_1,\dots,\sigma_{n-1}\>$, which is isomorphic to the Artin group of type $A_{n-1}$, i.e. the $n$-strand braid group. Since $A(A_{n-1})$ is a maximal parabolic subgroup of spherical type of the Euclidean Artin group $A(\tilde{A}_{n-1})$, the star of $x$ in $X$ identifies with the Deligne complex of the Artin group $A(A_{n-1})$. In particular, vertices in $X$ adjacent to $x$ may be identified with cosets of proper maximal parabolic subgroups of the braid group $A(A_{n-1})$. Furthermore, given two such cosets $gP_i$ and $hP_j$, for $g,h \in A(A_{n-1})$ and $1 \leq i,j \leq n-1$, the corresponding vertices of $X$ are adjacent if and only if $gP_i \cap hP_j \neq \emptyset$. As a consequence, the poset $L_{x,0}$ is isomorphic to the poset ${\cal B} \bs \{0,1\}$, and $L_x$ is isomorphic to the poset ${\cal B}$. According to Theorem~\ref{thm:cut_curves_parabolic}, we deduce that $L_x$ is a lattice.
\ep

\section{The thickening of a semilattice} \label{sec:semilattice}

We will now consider a generalization of Theorem~\ref{thm:lattice_orthoscheme_injective} to the case of a semilattice. This will be useful to consider Euclidean buildings, Deligne complexes in type $\tilde{C_n}$ and Artin groups of type FC. We start by recalling the definition of a flag poset.

\bdf
A poset $L$ is called \emph{flag} if any three elements which are pairwise upperly bounded have an upper bound.
\edf

\blem \label{lem:flag_implies_join_for_bounded}
Let $L$ denote a graded flag meet-semilattice with bounded rank. Then any family of elements of $L$ which are pairwise upperly bounded have a join.
\elem

\bp
We first prove that every finite family of $k$ elements $\{a_i\}_{1 \leq i \leq k}$ of $L$ which are pairwise upperly bounded have a join, by induction on the number of elements. By assumption, the property is true for $k=3$. Fix $k \geq 4$, assume that the property is true for $k-1$, and consider $k$ elements $\{a_i\}_{1 \leq i \leq k}$ of $L$ which are pairwise upperly bounded. Let $b$ denote the join of $a_{k-1}$ and $a_k$. The family $\{a_1,a_2,\dots,a_{k-2},b\}$ is pairwise upperly bounded, so by assumption it has a join $c \in L$. Then $c$ is the join of $\{a_i\}_{1 \leq i \leq k}$, which proves the induction.

\mk

Now consider an arbitrary family $A$ of elements of $L$ which are pairwise upperly bounded. Since $L$ has bounded rank, we may consider a finite subset $F \subset A$ such that the rank of the join $b$ of $F$ is maximal among all joins of finite subfamilies of $A$. For any $a \in A$, the join $b_a$ of $F \cup \{a\}$ satisfies $b \leq b_a$, and the rank of $b_a$ is at most the rank of $b$, hence $b_a=b$. We deduce that $b$ is the join of $A$.
\ep

\bthm \label{thm:semilattice_injective}
Let $L$ denote a graded poset with minimum $0$ and with bounded rank such that:
\bit
\item $L$ is a meet-semilattice.
\item $L$ is flag.
\eit
Then the $\ell^\infty$ orthoscheme realization of $|L|$ is injective.
\ethm

\bp
Let us consider $\ov{L}=L \cup \{1\}$: it is a bounded graded lattice. Let us denote by $\ov{M}=|\ov{L}|$ the geometric realization of $\ov{L}$, and $M=|L| \subset \ov{M}$ the geometric realization of $L$. Note that we consider $\ov{L}$ as a subset of $\ov{M}$, its vertex set. The orthoscheme realization of $M$ is endowed with the induced length metric as a subspace of the $\ell^\infty$ orthoscheme realization of $\ov{M}$.

\mk

If we consider the affine version $\ov{M}_\R$ of $\ov{L}$ over $\R$, then the elements $0_M=[(0,\dots,0),c] \in \ov{M}_\R$ and $1_M=[(1,\dots,1),c] \in \ov{M}_\R$ are such that $\ov{M}$ naturally identifies with the interval $I_{\ov{M}_\R}(0_M,1_M)$ in $\ov{M}_\R$, as in the proof of Theorem~\ref{thm:lattice_orthoscheme_injective}. In particular, $M$ and $\ov{M}$ are also posets.

\mk

According to Theorem~\ref{thm:affine_lattice_R_injective}, the space $\ov{M}_\R$, with its standard $\ell^\infty$ metric $\ov{d}$, is injective. According to Theorem~\ref{thm:lattice_orthoscheme_injective}, its interval $\ov{M}=I_{\ov{M}_\R}(0_M,1_M)$ is injective. Let us denote by $d$ the length metric of $M$.

\mk

Fix $\eps>0$, and consider the graph $\Gamma_\eps$ with vertex set $M$, and with an edge between $x,y \in M$ if $\ov{d}(x,y) \leq \eps$. We will prove that $\Gamma_\eps$ is Helly, by proving that it is clique-Helly and that its triangle complex is simply connected.

\mk

Let $\sigma \subset \Gamma_\eps$ denote a maximal clique: we claim that there exist $a,b \in \ov{M}$ such that $\sigma$ is the intersection of an interval $I_{\ov{M}}(a,b)$ in $\ov{M}$ with $M$. Since $\sigma$ is a subset of $\ov{M}_\R$ of diameter at most $\eps$, there exist $a_0,b_0 \in \ov{M}_\R$ such that $\sigma \subset I_{\ov{M}_\R}(a_0,b_0)$ and $\ov{d}(a_0,b_0) \leq \eps$. Since $\ov{M}$ is an interval, we deduce that there exist $a,b \in \ov{M} \cap I_{\ov{M}_\R}(a_0,b_0)$ such that $\sigma \subset I_{\ov{M}_\R}(a,b)=I_{\ov{M}}(a,b)$. In particular, we also have $\ov{d}(a,b) \leq \eps$. Since $\sigma$ is a maximal clique, we deduce that $\sigma = I_{\ov{M}}(a,b) \cap M$.

\mk

We will prove that $\Gamma_\eps$ is clique-Helly. Consider a family of pairwise intersecting maximal cliques $(\sigma_i)_{i \in I}$ in $\Gamma_\eps$. For each $i \in I$, according to the previous paragraph, there exists $a_i \leq b_i$ in $\ov{M}$ such that $\sigma_i=I_{\ov{M}}(a_i,b_i) \cap M$. For each $i \in I$, let $m_i \in L$ minimal such that $a_i \leq m_i$. For each $i \neq j$ in $I$, since $\sigma_i$ and $\sigma_j$ intersect, the intersection $I_{\ov{M}}(a_i,b_i)  \cap I_{\ov{M}}(a_j,b_j) \cap M$ is non-empty. In particular, the elements $m_i$ and $m_j$ have a common upper bound in $L$.

The family $(m_i)_{i \in I}$ is pairwise upperly bounded. Since $L$ is a graded flag semilattice with bounded rank, according to Lemma~\ref{lem:flag_implies_join_for_bounded}, this family has a join $m = \bigvee_{i \in I} m_i$ in $L$. In particular, the family $(\sigma_i \cap I_{\ov{M}}(0,m))_{i \in I}$ of intervals in the lattice $I_{\ov{M}}(0,m)=I_{M}(0,m)$ has a non-empty intersection. So the graph $\Gamma_\eps$ is clique-Helly.

\mk

We will now prove that the triangle complex of $\Gamma_\eps$ is simply connected. For each $t \in [0,1]$, consider the map $\pi_t : M \ra M$ sending each $x \in I_{\ov{M}}(0,m)$ to the point on the affine segment joining $0$ to $x$ at distance $td(0,x)$ from $0$. Note that $\pi_t$ is $1$-Lipschitz with respect to the distance $d$. Furthermore, if $t,t' \in [0,1]$ are such that $|t-t'| \leq \eps$, then for each $x \in M$ we have $d(\pi_t(x),\pi_{t'}(x)) \leq \eps$. As a consequence, any combinatorial loop $\gamma$ in $\Gamma_\eps$ may be homotoped in the triangle complex of $\Gamma_\eps$ to the loop $\pi_0(\gamma)$, which is the constant loop at $0$. So we have proved that the triangle complex of $\Gamma_\eps$ is simply connected.

\mk

According to Theorem~\ref{thm:clique_helly_global_helly}, we deduce that the graph $\Gamma_\eps$ is Helly. In particular, for each $\eps>0$, the metric space $M$ is $\eps$-coarsely injective. Since $M$ is complete according to Theorem~\ref{thm:orthoscheme_complete_length}, we deduce by Lemma~\ref{lem:epsilon_injective} that $M$ is injective.\ep

We will now give natural examples of such semilattices, which will be used in the sequel for Euclidean buildings and Deligne complexes of Euclidean type different from $\tilde{A_n}$.

\bpro \label{pro:semilattice_building_bn}
Let $L_0$ denote the vertex set of a (possibly non-thick) spherical building of type $B_n$, for some $n \geq 1$, and let $L=L_0 \cup \{0\}$. In other words, $L$ is the poset of subspaces of a polar space of projective dimension $n-1$. Then $L$ is a graded semi-lattice of rank $n$ with minimum $0$, such that any pairwise upperly bounded subset of $L$ has a join.
\epro

\bp
The semi-lattice property is obvious: if $S,S'$ are subspaces of a polar space $X$, then $S \wedge S'= S \cap S'$. If $A \subset L$ are pairwise upperly bounded in $L$, let $A_X=\bigcup_{S \in A} S \subset X$. For any $x,y \in A_X$, there exists a subspace of $X$ containing $x$ and $y$. According to~\cite[7.2.1]{tits_spherical_type}, the subset $A_X$ is contained in a subspace $S_0$ of $X$. Hence, for any $S \in A$, we have $S \subset S_0$. So the intersection of all subspaces of $X$ containing $\bigcup_{S \in A} S$ is the join of $A$.
\ep

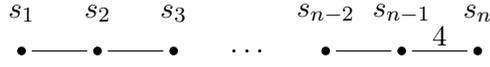
\begin{figure}
\begin{center}
\begin{tikzpicture}
\def \p {0.05}
\def \op {1}
\def \gris {black!10}

\draw[fill] (-3,0) circle (\p) node(s1) {};
\node (c1) at (-3,0.5) {\bfseries $s_1$};
\draw[fill] (-2,0) circle (\p) node(s2) {};
\node (c2) at (-2,0.5) {\bfseries $s_2$};
\draw[fill] (-1,0) circle (\p) node(s3) {};
\node (c3) at (-1,0.5) {\bfseries $s_3$};
\draw[fill] (1,0) circle (\p) node(s''n) {};
\node (c''n) at (1,0.5) {\bfseries $s_{n-2}$};
\draw[fill] (2,0) circle (\p) node(s'n) {};
\node (c'n) at (2,0.5) {\bfseries $s_{n-1}$};
\draw[fill] (3,0) circle (\p) node(sn) {};
\node (cn) at (3,0.5) {\bfseries $s_n$};

\node (dots) at (0,0) {\bfseries $\dots$};

\draw [-] (s1) edge (s2) (s2) edge (s3) (s''n) edge (s'n) (s'n) edge (sn);

\node (label) at (2.5,0.2) {$4$};

\end{tikzpicture}
\end{center}
\caption{The Dynkin diagram of type $B_n$}
\label{fig:dynkin_bn}
\end{figure}

Another application concerns the Artin group $A(B_n)$ of spherical type $B_n$, with the following Dynkin diagram (see Figure~\ref{fig:dynkin_bn}).

Let $s_1,\dots,s_n$ denote the standard generators of $A(B_n)$, with $s_{n-1}s_ns_{n-1}s_n=s_ns_{n-1}s_ns_{n-1}$. For each $1 \leq i \leq n$, let $P_i$ denote the maximal proper standard parabolic subgroup of $A(B_n)$
$$P_i = \<s_1,\dots,s_{i-1}\> \times \<s_{i+1},\dots,s_n\>.$$
Let $L_0=\{gP_i \st g \in A(B_n), 1 \leq i \leq n\}$ and $L=L_0 \cup \{0\}$. We define an order relation on $L$ as follows. First, define $0 \leq gP_i$ for each $gP_i \in L$. Otherwise, for $g_1,g_2 \in A(B_n)$ and $1 \leq k_1,k_2 \leq n$, define $g_1P_{k_1} \leq g_2P_{k_2}$ if $k_1 \leq k_2$ and $g_1P_{k_1} \cap g_2P_{k_2} \neq \emptyset$.

\blem \label{lem:order_relation_parabolics_type_A}
The relation $\leq$ is a partial order on $L$.
\elem

\bp
First note that $\leq$ is antisymmetric: if $g_1P_{k_1},g_2P_{k_2} \in L_0$ are such that $g_1P_{k_1} \leq g_2P_{k_2}$ and $g_2P_{k_2} \leq g_1P_{k_1}$, then $k_1=k_2$, and $g_1P_{k_1} \cap g_2P_{k_1} \neq \emptyset$ implies that $g_1P_{k_1} = g_2P_{k_1}$.

\mk

Now we show that $\leq$ is transitive: assume that $g_1P_{k_1},g_2P_{k_2},g_3P_{k_3} \in L_0$ are such that $g_1P_{k_1} \leq g_2P_{k_2}$ and $g_2P_{k_2} \leq g_3P_{k_3}$. We deduce that $k_1 \leq k_2 \leq k_3$, $g_1 \in g_2P_{k_2}P_{k_1}$ and $g_2 \in g_3P_{k_3}P_{k_2}$, so $g_1 \in g_3P_{k_3}P_{k_2}P_{k_1}$. Note that $\<s_1,\dots,s_{k_2-1}\> \subset P_{k_3}$ and $ \<s_{k_2+1},\dots,s_n\> \subset P_{k_1}$, so $P_{k_3}P_{k_2}P_{k_1}=P_{k_3}P_{k_1}$.

In particular $g_1 \in g_3P_{k_3}P_{k_1}$ and $k_1 \leq k_3$, so $g_1P_{k_1} \leq g_3P_{k_3}$.
\ep

Note that we are grateful to Luis Paris for his help in the following proof, notably the use of normal forms.

\bpro \label{pro:semilattice_artin_bn}
$L$ is a graded semi-lattice of rank $n$ with minimum $0$, such that any pairwise upperly bounded subset of $L$ has a join.
\epro

\bp
Let $t_1,\dots,t_{2n-1}$ the standard generators of the braid group $A(A_{2n-1})$. Consider the following morphism
\beq \phi : A(B_n) & \ra & A(A_{2n-1}) \\
\forall 1 \leq i \leq n-1, s_i & \mapsto & t_i t_{2n-i} \\
s_n & \mapsto & t_n.\eeq

This morphism $\phi$ is injective, see for instance~\cite{dehornoy_paris_gaussian}, \cite{michel_braid_monoid} and \cite{crisp_symmetrical}. Let $\sigma$ denote the involution of $A(A_{2n-1})$ defined by $\forall 1 \leq i \leq 2n-1, \sigma(t_i)=t_{2n-i}$. According to~\cite[Theorem~4]{crisp_symmetrical}, we know furthermore that the image of $\phi$ coincides with the fixed point set of $\sigma$.

\mk

For each $1 \leq i \leq 2n-1$, consider the standard proper maximal parabolic subgroup $Q_i=\<t_1,\dots,t_{i-1},t_{i+1},\dots,t_{2n-1}\>$ of $A(A_{2n-1})$. Let $M_0$ denote the poset of cosets of maximal proper parabolic subgroups of $A(A_{2n-1})$, and $M=M_0 \cup \{0,1\}$.  We will define a poset map $\psi:L_0 \ra M_0$. 

\mk

For each $g \in A(B_n)$ and each $1 \leq i \leq n$, let us define $\psi(gP_i)=\phi(g)Q_i$. Since $\phi(P_i) \subset Q_i$, the map $\psi : L_0 \ra M_0$ is well-defined. Assume that $gP_i \leq hP_j$, i.e. $i \leq j$ and $h^{-1}g \in P_jP_i$. Then $\phi(h)^{-1}\phi(g) \in \phi(P_j)\phi(P_i) \subset Q_jQ_i$, so $\phi(g)Q_i \leq \phi(h)Q_j$. As a consequence, $\psi$ is a rank-preserving injective poset map.

\mk

Note that the involution $\sigma$ extends naturally to an order-reversing involution on $M$, by letting $\sigma(gQ_i) = \sigma(g)Q_{2n-i}$. And for each $g \in A(B_n)$ and each $1 \leq i \leq n$, we have $\sigma(\psi(gP_i)) = \phi(g)Q_{2n-i}$.

\mk

We will now prove that $L$ is a meet-semilattice. Fix $a,b \in A(B_n)$ and $1 \leq i,j \leq n$, we will prove that $aP_i$ and $bP_j$ have a meet in $L$. Consider the elements $\phi(a)Q_i$ and $\phi(b)Q_j$ in $M$: they have a meet $\gamma Q_k$, for some $\gamma \in A(A_{2n-1})$ and $0 \leq k \leq i,j$. Then as $\phi(a)Q_i \leq \phi(a) Q_{2n-i}$ and $\phi(b)Q_j \leq \phi(b)Q_{2n-j}$, we deduce that $\gamma Q_k \leq \sigma(\gamma)Q_{2n-k}$.

So, up to the choice of $\gamma \in \gamma Q_k$, we may assume that $\sigma(\gamma) \in \gamma Q_{2n-k}$: let $q \in Q_{2n-k}$ such that $\sigma(\gamma) = \gamma q$. Since $\sigma$ is an involution, we have $\gamma = \gamma q \sigma(q)$, so $q=\sigma(q)^{-1} \in Q_k \cap Q_{2n-k}$. 

\mk

We claim that if $g \in A(A_{2n-1})$ is such that $\sigma(g)=g^{-1}$, there exists $h \in A(A_{2n-1})^+$ such that $g=h\sigma(h)^{-1}$. According to~\cite{charney_automation}, there exist unique $h,h' \in A(A_{2n-1})^+$ such that $g=hh'^{-1}$ and the right greatest common divisor of $h$ and $h'$ in the Garside monoid $A(A_{2n-1})^+$ is $1$. Note that $\sigma$ preserves the monoid $A(A_{2n-1})^+$, hence we have $g^{-1}=h'h^{-1}$ on one side, and $g^{-1}=\sigma(g)=\sigma(h) \sigma(h')^{-1}$ on the other side. By uniqueness of $h,h'$, we deduce that $h'=\sigma(h)$. Hence $g=h\sigma(h)^{-1}$.

Assume furthermore that $g \in Q_k \cap Q_{2n-k} = \<t_1,\dots,t_{k-1}\> \times \<t_k+1,\dots,t_{2n-k-1}\> \times \<t_{2n-k+1},\dots,t_{2n-1}\>$. Then we can decompose $g$ as a fraction $g=hh'^{-1}$ inside the parabolic subgroup $Q_k \cap Q_{2n-k}$: by uniqueness of $h,h'$, we deduce that $h,h' \in Q_k \cap Q_{2n-k}$.

\mk

According to the claim, there exists $q' \in Q_k \cap Q_{2n-k}$ such that $q=q' \sigma(q')^{-1}$. So, up to replacing $\gamma$ with $\gamma'=\gamma q' \in \gamma Q_k$, we have
$$\sigma(\gamma')=\sigma(\gamma)\sigma(q') = \gamma q \sigma(q') = \gamma q' = \gamma'.$$
So we may assume that $\gamma'$ is fixed by $\sigma$: according to~\cite[Theorem~4]{crisp_symmetrical}, we may consider $c \in A(B_n)$ such that $\phi(c)=\gamma'$. So we deduce that $cP_k \leq aP_i,bP_j$. Conversely, for any $c' \in A(B_n)$ and $k' \leq i,j$ such that 
$c'P_{k'} \leq aP_i,bP_j$, we have $\phi(c')Q_{k'} \leq \phi(a)Q_i \wedge \phi(b)Q_j = \phi(c)Q_k$, so $c'P_{k'} \leq cP_k$. We conclude that $cP_k$ is the meet of $aP_i$ and $bP_j$ in $L$.

\mk

We will now prove that any three pairwise upperly bounded elements of $L$ have an upper bound. Fix $g_1,g_2,g_3 \in A(B_n)$ and $1 \leq k_1,k_2,k_3 \leq n$ such that, for each $i \neq j$, the elements $g_iP_{k_i}$ and $g_jP_{k_j}$ have a join $g_{ij}P_{k_{ij}}$ in $L$. Consider the elements $\phi(g_{ij})Q_{k_{ij}}$ of $M$: they have a join $\gamma Q_k$ in $M$. Since, for each $1 \leq i,j \leq 3$, we have $\phi(g_i) Q_i \leq \phi(g_{ij}) Q_{k_{ij}} \leq \phi(g_{ij}) Q_{2n-k_{ij}} \leq \phi(g_j) Q_{2n-k_j}$, we deduce that $\gamma Q_k$ is inferior to the meet $\sigma(\gamma) Q_{2n-k}$ of the elements $\phi(g_{ij})Q_{2n-k_{ij}}$. We deduce that $k \leq n$. Also, as in the previous paragraph, we deduce that we can choose $\gamma \in A(A_{2n-1})$ such that $\gamma=\phi(c)$, for some $c \in A(B_n)$. Hence $cP_k$ is a common upper bound for the elements $(g_iP_{k_i})_{1 \leq i \leq 3}$.

\mk

Since $L$ is graded with finite rank, we deduce that any family of pairwise upperly bounded elements of $L$ have a join.
\ep

We believe that a similar statement holds for the Artin group of spherical type $D_n$, but we have no embedding into an Artin group of type $A$ in order to use a similar proof.

\section{Application to Euclidean buildings and the Deligne complex of other Euclidean types} \label{sec:type_affine_C}

In this section, we will prove that we can deduce from Theorem~\ref{thm:AN_extended_injective} an injective metric on Euclidean buildings and Deligne complexes of Euclidean types other than $\tilde{A_n}$. 

\mk

More precisely, let us consider a simplicial complex $X$ that is:
\bit
\item either a Euclidean building of type $\tilde{B_n}$, $\tilde{C_n}$ or $\tilde{D_n}$,
\item or the Deligne complex of Euclidean type $\tilde{C_n}$.
\eit

Note that the Coxeter groups of types $\tilde{B_n}$ and $\tilde{D_n}$ may be considered as subgroups of the Coxeter group of Euclidean type $\tilde{C_n}$, spanned by reflections with respect to subarrangements of hyperplanes, see below for more details. As a consequence, if $X$ is a Euclidean building of type $\tilde{B_n}$ or $\tilde{D_n}$, we will consider it as a (possibly non-thick) Euclidean building of type $\tilde{C_n}$.

\mk

Let $W(\tilde{C_n})$ denote the Euclidean Coxeter group of type $\tilde{C_n}$. Its Coxeter complex $\Sigma(\tilde{C_n})$ identifies with $\R^n$, and reflections of $W(\tilde{C_n})$ correspond to reflections with respect to hyperplanes
$$\{x_i = k \st 1 \leq i \leq n, k \in \Z\} \mbox{ and } \{x_i \pm x_j = 2k \st 1 \leq i \neq j \leq n, k \in \Z\}.$$

\mk

Note that we can see the following subarrangement has type $\tilde{B_n}$:
$$\{x_i = 2k \st 1 \leq i \leq n, k \in \Z\} \mbox{ and } \{x_i \pm x_j = 2k \st 1 \leq i \neq j \leq n, k \in \Z\}.$$
Also the following subarrangement has type $\tilde{D_n}$:
$$\{x_i \pm x_j = 2k \st 1 \leq i \neq j \leq n, k \in \Z\}.$$
This justifies that, in the case of Euclidean buildings, we may assume that $X$ is a (possibly non-thick) Euclidean building of type $\tilde{C_n}$.

\mk

The vertex set of $\Sigma(\tilde{C_n})$ identifies with $\Z^n$, and a strict fundamental domain of the action of $W(\tilde{C_n})$ on $\Sigma(\tilde{C_n})$ is given by the standard orthoscheme simplex $\sigma$ with vertices
$$v_0=(0,\dots,0), v_1=(1,0,\dots,0), \dots, v_n=(1,1,\dots,1).$$

\mk

Hence $v_0,v_1,\dots,v_n$ are representatives of the $n+1$ orbits of vertices of $\Sigma(\tilde{C_n})$ under the action of $W(\tilde{C_n})$. This enables to define a type function $\tau$ on vertices of $\Sigma(\tilde{C_n})$:
\beq \tau : \Sigma(\tilde{C_n})^{(0)} & \ra & \{0,1,\dots,n\} \\
g \cdot v_i & \mapsto & i.\eeq
More generally, we can define a partial order on vertices of $\Sigma(\tilde{C_n})$: say that $v < v'$ if $v,v'$ are adjacent vertices and $\tau(v) < \tau(v')$. We can therefore also view $\Sigma(\tilde{C_n})$ as the geometric realization of the poset of its vertices.

\mk

We will now see how to extend the type function and the partial order on vertices of $X$.

\bpro
Assume that $X$ is a Euclidean building of type $\tilde{C_n}$ or the Deligne complex of Euclidean type $\tilde{C_n}$. There exists a type function $\tau : X^{(0)} \ra \{0,1,\dots,n\}$ such that adjacent vertices of $X$ have different types. Moreover, let us define a partial order on vertices of $X$ by setting $v < v'$ if $v$ and $v'$ are adjacent in $X$ and $\tau(v) < \tau(v')$. Then $X$ is the geometric realization of the poset of its vertices.
\epro

\bp
Let us first consider the case where $X$ is a Euclidean building of type $\tilde{C_n}$. Given any vertex $v$ of $X$, consider any apartment $A \subset X$ containing $v$, and let us define $\tau(v)$ as defined with respect to the apartment $A \simeq \Sigma(\tilde{C_n})$. Since two apartments containing $v$ differ by an element of the Weyl group $W(\tilde{C_n})$, which preserves the type, we deduce that $\tau$ is well-defined on $X^{(0)}$. Similarly, given any two adjacent vertices $v,v'$ in $X$, say that $v < v'$ if $\tau(v) < \tau(v')$.

We will see that this relation is actually transitive on vertices of $X$: assume that three vertices $v_1,v_2,v_3$ of $X$ satisfy $v_1 < v_2$ and $v_2 < v_3$. Then the link of $v_2$ is isomorphic to the join of two spherical buildings of types $B_{\tau(v_2)}$ and $B_{n-\tau(v_2)}$. Hence we see that $v_1$ is adjacent to $v_3$ in $X$, so $v_1 < v_3$.

\mk

Let us now consider the case where $X$ is the Deligne complex of Euclidean type $\tilde{C_n}$. Let $s_0,s_1,\dots,s_{n-1},s_n$ denote the standard generators of the Artin group $A(\tilde{C_n})$. For each $0 \leq i \leq n$, consider the maximal spherical type standard parabolic subgroup
$$P_i = \<s_0,s_1,\dots,s_{i-1},s_{i+1},\dots,s_n\> \simeq A(B_i) \times A(B_{n-i}).$$
For each vertex $gP_i$ in $X$, where $g \in A(\tilde{C_n})$ and $0 \leq i \leq n$, let us define $\tau(gP_i)=i$. Given any two vertices $gP_i, hP_j$ of $X$, say that $gP_i < hP_j$ if they are adjacent in $X$ (i.e. $gP_i \cap hP_j \neq \emptyset$) and $i < j$. As in Lemma~\ref{lem:order_relation_parabolics_type_A}, one checks that it is a partial order.
\ep

We will endow $X$ with the piecewise $\ell^\infty$ orthoscheme metric $d$ given by the geometric realization of its poset of vertices. We will now describe the local structure of $X$ at any vertex, starting with a general statement about orthoscheme complexes of posets.

\blem \label{lem:local_splitting}
Let $L$ denote a graded poset, with an element $v \in L$ comparable to every element of $L$. Let $L^+=\{w \in L \st w \geq v\}$ and $L^-=\{w \in L \st w \leq v\}$. Then the $\ell^\infty$ orthoscheme realization of $L$ is locally isometric at $v$ to the $\ell^\infty$ product of the $\ell^\infty$ orthoscheme realizations of $L^+$ and $L^-$.
\elem

\bp
Since the $\ell^\infty$ orthoscheme realization of $L$ is obtained as a union of orthoschemes, it is sufficient to prove the result when $L$ is a chain.

In other words, consider the standard $\ell^\infty$ $n$-orthoscheme $C_n=\{x \in \R^n \st 1 \leq x_1 \leq x_2 \geq \dots \geq x_n \geq 0\}$, with vertices $v_i=(1,\dots,1,0,\dots,0)$ with $i$ ones and $n-i$ zeros, for $0 \leq i \leq n$. We fix a particular vertex $v=v_j$, for some $0 \leq j \leq n$, and we want to describe locally $C$ around $v$. The orthoscheme $C_n$ is locally isometric at $v$ to the space
$$T=\{x \in \R^n \st 1 \geq x_1 \geq x_2 \geq \dots \geq x_j \mbox{ and } x_{j+1} \geq x_{j+2} \geq \dots \geq x_n \geq 0\},$$
which is isometric to the $\ell^\infty$ product of $T^-=\{(x_1,\dots,x_j) \in \R^j \st 1 \geq x_1 \geq x_2 \geq \dots \geq x_j\}$ and $T^+=\{(x_{j+1},\dots,x_n) \in \R^{n-j} \st x_{j+1} \geq x_{j+2} \geq \dots \geq x_n \geq 0\}$.

\mk

The space $T^+$ is locally isometric at $v$ to the vertex $(1,1,\dots,1)$ in the $j$-orthoscheme $C_j$, and $T^-$ is locally isometric at $v$ to the vertex $(0,0,\dots,0)$ in the $(n-j)$-orthoscheme $C_{n-j}$.

In conclusion, if $L$ is a chain $v_0<v_1< \dots < v=v_j < v_n$, then the geometric realization of $|L|$ is locally isometric at $v$ to the product of the geometric realizations of the chains $L^-=(v_0<v_1< \dots < v_j=v)$ and $L^+=(v=v_j<v_{j+1}< \dots < v_n)$.
\ep

\bpro \label{pro:local_splitting}
Assume that $X$ is a Euclidean building of type $\tilde{C_n}$ or the Deligne complex of Euclidean type $\tilde{C_n}$. Fix a vertex $v \in X$ of type $\tau(v)=i \in \{0,1,\dots,n\}$.

There exist posets $L_{1,i}$ (resp. $L'_{0,n-i}$), which are either the poset of vertices of a spherical building of type $B_i$ (resp. $B_{n-i}$) or the poset of maximal proper parabolic subgroups of the Artin group of spherical type $B_i$ (resp. $B_{n-i}$). Let us define the posets $L_i=\{1\} \cup L_{1,i}$ and $L'_{n-i} = L'_{0,n-i} \cup \{0\}$, and let us consider the geometric realizations $|L_i|$ and $|L'_{n-i}|$, endowed with the piecewise $\ell^\infty$ orthoscheme metrics.

Then $X$ is locally isometric at $v$ to the $\ell^\infty$ direct product $|L_i| \times |L'_{n-i}|$.
\epro

\bp
The space $X$ is locally isometric at $v$ to the $\ell^\infty$ orthoscheme realization of the poset $L_v$ of vertices of $X$ comparable to $v$. The poset $L_v$ is the disjoint union of $L_{1,i} = \{w \in L_v \st w < v\} \sqcup \{v\} \sqcup L'_{0,n-i} = \{w \in L_v \st w > v\}$, such that $L_{1,i} < \{v\} < L'_{0,n-i}$. We may identify the poset $L_i=L_{1,i} \cup \{1\}$ with the interval $L_{1,i} \cup \{v\}$ in $L_v$. Similarly, we may identify the poset $L'_{n-i}=L'_{0,n-i} \cup \{0\}$ with the interval $L'_{0,n-i} \cup \{v\}$ in $L_v$. According to Lemma~\ref{lem:local_splitting}, the orthoscheme realization $|L_v|$ of $L_v$ is locally isometric at $v$ to the direct $\ell^\infty$ product $|L_i| \times |L'_{n-i}|$ of the orthoscheme realizations of $L_i$ and $L'_{n-i}$.

\mk

In case $X$ is a Euclidean building to type $\tilde{C_n}$, note that $L_{1,i}$ is the vertex set of a spherical building of type $B_i$ (with reversed order), and that $L'_{0,n-i}$ is the vertex set of a spherical building of type $B_{n-i}$.

In case $X$ is the Deligne complex of type $\tilde{C_n}$, note that $L_{1,i}$ is the poset of proper parabolic subgroups of the Artin group of type $B_i$ (with reversed order), and that $L'_{0,n-i}$ is the poset of proper parabolic subgroups of the Artin group of type $B_{n-i}$.
\ep

We can now prove the following.

\bthm \label{thm:building_typesBCD_injective}
Any Euclidean building of type $\tilde{B_{n}}$, $\tilde{C_{n}}$ or $\tilde{D_{n}}$, or the Deligne complex of type $\tilde{C_{n}}$, endowed with the piecewise $\ell^\infty$ metric $d$, is injective.
\ethm

\bp
According to Proposition~\ref{pro:local_splitting}, the space $X$ is locally isometric to a product of orthoscheme complexes of the type $|L|$, where $L=L_0 \cup \{0\}$ is a poset, and $L_0$ is either the poset of vertices of a spherical building of type $B_i$ or the poset of maximal proper parabolic subgroups of the Artin group of spherical type $B_i$, for some $0 \leq i \leq n$.

\mk

According to Propositions~\ref{pro:semilattice_building_bn} and \ref{pro:semilattice_artin_bn}, we know that in each case, $L$ is a graded meet-semilattice with minimum $0$ of rank $i$, such that any upperly bounded subset has a join. According to Theorem~\ref{thm:semilattice_injective}, we deduce that $|L|$ is injective.

\mk

We know that $X$ is uniformly locally injective. According to Theorem~\ref{thm:cartan_hadamard_injective}, we deduce that $X$ is injective.
\ep

We can also deduce that the thickening is a Helly graph.

\bcor \label{cor:building_typeBCD_Helly}
The thickening of the vertex set of any Euclidean building of type $\tilde{B_{n}}$, $\tilde{C_{n}}$ or $\tilde{D_{n}}$, or the extended Deligne complex of type $\tilde{C_{n}}$, is a Helly graph.
\ecor

\bp
Let $X$ denote either a Euclidean building of type $\tilde{B_{n}}$, $\tilde{C_{n}}$ or $\tilde{D_{n}}$, or the extended Deligne complex of type $\tilde{C_{n}}$. We will see that $X$ satisfies the assumptions of Theorem~\ref{thm:injective_orthoscheme_implies_helly}.

\mk

Since simplices of $X$ have a well-defined order, we see that $X$ is a simplicial complex with ordered simplices. For each maximal simplex $\sigma$ of $X$, the minimal and maximal vertices of $\sigma$ form a maximal edge in $X$. Therefore, $X$ has maximal edges as in Definition~\ref{def:ordered_simplices_maximal_edges}. According to Theorem~\ref{thm:injective_orthoscheme_implies_helly}, we deduce that the thickening of $X$ is a Helly graph.
\ep

One can also deduce another proof of the result of Huang and Osajda that FC type Artin groups are Helly (see~\cite[Theorem~5.8]{huang_osajda_helly}), with moreover explicit Helly and injective models.

\bthm \label{thm:FC_Artin_Helly}
Let $A=A(\Gamma)$ denote a FC type Artin group, with standard generating set $S$, and let $\leq_L$ denote the standard prefix order on $A$ with respect to $S$.
\bit
\item Consider the graph $Y$, with vertex set $A$, with an edge between $g,h \in A$ if there exists $a \in A$ and a spherical subset $T \subset S$ such that $a \leq_L g,h \leq_L a\Delta_T$, where $\Delta_T$ denotes the standard Garside element of $A(T)$. Then $Y$ is a Helly graph.
\item Consider the simplicial complex $X$, with vertex set $A$, with a $k$-simplex for each chain $g_0 <_L g_1 <_L \dots <_L g_k \leq_L g_0\Delta_T$, for some spherical subset $T \subset S$. Endow $X$ with the standard $\ell^\infty$ orthoscheme metric. Then $X$ is an injective metric space.
\eit
\ethm

\bp
Let us start by proving that $X$ is injective. Note that $X$ can be described as the union of (possibly not connected) complexes $X_T$, for $T \subset S$ spherical, by restricting the spherical subsets to be contained in $T$. So $X$ is locally isometric to the geometric realization of the poset $L=\bigcup_{T \subset S \mbox{ spherical}} [1,\Delta_T]$. This poset is graded, with minimum, with bounded rank, and is a meet-semilattice.

\mk

Moreover, the flag condition from the FC type Artin group $A$ translates into the flag condition for the poset $L$. Indeed, consider $x_1,x_2,x_3$ which are pairwise upperly bounded. Let $T_1,T_2,T_3 \subset S$ denote the supports of $x_1,x_2,x_3$ respectively. By assumption, for each $i \neq j$, the subset $T_i \cup T_j$ is spherical. As a consequence, the subset $T=T_1 \cup T_2 \cup T_3$ is a complete subset of $S$. According the FC type condition, we deduce that $T$ is spherical. Hence $x_1,x_2,x_3 \leq_L \Delta_T$. So $L$ is a flag poset.

\mk

According to Theorem~\ref{thm:semilattice_injective}, we deduce that $|L|$ is injective. So $X$ is uniformly locally injective and according to Theorem~\ref{thm:cartan_hadamard_injective}, we deduce that $X$ is injective.

\mk

We now turn to the proof that $Y$ is a Helly graph, as in the proof of Corollary~\ref{cor:building_typeBCD_Helly}. Since simplices of $X$ have a well-defined order, we see that $X$ is a simplicial complex with ordered simplices. For each maximal simplex $\sigma$ of $X$, the minimal and maximal vertices of $\sigma$ form a maximal edge in $X$. Therefore $X$ has maximal edges. According to Theorem~\ref{thm:injective_orthoscheme_implies_helly}, we deduce that the thickening $Y$ of $X$ is a Helly graph.
\ep

\section{Bicombings on Deligne complexes in types $\tilde{A_n}$ and $\tilde{C_n}$} \label{sec:corollaries_using_bicombing}

We now see that, in the Deligne complex of spherical types $A_n$, $B_n$ and Euclidean types $\tilde{A_n}$, $\tilde{C_n}$, we may find a convex bicombing.

\bthm \label{thm:deligne_bicombing}
Let $X$ denote the Deligne complex of the Artin group $A$ of spherical type $A_n$, $B_n$ or Euclidean type $\tilde{A_n}$, $\tilde{C_n}$. There is a metric $d_X$ on $X$ that admits a convex, consistent, reversible geodesic bicombing $\sigma$. Moreover, $d_X$ and $\sigma$ are invariant under $A$.
\ethm

\bp
The statement for types $A_n$ and  $\tilde{A_n}$ is Theorem~\ref{thm:AN_bicombing}.

\mk

Assume that $X$ is the Deligne complex of type $B_n$. We have seen in the proof of Proposition~\ref{pro:semilattice_artin_bn} that $X$ may be realized as the fixed point subspace of the Deligne complex $Y$ of type $A_{2n-1}$ for an involution $s$. Note that, since $s$ is order-reversing, $s$ induces an isometry of $Y$. According to Theorem~\ref{thm:AN_bicombing}, there exists a unique reversible, convex, consistent, geodesic bicombing $\sigma_Y$ on $Y$ for a metric $d_Y$. Since $\sigma_Y$ is unique, we deduce that $s$ preserves $\sigma_Y$, and that the fixed point subspace $X$ is $\sigma$-stable. Let $d_X$ denote the restriction of $d_Y$ to $X$, and let $\sigma_X$ denote the restriction of $\sigma_Y$ to $X$. We deduce that $\sigma_X$ is a convex, consistent, reversible geodesic bicombing on $(X,d_X)$. Moreover, $d_X$ and $\sigma_X$ are invariant under $A$.

\mk

Assume that $X$ is the Deligne complex of type $\tilde{C_n}$. According to~\cite[Theorem~5.2]{digne_garside_Cn} (it is also a consequence of Corollary~\ref{cor:symmetry_artin_system}), we see that $X$ may be realized as the fixed point subspace of the Deligne complex $Y$ of type $\tilde{A_{2n-1}}$ for an involution $s$. Following the same arguments as in the type $B_n$, we conclude that there exists a metric $d_X$ on $X$ and a convex, consistent, reversible geodesic bicombing $\sigma_X$ on $(X,d_X)$ that are both invariant under $A$.
\ep

We will describe several consequences we can derive from the fact that the Deligne complex has a consistent convex bicombing, that were usually known when the Deligne complex had a CAT(0) metric.

\mk

\bcor[Okonek]
The Deligne complex $X$ of Euclidean type $\tilde{A_n}$ or $\tilde{C_n}$ is contractible. In particular, the $K(\pi,1)$ conjecture holds in these cases.
\ecor

\bp
Metric spaces with a convex bicombing are contractible.
\ep

The proof of Morris-Wright (see~\cite{morris_wright_parabolic_FC}) for the intersection of parabolic subgroups in FC type Artin groups adapts directly to our situation. It relies mainly of the result by Cumplido et al. (see~\cite{cumplido_parabolic_spherical}) that in a spherical type Artin group, the intersection of parabolic subgroups is a parabolic subgroup.

\bcor \label{cor:intersection_of_parabolics}
Let $A$ denote the Artin group of Euclidean type $\tilde{A_n}$ or $\tilde{C_n}$. The intersection of any family of parabolic subgroups of $A$ is a parabolic subgroup.
\ecor

\bp
Note that any proper parabolic subgroup of an Artin group of Euclidean type has spherical type. Since the proof of~\cite[Theorem~3.1]{morris_wright_parabolic_FC} uses only the existence of an $A$-equivariant consistent geodesic bicombing on $X$, it adapts to these cases.
\ep

The results by Godelle (see~\cite{godelle_cat0}) describing centralizers and normalizers of parabolic subgroups also adapt to our case.

\bcor \label{cor:centralizers_normalizers}
Let $A$ denote the Artin group of Euclidean type $\tilde{A_n}$ or $\tilde{C_n}$. Then $A$ satisfies Properties $(\star)$, $(\star \star)$ and $(\star \star \star)$ from~\cite{godelle_cat0}, notably for any subset $X \subset S$, we have
$$Com_A(A_X) = N_A(A_X) =  A_X \cdot QZ_A(X),$$
where the quasi-centralizer of $X$ is $QZ_A(X) = \{g \in A \st g \cdot X = X\}$.
\ecor

Note that Godelle's Property $(\star \star)$ stating that a parabolic subgroup $P$ is contained in a parabolic subgroup $Q$, then $P$ is a parabolic subgroup of $Q$ has been proved by Blufstein and Paris (see~\cite{blufstein_paris}).

\bp
The proof of~\cite[Theorem~3.1]{godelle_cat0} uses only the following properties of the Deligne complex, which hold for the Deligne complex $X$:
\bit
\item There exists an $A$-equivariant consistent geodesic bicombing $\sigma_X$ on $X$.
\item Closed cells of $X$ are $\sigma_X$-stable.
\eit
These assumptions are satisfied.
\ep

Note that Cumplido uses these properties of parabolic subgroups to solve the conjugacy stability problem for parabolic subgroups, see~\cite[Theorem~14]{cumplido_conjugacy_stability}.

\mk

The results by Crisp about symmetrical subgroups of Artin groups also extend to our case (see~\cite{crisp_symmetrical}).

\bcor \label{cor:symmetry_artin_system}
Let $A$ denote the Artin group of Euclidean type $\tilde{A_n}$ or $\tilde{C_n}$. For any group $G$ of symmetries of the Artin system, the fixed point subgroup $A^G$ is isomorphic to an Artin group.
\ecor

For the explicit description of the Artin group $A^G$, we refer the reader to~\cite{crisp_symmetrical}.

\bp
The proof of~\cite[Theorem~23]{crisp_symmetrical} uses only the following properties of the Deligne complex $X$:
\bit
\item There exists an $A$-equivariant consistent geodesic bicombing $\sigma_X$ on $X$.
\item Subsets of $X$ which are $\sigma_X$-stable are contractible.
\eit
Let us check the last condition: if $C \subset X$ is a $\sigma_X$-stable subset, then $\sigma_X$ restricts to a convex bicombing on $C$, which implies that $C$ is contractible.
\ep

One particular case of interest is when $A$ is of Euclidean type $\tilde{A_{2n-1}}$, for some $n \geq 1$, and the group $G$ of symmetries of the Artin system (the $2n$-cycle) is generated by a symmetry of the cycle defining $A$. We recover the result from Digne (see~\cite{digne_garside_Cn}) that the fixed subgroup is isomorphic to the Artin group of type $\tilde{C_n}$.

\mk

For any divisor $1 < k < n$ of $n$, one may also consider the group $G$ of symmetries of the Artin system of type $\tilde{A_{n-1}}$ (the $n$-cycle) generated by a rotation of $k$. The fixed point subgroup is then isomorphic to the Artin group of Euclidean type $\tilde{A_{k-1}}$. When we see $A(\tilde{A_{n-1}})$ as the group of braids with $n$ strands on the annulus, we may think of the fixed subgroup as the subgroup of braids which are invariant by a rotation of the annulus of angle $\f{2\pi k}{n}$. 

\bibliographystyle{alpha}
\bibliography{../../../bibli}

\end{document}